\documentclass[11pt]{article}
\usepackage [dvips]{graphics}
\usepackage[centertags]{amsmath}
\usepackage{amsfonts}
\usepackage{amssymb}
\usepackage{amsthm}
\usepackage{newlfont}
\usepackage{xcolor}
\usepackage{tikz}

\newlength{\defbaselineskip}
\newcommand{\setlinespacing}[1]%
           {\setlength{\baselineskip}{#1 \defbaselineskip}}

\theoremstyle{plain}
\newtheorem{thm}{Theorem}[section]
\newtheorem{cor}[thm]{Corollary}
\newtheorem{lem}[thm]{Lemma}
\newtheorem{prop}[thm]{Proposition}
\newtheorem{exam}[thm]{Example}
\newtheorem{rem}[thm]{Remark}
\newtheorem{Def}[thm]{Definition}

\newtheorem{conjecture}[thm]{Conjecture}

\makeatletter\@addtoreset{equation}{section} \makeatother
\begin{document}
\title{A new perspective on  Riemann surfaces and
von Neumann algebras via
thin Blaschke products}



%
 \author{ Danni Guo, Hansong Huang, Shan Li, Xianfeng Zhao} 
\date{}

\maketitle \noindent\textbf{Abstract:}
   In 2014,
 K.~Guo and the second author
 constructed the first example of an infinite Blaschke product
 $B$  (a thin Blaschke product), for which the
 associated multiplication  operator $M_B$ is irreducible
   on  the Bergman space $L_a^2(\mathbb{D})$ over the unit disk $\mathbb D$.
  We take a substantial step forward toward  characterizing
    the irreducibility of the  operator $M_B$ on
    $L_a^2(\mathbb{D})$,   or equivalently, the trivialness of the von Neumann algebra
    $\{M_B,M_B^*\}'$. This characterization
    indicates  a close connection with the primeness of $B$ and with
   the two-component property of the Riemann surface $\mathcal{S}_B$ associated with $B$. The approach is built upon
    a nested family of
    connected components of the sublevel set
 $\{z\in \mathbb{D}:|B(z)|<t\}$ for $0<t<1,$
 which  captures the function-theoretic and geometric feature of
  $B$.







\section{Introduction}

%
%
%

 Let $\mathbb{D}$ denote the unit disk, and let $H^\infty(\mathbb D)$
  be the algebra of bounded analytic functions on $\mathbb{D}$.
    Let \(\mathcal H\) be the Hardy space \(H^2(\mathbb D)\) or the Bergman space \(L_a^2(\mathbb D)\), the latter consisting of all holomorphic functions \(f(z)=\sum\limits_{n=0}^\infty a_nz^n\) on \(\mathbb D\) satisfying
\[
\|f\|^2=\int_{\mathbb D}|f(z)|^2\,dA(z)<\infty,
\]
where \(dA\) denotes the area measure on \(\mathbb D\).
Each
   function $\varphi $ in $H^\infty(\mathbb{D})$ induces a bounded \emph{multiplication operator}
\[
M_\varphi: \mathcal{H} \to \mathcal{H}, \quad f \mapsto \varphi f.
\]
A closed subspace $\mathcal{M} \subseteq \mathcal{H}$ is \emph{reducing} for $M_\varphi$ if it is invariant under both $M_\varphi$ and its adjoint $M_\varphi^*$. The operator $M_\varphi$ is called \emph{reducible} if it admits a nontrivial reducing subspace (i.e., one other than $\{0\}$ and $\mathcal{H}$); otherwise it is \emph{irreducible}. Equivalently, $M_\varphi$ is irreducible if and only if its commutant $\{M_\varphi\}'$ consists only of scalar multiples of the identity; that is, the von Neumann algebra
\[
\mathcal{V}^*(\varphi) = \{M_\varphi, M_\varphi^*\}'
\]
satisfies $\mathcal{V}^*(\varphi)=\mathbb{C}I$.
It is known that $ \mathcal{V}^*(\varphi)$ is generated by all orthogonal projections
from $\mathcal{H}$ onto reducing subspaces of $M_\varphi.$

\subsection{Background}
A central theme at the interface of function theory and operator theory is the study of commutants of multiplication operators $M_\varphi$ ($\varphi \in H^\infty(\mathbb D)$), their reducing subspaces, and the von Neumann algebras they generate. A fundamental problem is: \emph{When does $M_\varphi$ have no nontrivial reducing subspace? That is, when is $M_\varphi$ irreducible?}
By the inner--outer factorization of $H^2(\mathbb D)$-functions, it is of particular interest to consider inner functions  $\varphi$;  By Frostman's
 theorem~\cite{Gar},
 this further reduces to
  Blaschke products. Though formulated
  in operator-theoretic language, this question
   lies at the heart of two seemingly distinct lines of inquiry.  

 First, the study of invariant subspaces,
and in particular reducing subspaces, of operators
 on function spaces has driven major developments in analysis.
  Beurling's theorem~\cite{Be48} fully characterizes $M_z$-invariant
  subspaces on $H^2(\mathbb D)$. On the Bergman space
  $L_a^2(\mathbb D)$,
 the full invariant-subspace structure remains
  open due to Aleman--Richter--Sundberg~\cite{ARS96}.
  Moreover, the Invariant Subspace Problem for Hilbert-space
  operators is equivalent to the saturation
  of the lattice $\operatorname{Lat} (M_z)$ of invariant subspaces of
  $M_z$ on $L_a^2(\mathbb D)$.
Attention has therefore turned to the study of reducing subspaces
 for multiplication operators  $M_\varphi$
  with $\varphi\in H^\infty(\mathbb D)$,
  which form a distinguished subclass of invariant subspaces.
   Considerable progress has been made in this direction,
   particularly for finite Blaschke  products $B$,
    as well as their counterparts on Hardy and Bergman spaces
    over domains in $\mathbb C^d$
    (see~\cite{BDGS,BSMS,DSZ,DPW,Gh22,GSZZ09,HZ25,LYZ23,SZZ}). 
       The study of reducing subspaces of $M_B$
       on $L_a^2(\mathbb D)$ was initiated
       in \cite{SW,Zhu2}, where $M_B$ was shown to have finitely
       many minimal reducing subspaces when $B$ is of order $2$.
       Subsequent work in \cite{DSZ,GSZZ09,SZZ} revealed a
       key connection: the dimension of $\mathcal{V}^*(B)$
       equals the number of connected components of the Riemann
       surface $\mathcal{S}_B$ associated with $B$. Douglas,
       Putinar and Wang \cite{DPW} later proved this number also
       equals the number of minimal reducing subspaces of $M_B$ by
        establishing commutativity of $\{M_B,M_B^*\}'$ for all finite
        Blaschke products $B$. Very recently, the first three authors
        and Luo \cite{GHLiLuo} introduced the notion of terminal
        branched components of $B$, via critical point analysis,
        to classify reducing subspaces of $M_B$.

Second, a parallel line of inquiry, rooted
in the group factor problem,
further motivates the study of reducing subspaces and the
von Neumann algebras $\mathcal{V}^*(\phi)$ associated with
multiplication operators $M_\phi$ induced by Blaschke products 
$\phi$. Recall that for $n \geq 1$, $\mathbb{F}_n$ denotes
the free group on $n$ generators; the Group Factor Problem
asks whether $\mathcal{L}(\mathbb{F}_n)$ is
 $C^*$-isomorphic to
$ \mathcal{L}(\mathbb{F}_m)$ for distinct $n,m \geq2$,
where $\mathcal{L}(\mathbb{F}_n)$ is the von Neumann factor
generated by the left regular representation of $\mathbb{F}_n$
on $\ell^2(\mathbb{F}_n)$.
For $k \geq2$, define $\Omega_k = \mathbb{D} \backslash
 \{a_1, a_2, \cdots, a_{k}\}$, where
 $a_1, a_2, \cdots, a_{k} \in \mathbb{D}  \backslash\{0\}$
  are distinct points. Then there exists a holomorphic covering map
  $\phi_k : \mathbb{D} \to \Omega_k$ that is an interpolating
   Blaschke product
  (see~\cite[Example 3.6]{GH11}). By \cite[Theorem 1.3]{GH11},
  $$\mathcal{V}^*(\phi_k) \simeq \mathcal{L}(\pi_1(\Omega_k))
   \simeq \mathcal{L}(\mathbb{F}_k).$$ Hence, the Group Factor Problem is equivalent to asking whether $\mathcal{V}^*(\phi_n) \simeq \mathcal{V}^*(\phi_m)$ for distinct $n,m \geq2$; see~\cite{GH11} for further discussion.

   These two lines of inquiry converge naturally on the following problem:

 \textbf{Problem.} \emph{For a given subclass $\mathcal{C}$ of Blaschke products,
can we characterize the irreducibility of the multiplication operator
  $M_B$ with $B$ in $\mathcal{C}$ in terms of the function-theoretic  or
 geometric properties of  $B$?}

In this paper we focus on a distinguished
   subclass, namely thin Blaschke products,
   and investigate how the irreducibility of
    $M_B$ on $L_a^2(\mathbb D)$ relates to the
    primeness of $B$ and the structure of the associated
    Riemann surface $\mathcal{S}_B$.

\subsection{Our main result}

We begin with some preliminary notations.

   For a Blaschke product  $B$,
 \emph{a critical point} of $B$ is defined as
  a point in the unit disk $\mathbb{D}$
 where its  derivative $B'$  vanishes. Define
 $$Z(B')=\{z\in \mathbb{D} :  B'(z)=0\}, $$
 and   its \emph{critical value set}  $ \mathcal{E}_B$
 by $$ \mathcal{E}_B  = \{B(z): z\in Z(B')\}. $$
  Recall that a thin Blaschke product is a Blaschke product
whose zero sequence $\{z_n\}_{n=1}^\infty$ satisfies
   $$\lim_{n\to \infty}\prod_{\substack{k\geq 1\\ k \neq n}}
   \Big|\frac{z_k-z_n}{1-\overline{z}_kz_n}\Big|=1.$$
   For a thin/Blaschke product $B$,
  $ \mathcal{E}_B$ is a discrete subset
  of $\mathbb{D}$ (see \cite[Proposition 2.1]{GH3} for example). We define
   $$\mathcal{S}_B: =\{(z,w)\in \mathbb{D}^2:
  B(z)=B(w)\not\in \mathcal{E}_B \}, $$
  which is a Riemann surface associated with $B $ \cite{GH3}.
   It is known that every thin Blaschke product $B$ has an infinite
   critical
 set $Z(B') $ (see Proposition \ref{infcritical}).
  We say $B$ is \emph{essentially injective} on its critical set
 if there exists a finite subset $F$ of $Z(B') $ such that $B$
 is injective on $Z(B')\backslash F $.

   K.~Guo and the second author \cite{GH3} gave a criteria for irreducibility of $M_B$ when $B$ is a thin Blaschke product .
 \begin{thm}    \cite[Theorem 4.1]{GH3} \label{plms}
 Suppose $B$ is a thin Blaschke product
 satisfying  the following:
  \begin{itemize}
\item[\emph{(1)}] $B$ is essentially injective on its critical set;
\item[\emph{(2)}] There is a finite subset $F$ of $Z(B')$
 such that each point in $Z(B')\backslash F$ is a simple zero of $B'$.
\end{itemize}
Then $\mathcal{V}^*(B)$  is trivial;  equivalently,
the multiplication operator $M_B$ on $L_a^2(\mathbb{D})$
 is irreducible.
  \end{thm}

  Following \cite{Step},
 an inner function $f$ is called \emph{prime}
if for any given composition $f=u\circ v$ with  inner functions $u$ and $v$, at least
 one of them is a M\"obius transform.  Motivated by Theorem \ref{plms},
we propose the following conjecture.

 \begin{conjecture} \label{Conj1}
Let $B$ be a thin Blaschke product.
Then the following are equivalent.
\begin{itemize}
\item[\emph{(1)}] $B$ is prime;
\item[\emph{(2)}] The von Neumann algebra $\mathcal{V}^*(B)$
 is trivial; equivalently, 
$M_B$ on $L_a^2(\mathbb{D})$ is irreducible;
\item[\emph{(3)}] The Riemann surface $\mathcal{S}_B$ has exactly
two components.
\end{itemize}
\end{conjecture} 

\begin{rem} We have some comments for Conjecture \ref{Conj1}.
\begin{itemize}
\item[\emph{(1)}] The essential injectivity of a
thin Blaschke product implies its primeness
(see Proposition \ref{onlyif}).
 Theorem \ref{plms} is therefore a special case of Conjecture \ref{Conj1}.
\item[\emph{(2)}]  Condition (2) is equivalent to the statement
 that every bounded operator on $L_a^2(\mathbb{D})$
  can be approximated in the strong operator topology by
  non-commutative polynomials in
$M_B$ and $M_B^*$, by the von Neumann Bicommutant Theorem. 
\item[\emph{(3)}] In Section \ref{sect4}, we also construct an example of a product $B$ of finite many
  thin Blaschke products,
  which is not prime, but $\mathcal{V}^*(B)$ is  trivial.
  This indicates that Conjecture \ref{Conj1} is sharp.
\item[\emph{(4)}]  For a thin Blaschke product $B$,
 \cite[Theorem 5.7.2]{GH4} established   that $M_B$ is irreducible
 if and only if $\mathcal{S}_B$ has no non-identity component
 of finite multiplicity; here the identity component is
  $\{(z,z)\in \mathbb{D}^2:B(z)\not\in \mathcal{E}_B\}$.
  This was first proved in \cite{GH3} under an additional mild
  assumption.

  Conjecture~\ref{Conj1} refines the existing criterion in
   two respects: it strengthens the condition on
    $\mathcal{S}_B$ to exactly two components, and it introduces
  a purely function-theoretic condition, primeness.
  \end{itemize}
\end{rem}

We propose Conjecture \ref{Conj1}  partially based on
 the implication chain
  $(3) \Rightarrow (2) \Rightarrow (1)$
  (see the proof of Proposition \ref{eqq}).
 It remains open whether the reverse implications
 $(1) \Rightarrow (3)$ or
  $(1) \Rightarrow (2)$   holds.
  The following conjecture, which asserts
   $(1) \Rightarrow (3)$, is therefore
    equivalent to Conjecture~\ref{Conj1}.
  \begin{conjecture} \label{Conjt}
Let $B$ be a thin Blaschke product. If $B$ is prime, then
  the Riemann surface $\mathcal{S}_B$ has exactly two
  components.
 \end{conjecture}
  Our main result  provides substantial support
   for   Conjecture \ref{Conjt}.
   Moreover, it essentially  generalized Theorem \ref{plms}.
 \begin{thm}  \label{main}
   Suppose $B$ is a thin Blaschke product
 such that $B$ is essentially injective on its critical set.
 Then the Riemann surface
    $ \mathcal{S}_B$ has exactly  two components.
Consequently, $\mathcal{V}^*(B)$   is trivial; equivalently,
the multiplication operator $M_B$ on $L_a^2(\mathbb{D})$
 is irreducible.
  \end{thm}

 For thin Blaschke products, primeness and essential injectivity are closely related. In fact, if
  $B$ is a thin Blaschke product and $B$
  is essentially injective on
its critical points, then $B$ is prime (see Proposition \ref{onlyif}).
It seems that the two properties are not far apart.
 One reason is that essential injectivity
 carries a clear geometric intuition.
In fact, we will show later  (see Proposition \ref{EIII}) that for any thin Blaschke product
  $B$, the following are equivalent:

 (I) \emph{The mapping B is not essentially
  injective on its critical set.}\\

(II) \emph{There exist  two sequences $\{z_n\}_{n=1}^\infty$ and $\{w_n\}_{n=1}^\infty $
in $\mathbb{D}$ such that}
\begin{itemize}
\item[(1)] $B(z_n) = B(w_n)  $  and $z_n \neq w_n$ for all $n$;
\item[(2)] $B'(z_n) = B'(w_n) = 0;$
\item[(3)] $B(z_n)\neq B(z_m)$ if $m\neq n$.
\end{itemize}
 By  \cite[Theorem 4.1]{GLMR94}, if  $B$ is a thin
  Blaschke product with the decomposition
$B=u\circ v$,  $u$ and $v$ being inner functions,
then either $u$ is a M\"obius map or $v$ is
 a finite Blaschke product.
 Using this result,
   we can show that a non-prime thin Blaschke product $B$ satisfies Condition (II); that is,
   $B$ is not essentially injective on its critical set.
    This observation, along with Theorem \ref{main}, lends
    some support to Conjecture \ref{Conjt}.

It is tempting to conjecture that
 primeness and essential injectivity are
 equivalent for thin Blaschke products.
 If this were true, Conjectures~\ref{Conj1}
 and~\ref{Conjt} would follow immediately
 from Theorem~\ref{main}.
 However, in  Section~\ref{6cexample},  a prime thin Blaschke product that fails to be essentially injective on its critical
 set is constructed.
 This example does
  not contradict Conjectures~\ref{Conj1} or~\ref{Conjt}. 



\subsection{The strategy for our main theorem}

 Our main idea is as follows.
 A thin Blaschke product $B$ can be viewed
 as the union of countably many finite Blaschke products.
 In fact, we can choose   increasing $t_n \in (0,1)$ and
  simply connected components $\Omega_n$ of $\{ |B| < t_n \}$
   with  $\Omega_n \subsetneqq \Omega_{n+1}$ and
   $\bigcup\limits_{n=1}^\infty \Omega_n = \mathbb{D}$.
   Each restriction $B|_{\Omega_n}: \Omega_n \to \mathbb{D}$ is
    a proper holomorphic map, factoring as a finite Blaschke
     product conjugated by conformal maps.
     This reduces the analysis of $B$ to that of
     the sequence $\{B|_{\Omega_n}\}$.
 In this way,  the proof for the following theorem
 (the case of finite Blasche product) serves as the core ingredient of  Theorem \ref{main}.

 \begin{thm} \label{submain}
 Suppose $B$ is a finite Blaschke product with
  at least two distinct critical points and
   $B|_{Z(B')}$ is injective.
 Then  $\mathcal{S}_B$ has exactly two components. In this case,
  the multiplication operator $M_B$ on $L_a^2(\mathbb{D})$
  has exactly two minimal reducing subspaces.
   \end{thm}


For a thin Blaschke product $B$,
Theorem~\ref{main} gives a function-theoretic criterion for
 the irreducibility of $M_B$. Our perspective is
 to view $B$ as the limit of an exhaustion by finite Blaschke
 products. The proof of Theorem~\ref{main} rests on two key ideas.

First, the proof of Theorem~\ref{main}
proceeds by analyzing the critical points of $B$.
 Unlike the earlier result in Theorem \ref{plms}, Theorem \ref{main}
 allows critical points of higher order. For finite/thin
 Blaschke products, the Riemann surface $\mathcal{S}_B$
  can be studied via local inverses, where points in a
  fiber $B^{-1}(w)$ are connected by analytic continuation.
  In the presence of higher-order critical points, however,
  it becomes difficult to determine which local inverses
  lie in the same equivalence class.
  A global argument 
   is therefore needed, and finally we prove that $\mathcal{S}_B$ has exactly
   two components under the stated conditions. This global
   analysis is the the key new contribution of  Theorem~\ref{main},
    compared with earlier results in~\cite{GH3}.

Second, our approach differs substantially from that of
\cite{GHLiLuo},
where only a single value $t_0\in(0,1)$ and the components of $\{z\in \mathbb{D}:|B(z)|<t_0\}$
 (the terminal branched components) were needed.
  Here, by contrast, we must treat $\{z\in \mathbb{D}:|B(z)|<t\}$
   for every $t\in(0,1]$. To this end,
    we partition $(0,1]$ into finitely many
    subintervals and carry out an induction
    on $t$ over these subintervals.
    We then prove that the sub-Riemann surfaces
 of $\mathcal{S}_B$ lying on $\{z\in \mathbb{D}:|B(z)|<t\}$
 admit exactly two components under suitable conditions,
 and finally show that
     $\mathcal{S}_B$ itself has precisely two components.

The paper is organized as follows.
 Section~\ref{sect2} presents preliminaries on
 finite and thin Blaschke products and establishes
  several lemmas. 
  Section~\ref{sect3} gives the proofs of
  Theorems~\ref{main} and \ref{submain}. Section~\ref{sect4} shows that the
   essential injectivity of a thin Blaschke
    product implies its primeness, and also
    establishes the implications
  $(3) \Rightarrow (2) \Rightarrow (1)$ in
     Conjecture~\ref{Conj1}.
     Section~\ref{sect5}
  constructs a product \(B\) of two thin Blaschke
  products, showing that Conditions~(1)--(3)
  in Conjecture~2 are not equivalent for this
  \(B\).
     Finally,
      Section~\ref{6cexample} constructs a
       prime thin Blaschke product that is not
 essentially injective on its critical set.

Throughout this paper,
$B$ denotes a finite/thin Blaschke product.
For $0<t<1,$ the set $\{|B|<t\}$ denotes
$\{z\in \mathbb{D}:|B(z)|<t\}$, and $\{|B|\leq  t\}$ denotes
$\{z\in \mathbb{D}:|B(z)|\leq  t\}$.

\section{Some preliminaries} \label{sect2}
\subsection{Finite and thin Blaschke products}
In this subsection, we collect some
facts about finite and thin Blaschke products.

By direct computations, the derivative of
a finite Blaschke product  $B$   vanishes nowhere
 on the unit circle. Remind  that  a critical point of a Blaschke product $B$ refers to a point in $\mathbb{D}$
 where $B'$  vanishes. We   have the following statement  \cite{GMR18}.
\begin{lem}\label{critical}
A Blaschke product of order $n$ contains exactly $n-1$
critical points counting multiplicity
  in the unit disk $\mathbb{D}$.
\end{lem}

 The following lemma is from \cite[Lemma 2.2]{GH3}.
 Given two domains $\Omega_1, \Omega_2$,
 a holomorphic map $f: \Omega_1 \to \Omega_2$
 is called \emph{proper} if the preimage
 of every compact subset of $\Omega_2$ is
  compact in $\Omega_1$. Equivalently, for any sequence $\{z_n\}_{n=1}^\infty$ satisfying $z_n \to \partial \Omega_1$, we have $f(z_n) \to \partial \Omega_2$.
  \begin{lem}\label{proper}
Suppose that $\phi$ is a holomorphic function over
 $\mathbb{D}$ and $\Delta$ is a planar domain.
 If $\Delta_0$ is a component of $\phi^{-1}(\Delta)$
such that $\overline{\Delta}_0\subseteq \mathbb{D}$,
then the restriction $\phi|_{\Delta_0}$ is a proper map onto $\Delta$.
\end{lem}

 Lemma \ref{proper} directly implies Proposition \ref{sproper} below;
   see also \cite[Proposition 2.7]{GHLiLuo}.
   \begin{prop} \label{sproper}
Let $B$ be a finite Blaschke product, and
let $\Omega$ be a component of $\{|B|<t\}$ for $t\in (0,1)$.
Then $\Omega$ is a simply connected domain.
In this case, there exists  a conformal map $G$ from $\Omega $ onto
   the unit disk, and a finite Blaschke product $\varphi$
   such that  $$B|_{\Omega}=  t \varphi\circ G .  $$
\end{prop}
  We include a proof of the following  fact
for completeness; see also \cite[Remark 3.9(e)]{GHLiLuo}.

 \begin{lem}\label{entirecompo}
Let $B$ be a finite Blaschke product.
   For each $s$ in $(0,1)$, write  $$\Omega_s:= \{z\in \mathbb{D}:|B(z)|<s\}.$$
 There exists a number $t$ in $(0,1)$ such that $\Omega_s$
  is connected for all $s\geq t.$
\end{lem}
\begin{proof} Let $r$ be the maximum module of
 zeros of $B$, and let $t$ be the maximum of $B$ on the closed disk
 $\{z\in \mathbb D: |z|\leq   \frac{r+1}{2}\}$.
 For
 $s\geq  t $, let $ \Omega_1$  be the component of $\Omega_s$
  containing all zeros of $B.$
 We will prove that  $\Omega_s=\Omega_1 $ to
 complete the argument.
Suppose for contradiction that
 there exists another distinct
  component  $\Omega_2 \subseteq \Omega_s$.
 Then by Lemma \ref{proper}
  $B|_{\Omega_2}:\Omega_2\to s\mathbb{D}$ is a proper map,
   and onto. In particular,
  $B$ has at least one zero in $\Omega_2.$
  This leads to a contradiction with
  the fact that $\Omega_1$ contains all zeros of $B$.
\end{proof}

 The next result is included in \cite[Theorem 2.6]{GHLiLuo}.
\begin{prop} \label{Jordancurve}
Let $\phi$ be a holomorphic map from $\mathbb{D}$ to
 $\mathbb{D}$, and fix $0<t<1$. If $\Omega$ is a component
 of $\{|\phi|<t\}$    and $\overline{\Omega}\subseteq \mathbb{D}$,
 then $\Omega$ is bounded by a Jordan curve.
\end{prop}
\begin{proof} The idea of the proof comes from
that in \cite[Theorem 2.6]{GHLiLuo}. Since $\phi$ is holomorphic on $\mathbb{D}$,
  by the Maximal Modulus Principle,
   it is routine to show that
   each component $\Omega$ of the set $\{|\phi|<t\}$
  is simply connected.
It remains to prove that $\partial\Omega$ is a Jordan curve.

   Assume for contradiction that $\partial\Omega$ is not a Jordan curve. Then $\partial\Omega$ contains a Jordan curve $\gamma$; let $\Omega_\gamma$ denote its interior. Choose $\Gamma$ among all such curves so that $\Omega_\Gamma$ is maximal with respect to inclusion.
   Since $|\phi| = t$ on $\Gamma$ and $|\phi| < t$
   on $\Omega_\Gamma$, and
    since $\Omega  \subseteq \{|\phi|<t\}$,
    it follows that $\Omega \subseteq \Omega_\Gamma$.
   However, as $\Omega$ is a component of $\{|\phi|<t\}$,
     we must have $\Omega = \Omega_\Gamma$. This implies that
      $\partial\Omega = \Gamma$,
      which contradicts the assumption that
      $\partial\Omega$ is not a Jordan curve.
Therefore, $\partial\Omega$ must be a Jordan curve, as desired.
\end{proof}

The following thin Blaschke product version of Proposition
 \ref{sproper} appears in \cite[Proposition 2.3]{GH3}; we include a proof for completeness.
 \begin{prop} \label{atomz}
If $B$ is a thin Blaschke product, then for each
$t \in (0,1)$, every   component $\Omega$ of
$\{ |B| < t \}$ is compactly contained in $\mathbb{D}$,
and $B: \Omega \rightarrow t\mathbb{D}$ is a proper map
between simply connected domains.
\end{prop}

 Define  the pseudo-hyperbolic distance
$$ d(z,w) =  \left| \frac{z-w}{1-\overline{z}w} \right|. $$
To prove Proposition \ref{atomz},  we need the following lemma
(see \cite{GM05,Lawson69}, and \cite[Lemma 4.2]{Gar}).
\begin{lem} \label{atomic}
For a thin Blaschke product $B$ with zero sequence
$\{z_n\}$ and every $\delta \in (0,1)$, there exists
 an $\varepsilon \in (0,1)$ such that $|B(z)| \geq\delta$
 whenever $d(z,z_n) \geq \varepsilon$ for all $n$.
\end{lem}
\noindent \textbf{Proof of Proposition \ref{atomz}. }
The idea of the proof is from that of
\cite[Proposition 2.3]{GH3}.

For each
$t \in (0,1)$, let  $\Omega$ be a connected component of
$\{ |B| < t \}$.  It is straightforward to
  verify that $\Omega$ is simply connected.
  Our goal is to prove that
  $B: \Omega \to t\mathbb{D}$ is a proper map.

Given that $B$ is thin with the zero sequence $\{z_n\}_{n=1}^\infty$,
by Lemma \ref{atomic},
  there exists $\varepsilon \in (0,1)$ satisfying
 $$ B^{-1}(\Delta(0,t)) \subseteq \bigcup_{n=1}^\infty
 \Delta(z_n, \varepsilon). $$
  For intersecting disks
  $\Delta(z_n,\varepsilon) \cap \Delta(z_k,\varepsilon)$
  containing a point $z'$, we apply [16, Lemma 1.4]
  to obtain  \begin{equation} \label{omm}
  d(z_n,z_k) \leq    \frac{d(z_n,z') +
   d(z',z_k)}{1 + d(z_n,z')d(z',z_k)} \leq     d_\varepsilon,
   \end{equation}
    where $$d_\varepsilon = \max  \Big\{ \frac{x+y}{1+xy}:
   0\leq     x,y \leq    \varepsilon \Big\}  < 1.$$

By the thinness condition of $B$, we have
 $$\lim_{n\to\infty} \prod_{\substack{k\geq 1\\ k \neq n}} d(z_n,z_k) = 1.$$
There exists  $N \geq 1$ such that for all $n \geq  N$,
$$d(z_n,z_k)\geq \dfrac{1+d_\varepsilon}{2}>d_\varepsilon  .$$
Consequently, by
(\ref{omm}),
$\Delta(z_n,\varepsilon)$ has no intersection with other
disks $\Delta(z_m,\varepsilon)$ with $m\neq n$.
The union $K = \bigcup\limits_{n=0}^N \overline{\Delta(z_n,\varepsilon)}$
is compact. For every component $\Omega$ of ${ |B| < \delta }$,
either  $\Omega \subseteq K$,
or  $\Omega \subseteq \Delta(z_n,\varepsilon)$ for some
$n \geq  N$.
In both cases, $\overline{\Omega}$ is compact.
Furthermore, standard analysis shows
 $B(\partial \Omega) \subseteq \partial(t\mathbb{D})$.
 Therefore, by continuity of $B$ on $\overline{\Omega}$,
  $B$ is a proper map as desired. $\hfill \square$

The following result, while likely known in the literature,
is included here for completeness.
\begin{prop}\label{infcritical}
 Every thin Blaschke product has
 infinitely many critical points in $\mathbb{D}$. \end{prop}

\begin{proof}Let $B$ be a  thin Blaschke product
with the zero sequence $ z_1,z_2,\cdots $.  For every $n \geq 2$,
 consider a disk $D$ containing the first $n$
 zeros ${z_1,\cdots,z_n}$.
 There exists $t \in (0,1)$ such that $|B| < t$ on $\overline D$. Let $\Omega$ be the component of $\{|B|<t\}$ containing $\overline D$. Then the restriction $B|_\Omega: \Omega \to \mathbb D$ is a finite Blaschke product of degree at least $n$.
   More precisely, let $F$ be a conformal map from the unit
   disk $\mathbb{D}$ onto $\Omega.$
 Defining $h=\frac{1}{t}B|_\Omega \circ F  $, we observe that
$h$ is a proper map from $\mathbb{D}$ onto $\mathbb{D}$.
Then by \cite[7.3.1]{Ru1}, 
$h$ must be a finite Blaschke product, and $B|_\Omega=th\circ F^{-1}$.
  Therefore, $B|_\Omega$ inherits at least $  n-1$
   critical points from $h$, counting multiplicity.   Since $n$ is
    arbitrary, we conclude that $B$ possesses infinitely
    many critical points in $\mathbb{D}$.
\end{proof}

Thin Blaschke products have many good properties.
The following characterization
is due to \cite[Theorem 2.3]{Tol91}.
\begin{thm} \label{thin} A Blaschke product $B$ is thin if and only if
$$\lim_{|z|\to 1^-} \big[|B(z)|^2+(1-|z|^2)|B'(z)|\big]=1.$$
\end{thm}

For $w\in \mathbb{D}$, define $\varphi_w(z)=\dfrac{w-z}{1-\overline{w}z},~$ $z\in \mathbb{D}.$
As a consequence of the above theorem, we have the following
(see \cite[Proposition 4.3]{GLMR94} and \cite[Proposition 4]{CGP2013}).
\begin{prop} \label{thin2} Let $B$ be a thin Blaschke product, $w\in \mathbb{D}$
and $\phi$ be a finite Blaschke product. Then both
$\varphi_w\circ B$ and $B\circ \phi$ are also thin Blaschke products.
\end{prop}
\subsection{Local inverse and gluable property}
 In this subsection, we present the definitions of local inverse and
analytic continuation. Also, we recall the gluable property of
 finite and thin Blaschke products from \cite{GH3}.

To begin with, we recall the notions of analytic continuation and local inverse (see \cite[Chapter 16]{Ru3} and \cite{T1}, respectively).
 A \emph{function element}   is an ordered pair $(f,D)$, where $D$
is an open disk   and $f$ is a holomorphic function on
$D$.
A \emph{curve} is   defined
 as a continuous map  $\gamma: [0,1] \to \mathbb{C} $.
 Consider a function element $(f_0,D_0)$
 and a curve $\gamma$ with $\gamma(0) \in D_0$. If there exists
  a partition $$0 = s_0 < s_1 < \cdots < s_n = 1$$ of $[0,1]$,
  and  a collection
    of function elements $\{(f_i,D_i)\}_{0 \leq    i \leq    n}$  so that
         $(f_j,D_j)$ and $(f_{j+1},D_{j+1})$ are direct continuations,
          and $\gamma([s_j,s_{j+1}]) \subseteq D_j$
           for $0 \leq   j \leq   n-1$ with $\gamma(1) \in D_n,$
  then $(f_n,D_n)$ is called
 \emph{an analytic continuation of $(f_0,D_0)$ along $\gamma$},
  and we say $(f_0,D_0)$ \emph{admits analytic continuation
  along $\gamma$}. This relationship
  is denoted by $f_0 \sim f_n$,
  which establishes an equivalence relation.
For a finite/thin Blaschke product $B$,
and a subdomain $\Delta $ of $\mathbb{D} \backslash B^{-1}(\mathcal{E}_B)$,
a holomorphic map $\rho$ on $\Delta $ is called
 \textbf{a local inverse}  of $B$ if $B\circ \rho =B.$
 For  two local inverses $\rho_1 $ and $\rho_2$ of $B$,
 we say $\rho_1 \sim \rho_2 $  if
   $\rho_1$ is an analytic continuation
 of $\rho_2$ along a curve in
 $\mathbb{D} \backslash B^{-1}(\mathcal{E}_B).$

 For every  $a,b\in \mathbb{D} \backslash B^{-1}(\mathcal{E}_B),$
 there are small  neighborhoods $V_a $ of $a$, $V_b$ of $b$,
  and a disk $D$ such that
   both $  B:V_a\to D  $ and $B:V_b\to D$ are
 biholomorphic maps. Then define
 $$\rho_{a,b} (z)=(B|_{V_b})^{-1}\circ B (z), \ \ \ \ z\in  V_a.$$
 We have $\rho_{a,b} (a)=b,$ and $B(\rho_{a,b})=B $.
 If $\rho$ is a local inverse  defined on a neighborhood $V$ of $a$
 and $\rho(a)=b,$ then it is routine to check that $$\rho=\rho_{a,b}$$ on
 $V\cap V_a.$ In this way, $\rho \sim \rho_{a,b}$.

The following notion introduced in \cite{GH3} plays
 an important role in the proof of Theorem \ref{submain}.
\begin{Def} \label{def1}
Let $\psi$ be an analytic function on a domain $\Omega.$ For two distinct points
 $z_1,z_2\in \Omega \backslash \psi^{-1}(\mathcal{E}_\psi)$ with
$\psi(z_1)=\psi(z_2),$ if there exist curves $\sigma_1 $ and $\sigma_2$ such that
\begin{itemize}
\item[\emph{(1)}]  $\sigma_i(0)=z_i$ for $i= 1,2,$ and $\sigma_1(1)=\sigma_2(1)$;
\item[\emph{(2)}]  $ \psi(\sigma_1(t))=\psi(\sigma_2(t))$ for $0\leq  t\leq    1;$
\item[\emph{(3)}] $ \sigma_i(t) \in  \Omega \backslash \psi^{-1} (\mathcal{E}_\psi) $,
 $i=1,2  $ and $0\leq     t< 1,$
\end{itemize}
then $z_1$ and $z_2$ are said to be directly $\psi$-glued in $\Omega,$ denoted by
$z_1 \overset{\psi }{ \backsimeq} z_2. $
\end{Def}
Roughly, two points are directly $\psi$-glued if
they are connected by a $\psi$-equal curve.
For example, set $\psi(z)=z^n~(n>1)$ and $\Omega=\mathbb{D}$.
 If $z_1\neq z_2$  and $z_1^n=  z_2^n$,
 then we have $z_1 \overset{\psi }{\backsimeq} z_2  $. In fact,
   defining
   $\sigma_1(t)=(1-t)z_1$ and $\sigma_2(t)=(1-t)z_2$,
  we have $\psi(\sigma_1(t))=\psi(\sigma_2(t))  $
  for all $ 0\leq     t \leq    1.$

 The following is from \cite[Definition 4.4]{GH4}.
\begin{Def} Let $\Omega$ be a subdomain of $\mathbb{D}$.
If there exists a finite sequence $z_0,z_1,\cdots,z_k$ in
$\Omega$ such that $z_j$ and $z_{j+1}$ are
directly $\psi$-glued in $\Omega $ for
$0\leq    j\leq     k-1$, and all related curves $\sigma_i$
 in Definition \ref{def1} can be chosen in $\Omega$, then
  $z_0$ and $z_k$ are said to be  $\psi$-glued on $\Omega.$

   A function $\psi$ is called gluable
   if for every $w'\in \mathbb{D}\backslash \psi^{-1}(\mathcal{E}_\psi)$,
 and every $w''\neq w'$ satisfying  $\psi(w'')=\psi(w'),$ the points
  $w'$ and $w''$ are  $\psi$-glued in $ \mathbb{D} $.

\end{Def}
The following is from  \cite[Propositions 4.5 and 4.6]{GH4}.
\begin{prop}\label{45}
A finite/thin Blaschke product is gluable.
Therefore, a proper map between two simply connected domains  is gluable.
\end{prop}

    For a holomorphic proper map $\phi: \Omega_1 \to \Omega_2$, it is a standard fact that there exists an integer $k \geq1$ such that for every $w \in \Omega_2$, the equation $\phi(z)=w$ has exactly $k$ solutions in $\Omega_1$, counting multiplicities. In this case, we call
    $\phi:\Omega_1\to \Omega_2 $
   is of  multiplicity $k$. Following \cite{GHLiLuo}, we present the following notion.

\begin{Def} For a holomorphic proper map $\phi:\Omega_1\to \Omega_2,$
   we say $\Omega_1$ is \emph{of multiplicity $k$} with respect to $\phi$ if $\phi:\Omega_1\to \Omega_2 $
   is of multiplicity $k$.\end{Def}

The following result is from \cite[Lemma 2.12]{GHLiLuo}.
\begin{lem} \label{continuation}
Let $\phi$ be a holomorphic function on a planar domain $V$.
Let $\sigma_1, \sigma_2: [0, 1) \to V$
be two disjoint   curves satisfying:
\begin{align*}
\phi(\sigma_1(s)) &= \phi(\sigma_2(s)),
\quad 0 \leq    s < 1, \quad \text{and} \\
\phi'(\sigma_1(s)) &\neq 0 \quad \text{and}
\quad \phi'(\sigma_2(s)) \neq 0, \quad 0 \leq     s < 1.
\end{align*}
For each $s \in [0, 1)$, let $\rho_s$ be the local inverse of $\phi$
 defined in a neighborhood of   $ \sigma_1(s) $  such that
$$
\rho_s(\sigma_1(s)) = \sigma_2(s) \quad \text{for all } 0 \leq    s < 1.
$$
Then all these $\rho_s$  extends analytically to each other along the curve
$\sigma_1$, that is,
 $\rho_0 \sim \rho_s$, $0 \leq     s < 1$.
\end{lem}

We need Lemma 5.8.3 in \cite{GH4}, restated as follows.
 \begin{lem}\label{gluel}
Given $z_0,z_1$ and $z_2$ in $\mathbb{D} \backslash B^{-1} (\mathcal{E}_B)$
with $ B(z_0)=B(z_1)=B(z_2)$, suppose there are three curves $\sigma_0,\sigma_1,\sigma_2$
such that
\begin{itemize}
\item[\emph{(1)}]  $\sigma_i(0)=z_i$ for $i=0,1,2,$ and $\sigma_1(1)=\sigma_2(1)$;
\item[\emph{(2)}]  $B(\sigma_0(t))=B(\sigma_1(t))=B(\sigma_2(t))$ for
$0\leq  t\leq    1 $,
 and $ \sigma_i(t) \in  \mathbb{D} \backslash B^{-1} (\mathcal{E}_B) $,
 $i=0,1,2, $ $0\leq     t< 1;$
\item[\emph{(3)}] $B'(\sigma_0(1))\neq 0.$
\end{itemize}
Then $\rho_{z_0,z_1} \sim \rho_{z_0,z_2}$.
\end{lem}
Remind that for a  thin Blaschke product $B$, both $ \mathcal{E}_B $
and $ B^{-1}(\mathcal{E}_B)$
are discrete in $\mathbb{D}$ (see \cite[Proposition 2.1]{GH3}).
For a thin/finite Blaschke product $B$,
 and for a  subdomain $\Omega$ of $\mathbb{D}$,
we define
   $$ \mathcal{S}_{B,\Omega}=\{ (w,z)\in  \Omega^2: B(z)=B(w)\notin
    \mathcal{E}_B \}.$$
  Specifically, for $\Omega=\mathbb{D}$,
    the associated Riemann surface
 $$\mathcal{S}_{B, \mathbb{D}}=\mathcal{S}_B
  =\{(z,w)\in \mathbb{D}^2: B(z)=B(w)\not\in \mathcal{E}_B \}. $$

We also need the
 following. The former part is directly from Proposition 5.3.2
 in \cite{GH4},  and its later part
 is a consequence of  Proposition \ref{sproper}.
    \begin{prop}\label{inversecomponent}
   Let $B$ be a finite Blaschke product. If $\mathcal{S}$ is a
   component of $\mathcal{S}_B$, then
   $\mathcal{S}^-:=\{ (w,z)\in  \mathbb{D}^2: (z,w)\in \mathcal{S}\}$
   is also a component of $\mathcal{S}_B$.
   Consequently, for each $t\in (0,1),$ if $\Omega$
   is a component of $\{|B|<t\}$  and   $\mathcal{T}$
   is a component   of $ \mathcal{S}_{B,\Omega}$, then
   $\mathcal{T}^-:=\{ (w,z)\in  \Omega^2: (z,w)\in \mathcal{T}\}$
  is a component of  $ \mathcal{S}_{B,\Omega}$.
     \end{prop}

For every subdomain $\Omega$ of $\mathbb{D}$,
 we define \textbf{the identity component} of
 $$\{(z,w)\in \Omega^2: B(z)=B(w)\not\in \mathcal{E}_B\}, $$
  to be the set $\{(z,z):z\in
  \Omega \backslash B^{-1}(\mathcal{E}_B)  \}$.
\begin{rem} \label{remxyyx}
 Let $B$ be a finite/thin Blaschke product.
We   have the following.
 \begin{itemize}
\item[\emph{(1)}]For two local inverses $\rho $ and $\tau$ of $B$,
  $\rho \sim \tau$
 if and only if there exists $z_0$ and $z_1$ in
 $\mathbb{D}\backslash B^{-1}(\mathcal{E}_B)$, $(z_0,\rho(z_0))$
 and $(z_1,\tau(z_1))$ belong to the same
 component of $\mathcal{S}_B$ (see \cite[Above Lemma 3.6]{GH4}).
 If $(z_0,\rho(z_0))$ belongs to a
 component $\mathcal{S}$ of $\mathcal{S}_B$, then for every
 $\widetilde{\rho}\sim \rho,$ $(z,\widetilde{\rho}(z)) \in \mathcal{S}$
 for every $z$  in the domain of $\widetilde{\rho}.$

 Consequently, the Riemann surface
 $\mathcal{S}_B$  has exactly two components
 if and only if every pair of non-identity local inverses
   $\rho $ and $\tau$   of $B$ admit analytic continuation to
   each other.
 \item[\emph{(2)}]  The equivalence in (1) still holds if $\mathbb{D}$
 is replaced by any   (necessarily  simply connected)
  component $\Omega$ of $\{|B|<t\}$, in which case we say that
   $\rho  $ can be analytically continuable to
 $ \tau$ along a curve lying entirely in
   $\Omega\backslash B^{-1}(\mathcal{E}_B)$.

      As a consequence, the second half of Proposition \ref{inversecomponent}
      yields the following.
        Let $x,y, z,w\in  \Omega\backslash B^{-1}(\mathcal{E}_B)$
satisfy $B(x)=B(y)$ and  $B(z)=B(w)$. Then
 $\rho_{x,y} $ is analytically continuable to
 $ \rho_{z,w}$
along some curve lying entirely in $\Omega\backslash
 B^{-1}(\mathcal{E}_B)$
 if and only if
$\rho_{y,x} $  is analytically continuable to $ \rho_{w,z} $
along some curve in $\Omega\backslash B^{-1}(\mathcal{E}_B)$.
\item[\emph{(3)}] If $\mathcal{S}$ is a component of  $\mathcal{S}_B$, then
for  each $z\in \mathbb{D}\backslash B^{-1}(\mathcal{E}_B),$
the cardinality $n$ of the fiber $\{w: (z,w)\in \mathcal{S}\}$ is a constant
 (see \cite[Section 3]{GH3}), 
 called the multiplicity of  $\mathcal{S}$,
  denoted by $n=\sharp \mathcal{S}.$
  If  $\mathcal{S}$  contains a
  point  $(z, \rho(z))$ for some
  $z\in \mathbb{D}\backslash B^{-1}(\mathcal{E}_B)$, then
   this multiplicity $ \sharp \mathcal{S} $ coincides  with
 the cardinality of $$  \{w: (z,w)\in \mathbb{D}^2,  B(z)=B(w)\not\in \mathcal{E}_B  ~~
 \mathrm{and} ~~ \rho_{z,w}\sim \rho \}. $$
\end{itemize}
  \end{rem}

\begin{Def} \label{pure}
 Let $B$ be a finite/thin Blaschke product, $0<t<1$ and
let $\Omega$ be a component of $\{|B|<t\}$. If
$$\{(z,w)\in \Omega^2: B(z)=B(w)\not\in \mathcal{E}_B\}$$
contains exactly two components, then we say  $\Omega$ is \emph{purely gluable}
with respect to $B$.
    \end{Def}

We have the following result.
\begin{prop} \label{newequi}
Let $0<t<1.$  For a
 component $\Omega$   of $\{|B|<t\}$,
 $\Omega$ is purely gluable with respect to
$B$ if and only if excluding the identity map,
every local inverse defined on a
subdomain of  $\Omega$ can be analytically continuable  to any other
such non-identity local inverse
 along some curve in $\Omega\backslash B^{-1}(\mathcal{E}_B)$.\end{prop}
 \begin{proof} In fact, by Propositions \ref{sproper} and \ref{atomz},
   the component $\Omega$   of $\{|B|<t\}$   is  simply connected
   and there exists  a conformal map $G:\Omega\to \mathbb{D}$
    such that
    $$B|_{\Omega}=  t \varphi\circ G(z), \ \ \ \  z\in \Omega $$
   for some finite Blaschke product $\varphi$.
   In this case, the Riemann surface
   \begin{align*}
   &\{(z,w)\in \Omega^2: B(z)=B(w)\not\in \mathcal{E}_B\}\\
   &= \{(z,w)\in \Omega^2: t\varphi(G(z))=t\varphi(G(w))
   \not\in \mathcal{E}_B\}
   \end{align*}
   is biholomorphic to $$\{(G(z),G(w))\in  \mathbb{D}^2: t\varphi(G(z))=t\varphi(G(w))
   \not\in \mathcal{E}_B\},$$ which simplifies to
   $\{(\widetilde{z},\widetilde{w})\in  \mathbb{D}^2: \varphi(\widetilde{z})= \varphi(\widetilde{w})
   \not\in \mathcal{E}_\varphi\}$. This shows that
   $\{(z,w)\in \Omega^2: B(z)=B(w)\not\in \mathcal{E}_B\}$ is
   biholomorphic to $\mathcal{S}_\varphi.$
   Combining this with Remark \ref{remxyyx}(1),
   we immediately obtain the desired conclusion.
\end{proof}

 The following proposition comes   from \cite{DPW,DSZ} and \cite{GH4}.
  \begin{prop}\label{thinfin}
   For  thin/finite Blaschke products, we have the following properties.
\begin{itemize}
\item[\emph{(1)}] If $B$ is a thin/finite Blaschke product and $\mathbb{D}$
is purely gluable, then $\mathcal{S}_B$ has exactly two components.
\item[\emph{(2)}] If $B$ is a finite Blaschke product and $\mathcal{S}_B$ has
exactly two components,
then $M_B$ on $L_a^2(\mathbb{D})$ has exactly two minimal reducing subspaces.
\item[\emph{(3)}]  If $B$ is a thin Blaschke product
and $\mathcal{S}_B$ has exactly two components,
then $M_B$ on $L_a^2(\mathbb{D})$ has no nontrivial reducing subspaces.
\end{itemize}
\end{prop}
\begin{proof} By definition,   $\mathcal{S}_B$ has only two components
 provided that $\mathbb{D}$ is purely gluable. Thus $(1)$ is straightforward.

 For (2), just note that for a finite Blaschke product $B$,
 the number of  components of $\mathcal{S}_B$ equals that of minimal reducing subspaces
 for $M_B$ on $L_a^2(\mathbb{D})$ (see \cite{DPW,DSZ}).

 It remains to deal with (3). Now assume $B$ is a thin Blaschke
product.
By Theorem 5.3.1 in
\cite{GH4}, for a thin Blaschke
product $B$, $M_B$ on $L_a^2(\mathbb{D})$ is irreducible if and only if
$\mathcal{V}^*(B)=\mathbb{C}I $,  if and only if
 the Riemann surface
$\mathcal{S}_B$ has a component of finite multiplicity,
distinct from the identity component $ \{(z,z): B(z)\not\in \mathcal{E}_B\}.$
When $\mathcal{S}_B$ has exactly
two components, the non-identity component $\mathcal{S}$  must
be of infinite multiplicity. Then  $M_B$ is irreducible.
\end{proof}

\subsection{Critical points of finite Blaschke products}
This subsection presents properties related to
the critical points
 of finite Blaschke products.

We say a curve
$\gamma$ is    \emph{non-self-intersecting} if $\gamma$ is injective, that is,
$\gamma(s)\neq \gamma(t)$   for $0\leq     s<t\leq    1 $.
Part (1) of the following theorem appears in
   \cite[Lemma 2.8]{GHLiLuo}.

\begin{thm} \label{e1}
Let $ B $ be a finite Blaschke product
  of order $ n $ with its critical set  $  \mathcal{E}_B $.
    Then we have the following.
\begin{itemize}
\item[\emph{(1)}]  If $ \gamma $ is a non-self-intersecting
    curve such that $ \gamma([0,1]) \subseteq \overline{\mathbb{D}}
    \backslash     \mathcal{E}_B  $, then there exist  exactly $ n $ curves
    $ \sigma_1, \cdots, \sigma_n $
    satisfying:
    $$
    B(\sigma_i(t)) = \gamma(t) \quad \text{and}
     \quad \sigma_i(t) \neq \sigma_j(t)  \ \  \text{for } \ \
      0 \leq     t \leq    1;\
      \ i \neq j.
    $$
    \item[\emph{(2)}] If $ \gamma $ is a non-self-intersecting curve
    such that
    $$
    \gamma[0,1) \subseteq \overline{\mathbb{D}} \backslash
    \mathcal{E}_B  \quad \text{and} \quad \gamma(1) \in  \mathcal{E}_B,
    $$
    then there exist curves $ \sigma_i ~(1\leq    i\leq   n) $
    satisfy $ B(\sigma_i(t)) = \gamma(t) $ for $ 0 \leq   t \leq   1 $
    and $\sigma_i(t) \neq \sigma_j(t) \ \text{for }
      0 \leq     t < 1,
      \, i \neq j$. Furthermore,
    if $ p $ is a critical point of $ B $ with $ B(p) =
     \gamma(1) $ and of order $ k-1 $, there exist exactly $k$ distinct
      indices $ i_1, \cdots, i_{k } $ such that:
    $$
    \sigma_{i_j}(1) = p, \ \ \ \  1 \leq   j \leq     k.
    $$
\end{itemize}
\end{thm}

Recall that if a nonconstant holomorphic function \(f\) on a domain \(\Delta\) has a zero at \(a \in \Delta\), then there exist a holomorphic function \(g\) on \(\Delta\) with \(g(a)\neq 0\) and a positive integer \(m\) such that \(f(z)=(z-a)^m g(z)\) for \(z\in\Delta\). The integer \(m\) is the order of the zero of \(f\) at \(a\), denoted by \(\mathrm{order}_f(a)\).
For the proof of Theorem \ref{e1},
we need the following form of B\"ottcher's theorem, see
 \cite[Theorem 9.1]{Mi} for example.
\begin{thm}[B\"ottcher] \label{bott}
 Let $\psi$ be a holomorphic function
 on a neighborhood of a point $p$.  If
$\psi'(p)=0$ and $\psi \not\equiv 0,$ then
there exist a simply connected domain $\Delta$ containing $p$, a number
$\varepsilon>0$, and a biholomorphic map
 $g: \Delta \to \varepsilon \mathbb{D}$
 so that
 $$\psi(z)-\psi(p)= [g(z)]^k, \ \ \ z\in \Delta,$$
 where $k=\mathrm{order}_{\psi'}(p)+1. $ In this case,
 $\psi|_{\Delta}:\Delta \to \psi(\Delta) $
  is a holomorphic proper map of multiplicity $k.$
\end{thm}
\begin{proof}[\textbf{\emph{Proof of Theorem \ref{e1}.}}]
    (1) First we assume that $ \gamma $ is a non-self-intersecting
    curve such that $ \gamma[0,1] \subseteq \overline{\mathbb{D}}
    \backslash   \mathcal{E}_B  $. For each $ t_0 \in [0,1] $ with $ \gamma(t_0)
    \notin \mathcal{E}_B $, there exist $ n $ points $ z_1,
     \cdots, z_n $ such that $ B(z_i) = \gamma(t_0) $
     for $ 1 \leq     i \leq     n $. Since $ B'(z_i) \neq 0 $ for each $i,$
     by the Inverse Function Theorem, there exist
     neighborhoods $ U_i $ of $ z_i $ such that the map
     $ B|_{U_i}: U_i\to B(U_i) $ is biholomorphic.
     On each $ U_i $, we can define $ \sigma_i $ as:
    $$
    \sigma_i(t) = B^{-1}(\gamma(t)) \ \ \ \ (t\in   [t_0-\delta,   t_0+\delta])
    $$
    for some $\delta>0.$
  By compactness of \([0,1]\), we can cover it with finitely many intervals  $ \Delta_1, \cdots, \Delta_N $ of the form \([t_0-\delta, t_0+\delta]\).
     For each $ \Delta_j $, there correspond $ n $ curves $ \sigma_{i,j} $
     satisfying $ B(\sigma_{i,j}(t)) = \gamma(t) $ for $ t \in \Delta_j $ with
     $i=1,\cdots,n$.
     By pasting  these curves over overlapping intervals,
      we obtain the desired curves $ \sigma_i $.
      It is routine to verify that the values
      $ \sigma_i(0) $ are distinct for $i=1,2,\cdots,n$,
      and by continuity, the curves $ \sigma_i(t) $ remain distinct
       for all $ t \in [0,1] $.

    (2)   Next we assume $ \gamma $ is
a non-self-intersecting curve
 such that
 $$   \gamma[0,1) \subseteq \overline{\mathbb{D}} \backslash
    \mathcal{E}_B  \quad \text{and} \quad \gamma(1) \in  \mathcal{E}_B.
    $$  For every $t^*\in (0,1),$
    by the discussions in (1), there exist $n$ curves $\sigma_i$ defined
    on $[0,t^*]$ such that $B(\sigma_i(s))=\gamma(s),~0\leq    s\leq     t^*.$ The idea is to
    find the preimages of $\gamma(s)~(t^*\leq     s \leq     1)$ when $t^*$ is enough close to $1.$

  Before continuing, let us make an observation.  Let $p$ be a critical point with
   $B(p)=\gamma(1)$ and
   critical order $k-1$.
    By  B\"ottcher's theorem, there exist a simply connected domain $\Delta$ containing $p$,
$\varepsilon>0$, and a biholomorphic map
 $g: \Delta \to \varepsilon \mathbb{D}$
 so that $g(0)=0$ and
 $$B(z)-B(p)= [g(z)]^k \ \ \ \  z\in \Delta.$$
For any curve $\gamma_0:[0,1]\to \mathbb{C}$  in the
 $ \frac{\varepsilon}{2} $-neighborhood of  $ B(p)  $
 with $\gamma_0(1)=B(p),$
 the preimage $B^{-1}(\gamma_0) \cap \Delta$ consists of $n$ distinct
 analytic arcs $\{\Gamma_i\}_{1\leq     i\leq     n}$, and
   exactly $k$ of these arcs terminate at $p$.
          Note that if $p$ is a regular point,
          then the above analysis remains valid with $k=1$.

    Let $\gamma_0 = \gamma|_{[t^*,1]}$ where $t^*$ is sufficiently close to 1. Then
  $B^{-1}(\gamma|_{[0,t^*]})$ splits into $n$ curves,
  and each curve connects
   to one of the $\Gamma_i$.   The conclusions in (2) hold,
    thereby completing the proof.
 \end{proof}

We also have a consequence of Theorem \ref{e1}.
\begin{cor}\label{e2}
Let $B$ be a finite Blaschke product
 and fix  $0<t<1$. Suppose
$\Omega_1,\cdots,\Omega_N $ are all the
 distinct components of $\{|B|<t\}$.
For every point  $p\in \partial \Omega_{1} $,
 $p  $ is a critical point of $B$ if and only if
 $p$   lies in the union
 $\bigcup\limits_{2\leq   j\leq   N}\partial \Omega_{1}
 \cap \partial \Omega_{j} .$
\end{cor}

\begin{figure}[htbp]
\centering
\tikzset{every picture/.style={line width=0.75pt}} 

\begin{tikzpicture}[x=0.75pt,y=0.75pt,yscale=-1,xscale=1]

\draw  [line width=3] [line join = round][line cap = round] (204.53,183.89) .. controls (204.53,183.89) and (204.53,183.89) .. (204.53,183.89) ;
 \draw   (200,185.58) .. controls (219.5,175.08) and (244.5,179.58) .. (279,210) .. controls (313.5,240.42) and (309.5,240.08) .. (279,270) .. controls (248.5,299.92) and (211.5,309.58) .. (191.5,279.58) .. controls (171.5,249.58) and (180.5,196.08) .. (200,185.58) -- cycle ;
\draw  [line width=3] [line join = round][line cap = round] (199.53,288.89) .. controls (199.53,288.89) and (199.53,288.89) .. (199.53,288.89) ;
\draw    (204.76,183.79) .. controls (206.15,191.19) and (210.48,193.62) .. (213.41,196.06) .. controls (216.35,198.49) and (224.74,197.8) .. (227.68,211.59) ;
\draw    (267.83,236.43) .. controls (274.21,240.45) and (278.93,238.94) .. (282.71,238.45) .. controls (286.49,237.95) and (291.7,231.34) .. (303.8,238.58) ;
\draw    (199.82,288.87) .. controls (207.26,287.74) and (209.84,283.51) .. (212.38,280.66) .. controls (214.91,277.81) and (214.52,269.4) .. (228.41,266.95) ;
\draw  [line width=3] [line join = round][line cap = round] (303.5,238.45) .. controls (303.5,238.45) and (303.5,238.45) .. (303.5,238.45) ;
\draw  [line width=3] [line join = round][line cap = round] (400.46,296.9) .. controls (400.46,296.9) and (400.46,296.9) .. (400.46,296.9) ;
\draw   (405.04,295.36) .. controls (385.19,305.19) and (360.36,299.84) .. (326.91,268.27) .. controls (293.47,236.7) and (297.45,237.16) .. (328.96,208.3) .. controls (360.46,179.44) and (397.76,171.04) .. (416.73,201.7) .. controls (435.7,232.36) and (424.88,285.53) .. (405.04,295.36) -- cycle ;
\draw  [line width=3] [line join = round][line cap = round] (409.03,192.13) .. controls (409.03,192.13) and (409.03,192.13) .. (409.03,192.13) ;
\draw    (400.22,296.99) .. controls (399.08,289.54) and (394.84,286.97) .. (391.99,284.44) .. controls (389.14,281.9) and (380.73,282.31) .. (378.26,268.42) ;
\draw    (339.68,241.95) .. controls (333.44,237.71) and (328.67,239.06) .. (324.88,239.43) .. controls (321.08,239.79) and (315.65,246.23) .. (303.8,238.58) ;
\draw    (408.74,192.14) .. controls (401.26,193.01) and (398.53,197.16) .. (395.9,199.92) .. controls (393.27,202.68) and (393.37,211.09) .. (379.41,213.07) ;

\fill  (204.53,183.89) circle (1.5pt); %
\fill  (199.53,288.89) circle (1.5pt); %
\fill  (303.5,238.45) circle (1.5pt); %
\fill  (400.46,296.9) circle (1.5pt); %
\fill (409.03,192.13) circle (1.5pt); %
\draw (268.23,217.5) node [anchor=north west][inner sep=0.75pt]  [rotate=-358.86]  {$\sigma_1$};
\draw (324.23,221.5) node [anchor=north west][inner sep=0.75pt]  [rotate=-358.86]  {$\sigma_2$};
\draw (373.3,310.97) node [anchor=north west][inner sep=0.75pt]  [rotate=-358.86]  {$\Omega_1$};
\draw (227.23,310.5) node [anchor=north west][inner sep=0.75pt]  [rotate=-358.86]  {$\Omega_2$};
\end{tikzpicture}
\caption{Lemma \ref{component3} }
\label{Fig-2}
\end{figure}

The following lemma is included in \cite[Corollary 2.9]{GHLiLuo}.
\begin{lem}\label{component3}
Let $B$ be a finite Blaschke product with the critical set
 $\mathcal{E}_B,$
 and fix  $0<t<1$. Suppose
$\Omega_1, \Omega_2$ are distinct components of $\{|B|<t\}$.
For any point  $p\in \partial \Omega_{1} \cap \partial \Omega_{2} $,
  let $\gamma$ be a curve  such that
$ \gamma [0,1) \in t\mathbb{D}\backslash \mathcal{E}_B,  $
and $\gamma(1)=B(p).$  Then there exists
 a unique pair of curves $\sigma_1,\sigma_2$
  such that for $0\leq s<1$, $\sigma_i(s)\in \Omega_i$,
$$B(\sigma_i(s))=\gamma(s) \quad \text{and} \quad
 \sigma_i(1)=p,\quad i=1,2.$$
 Consequently, $p$ is a  critical point of $B$.

Furthermore, if in addition, \(\gamma\) is non-self-intersecting and there exists a curve \(\sigma\) such that \(B(\sigma(s))=\gamma(s)\) with \(\sigma(s)\in\Omega_1\) for \(0\leq  s<1\), then
  either $\sigma$ is disjoint from $\sigma_1$,
 or $\sigma_1(s)=\sigma(s),0\leq     s\leq     1.$
\end{lem}

The following describes the structure of the boundaries of components  of $\{|B|<t\}$ for
 $0<t<1$.  This theorem will be used in our induction step in proving Theorem \ref{submain}.
\begin{thm}\label{component2}
Let $B$ be a finite Blaschke product, and fix  $0<t<1$. Denote by
$\{\Omega_i \}_{i=1}^k$  the components of the sublevel set $\{|B|<t\}$ where
$k\geq  2$.
Then these components satisfy the following  properties.
\begin{itemize}
\item[\emph{(1)}] \textbf{Critical adjacency.} For $i\neq j$, the intersection
 $\partial \Omega_i \cap \partial \Omega_j$ contains at most one point.
 Moreover, if $\partial \Omega_i \cap \partial \Omega_j$ is nonempty,
  then it  consists precisely of a critical point of $B$.
\item[\emph{(2)}] \textbf{Cyclic connectivity.}
Under a reordering  $\{\Omega_i\}_{1\leq    i\leq     k}$
  satisfying
$\partial \Omega_{i} \cap \partial \Omega_{i+1} \neq \emptyset$
 for $1\leq    i\leq     k-1,$ exactly one of the following holds:

 \emph{(a)}    $\partial  \Omega_k\cap \partial  \Omega_1=\emptyset$
 (open chain configuration);

  \emph{(b)}   there exists a unique point $p$ in
 $ \bigcap\limits_{1 \leq  i <j\leq    k} \partial \Omega_i \cap \partial \Omega_j $  
 (star-shaped intersection).
\end{itemize}
\end{thm}

\begin{figure}[htbp]
\centering
\tikzset{every picture/.style={line width=0.75pt}}
\begin{tikzpicture}[x=0.75pt,y=0.75pt,yscale=-1,xscale=1]
\draw  [line width=3] [line join = round][line cap = round] (305.95,283.39) .. controls (305.95,283.39) and (305.95,283.39) .. (305.95,283.39) ;
\draw   (250.53,187.7) .. controls (270.53,177.7) and (303.07,193.2) .. (314.73,212) .. controls (326.4,230.8) and (300.58,230.74) .. (293.73,244) .. controls (286.89,257.25) and (295.75,266.02) .. (303.73,278) .. controls (311.72,289.98) and (291.29,290.9) .. (270.93,291.7) .. controls (250.58,292.5) and (233.38,283.98) .. (223.33,268.9) .. controls (213.28,253.82) and (230.53,197.7) .. (250.53,187.7) -- cycle ;
\draw   (391.93,297.4) .. controls (381.03,313.51) and (378.93,325.4) .. (367.93,336.4) .. controls (356.93,347.4) and (353.55,350.72) .. (340.33,352.28) .. controls (327.12,353.84) and (313.87,357.18) .. (300.54,347.37) .. controls (287.2,337.57) and (289.03,337.28) .. (285.83,333.28) .. controls (282.63,329.28) and (264.53,311.22) .. (272.93,295.4) .. controls (281.33,279.57) and (295.58,315.48) .. (306.93,308.5) .. controls (318.29,301.52) and (303.73,287.38) .. (306.53,282.1) .. controls (309.33,276.82) and (348.99,270.53) .. (326.74,242.16) .. controls (304.49,213.79) and (328.95,161.6) .. (350.93,159.7) .. controls (372.92,157.8) and (412.12,194.17) .. (416.47,235.22) .. controls (420.83,276.28) and (402.83,281.28) .. (391.93,297.4) -- cycle ;
\draw  [line width=3] [line join = round][line cap = round] (318.03,219.89) .. controls (318.03,219.89) and (318.03,219.89) .. (318.03,219.89) ;
\draw  [line width=3] [line join = round][line cap = round] (277.87,291.44) .. controls (277.87,291.44) and (277.87,291.44) .. (277.87,291.44) ;
\draw (242.12,222.99) node [anchor=north west][inner sep=0.75pt]  [rotate=-1.67]  {$\Omega_1$};
\draw (275.46,240.99) node [anchor=north west][inner sep=0.75pt]  [rotate=-1.67]  {$\gamma_1$};
\draw (308.52,246.06) node [anchor=north west][inner sep=0.75pt]  [rotate=-1.67]  {$\Omega$};
\draw (374.66,241.26) node [anchor=north west][inner sep=0.75pt]  [rotate=-1.67]  {$\Omega_2$};
\draw (338.52,247.12) node [anchor=north west][inner sep=0.75pt]  [rotate=-1.67] {$\gamma_2$};
\draw (314.46,278.66) node [anchor=north west][inner sep=0.75pt]  [rotate=-1.67]  {$q$};
\draw (325.46,207.66) node [anchor=north west][inner sep=0.75pt]  [rotate=-1.67]  {$p$};

\fill (277.87,291.44) circle (1.5pt);
\fill  (305.95,283.39) circle (1.5pt); 
\fill  (318.03,219.89) circle (1.5pt); 
\end{tikzpicture}
\caption{Excluded cases in Theorem \ref{component2}}
\label{Fig-1}
\end{figure}
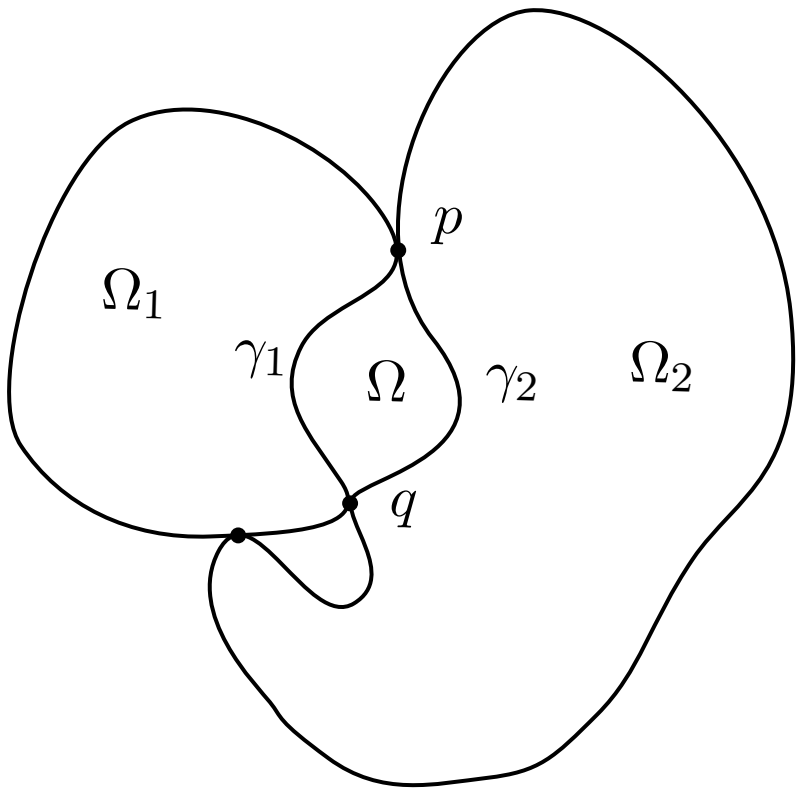

\begin{proof} By Lemma \ref{component3}, if
$\partial \Omega_i \cap \partial \Omega_j\neq \emptyset $
for $i\neq j,$ then every point in $\partial \Omega_i \cap \partial \Omega_j$
is a critical point.  We now  prove that
$\partial \Omega_1 \cap \partial \Omega_2$ contains at most one point.
  To reach a contradiction, we assume that $\partial \Omega_1
   \cap \partial \Omega_2$ is not empty, and
    $\partial \Omega_1 \cap \partial \Omega_2$ contains
 two distinct points $p$ and $q$.
By Lemma \ref{Jordancurve}, \(\partial \Omega_1\) is a Jordan curve homeomorphic to the unit circle, and then we may choose two adjacent critical points \(p\) and \(q\) on \(\partial \Omega_1\). Since both \(\partial \Omega_1\) and \(\partial \Omega_2\) are homeomorphic to the unit circle, there exist curves \(\gamma_1 \subseteq \partial \Omega_1\) and \(\gamma_2 \subseteq \partial \Omega_2\), each joining \(p\) to \(q\). We may further assume that the domain \(\Omega\) bounded by the closed curve \(\gamma_1 \cup \gamma_2\) is simply connected and satisfies the condition that
$$ \Omega \cap \Omega_1 = \Omega \cap \Omega_2 = \emptyset
\ \ \ \ (\text{see \  Figure}\  \ref{Fig-1}).$$
Since $|B| = t$ holds both on $\gamma_1$ and $\gamma_2$,  by the maximum principle, on $\Omega$ we have $|B| < t$. Note that $\partial \Omega = \gamma_1 \cup \gamma_2$ and $\partial \Omega \cap \partial \Omega_1$ contains the curve $\gamma_1$.
Noting that $\Omega$ is also a component of $\{|B|<t\}$,
by Lemma \ref{component3}   each point in
 $\partial \Omega \cap \partial \Omega_1$ is a critical point.
  In particular, each point on $\gamma$ is a critical point.
   This gives that $B$ has infinitely many critical points,
    which contradicts to  Lemma  \ref{critical}.
Therefore  $\partial \Omega_1 \cap \partial \Omega_2$
 contains at most one point.

\begin{figure}[htbp]
\centering
\tikzset{every picture/.style={line width=0.75pt}} 
\begin{tikzpicture}[x=0.75pt,y=0.75pt,yscale=-1,xscale=1]
\draw   (343,297.78) .. controls (353.39,306.44) and (371.56,343.61) .. (351.8,343.19) .. controls (332.04,342.76) and (317.98,333.42) .. (302.86,315.82) .. controls (287.74,298.21) and (298.68,290.42) .. (307.86,290.33) .. controls (317.04,290.23) and (332.61,289.13) .. (343,297.78) -- cycle ;
\draw   (381.2,197.79) .. controls (378,202.15) and (377.34,203.73) .. (370.8,211.44) .. controls (364.27,219.15) and (362.8,220.71) .. (361.29,216.97) .. controls (359.78,213.24) and (350.34,198.25) .. (354.56,185.68) .. controls (358.78,173.11) and (375.29,156.76) .. (380.22,159.95) .. controls (385.15,163.13) and (384.56,166.56) .. (384.96,172.55) .. controls (385.37,178.55) and (384.41,193.43) .. (381.2,197.79) -- cycle ;
\draw   (322.42,262.87) .. controls (318.71,271.08) and (307.36,291.52) .. (300.81,291.55) .. controls (294.26,291.58) and (287.01,285.01) .. (279.97,267.6) .. controls (272.93,250.18) and (283.34,232.69) .. (290.22,222.66) .. controls (297.11,212.63) and (314.07,205.92) .. (316.89,209.71) .. controls (319.71,213.5) and (327.23,221.53) .. (326.17,229.27) .. controls (325.11,237) and (326.13,254.65) .. (322.42,262.87) -- cycle ;
\draw  [line width=3] [line join = round][line cap = round] (362.38,219.21) .. controls (362.38,219.21) and (362.38,219.21) .. (362.38,219.21) ;
\draw   (339.47,186.4) .. controls (348.59,203.71) and (367.45,217.09) .. (360.57,222.79) .. controls (353.7,228.48) and (326.49,216.36) .. (321.23,212.71) .. controls (315.97,209.06) and (306.69,202.44) .. (306.13,190.08) .. controls (305.58,177.73) and (330.35,169.1) .. (339.47,186.4) -- cycle ;
\draw   (413.97,204.65) .. controls (431.19,211.12) and (436.28,225.37) .. (431.49,232.9) .. controls (426.7,240.43) and (412.5,247.37) .. (399.77,246.62) .. controls (387.04,245.87) and (365.21,229.4) .. (361.17,223.48) .. controls (357.14,217.57) and (396.75,198.17) .. (413.97,204.65) -- cycle ;
\draw   (398.92,322.58) .. controls (412.11,325.53) and (433.54,343.6) .. (419.5,356.28) .. controls (405.46,368.97) and (382.86,366.08) .. (371.43,356.93) .. controls (360,347.78) and (356.12,336.24) .. (364.25,331.97) .. controls (372.37,327.69) and (385.72,319.62) .. (398.92,322.58) -- cycle ;
\draw   (251.14,260.81) .. controls (244.52,255.28) and (230.78,229.54) .. (241.8,228.37) .. controls (252.82,227.2) and (261.54,232.96) .. (271.66,244.66) .. controls (281.79,256.37) and (276.48,262.86) .. (271.39,263.62) .. controls (266.3,264.38) and (257.76,266.34) .. (251.14,260.81) -- cycle ;
\draw   (217.65,246.92) .. controls (210.04,245.75) and (196.36,234.18) .. (202.92,223.89) .. controls (209.47,213.59) and (222.3,214) .. (229.55,219.82) .. controls (236.79,225.63) and (240.08,233.76) .. (235.99,237.48) .. controls (231.9,241.21) and (225.27,248.09) .. (217.65,246.92) -- cycle ;
\draw   (463.89,268.23) .. controls (460.16,281.23) and (433.45,312.84) .. (425.88,294.58) .. controls (418.31,276.32) and (413.51,253.87) .. (421.64,242.88) .. controls (429.77,231.9) and (454.56,235.28) .. (458.59,241.04) .. controls (462.63,246.8) and (467.62,255.24) .. (463.89,268.23) -- cycle ;
\draw   (461.73,329.9) .. controls (463.98,348.4) and (448.25,351.72) .. (436.45,349.55) .. controls (424.65,347.37) and (415.44,332.73) .. (419.21,318.58) .. controls (422.98,304.44) and (431.98,296.23) .. (439.16,301.95) .. controls (446.35,307.66) and (459.48,311.4) .. (461.73,329.9) -- cycle ;

\fill (362,220.24) circle (1.5pt); 
\draw (310.55,152.15) node [anchor=north west][inner sep=0.75pt]  [rotate=-1.53] [align=left] {$\Omega_1$};
\draw (373.15,143.24) node [anchor=north west][inner sep=0.75pt]  [rotate=-359.39] [align=left] {$\Omega_2$};
\draw (340,249) node [anchor=north west][inner sep=0.75pt]  [rotate=-2.85] [align=left] {$\Omega_\wp$};
\draw (290.77,249.71) node [anchor=north west][inner sep=0.75pt]  [rotate=-2.85] [align=left] {$\Omega_k$};
\draw (283.65,323.24) node [anchor=north west][inner sep=0.75pt]  [rotate=-359.39] [align=left] {$\Omega_{k-1}$};
\end{tikzpicture}
\caption{A case excluded by Cyclic connectivity
 (Theorem \ref{component2})}
\label{Fig-13}
\end{figure}
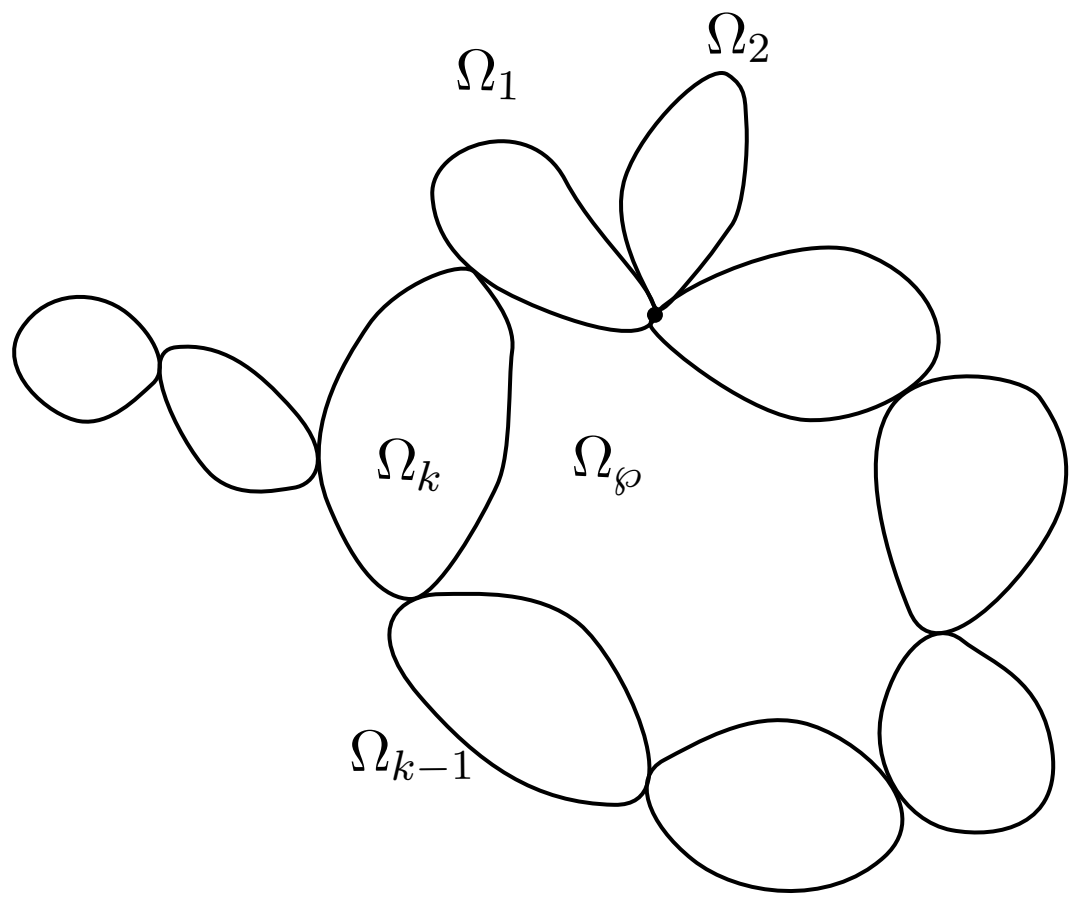

Finally, we turn to
 the proof of  (2). For this,
  assume that $\partial \Omega_i \cap \partial \Omega_{i+1} \neq
\emptyset $ ($1 \leq     i \leq    k-1$) and (b) fails;
that is, for $i \neq j$, not all of $\partial \Omega_i \cap \partial
\Omega_j$ equals a singlet $\{p\}$ for a fixed point $p$.
Suppose conversely that $\partial \Omega_k \cap \partial
\Omega_1 \neq \emptyset$. Since each curve $\partial \Omega_i$
($1 \leq     i \leq    k$) is topologically a circle,
a simple geometric observation shows that there exists a contour $\wp$ contained
 in $\bigcup\limits_{1 \leq     i \leq     k} \partial \Omega_i$ so that
the  interior $\Omega_\wp$ of $\wp$ has no intersection
 with any $\Omega_i$.
But $\Omega_\wp$ is a component of $\{|B| < t\}$.

To complete the proof of Theorem \ref{component2},
we now show that $B$ must have infinitely many critical
 points, which yields the desired contradiction.
In fact, there exists some $\Omega_{i_0}$ such that
 $\partial \Omega_\wp \cap \Omega_{i_0}$ contains
 a non-degenerate curve, and thus
 infinitely many distinct points.
 But by Lemma \ref{component3}, every point
 in this intersection is a critical point of $B$. However, $B$ has only finitely many critical points, a contradiction.
This completes the proof of Theorem \ref{component2}.
\end{proof}

  As we will show, the following lemma, which establishes
 a key property for finite Blaschke products,
  plays a critical role in the proof of Theorem~\ref{main}
  by serving as the foundational step in
  our approximation argument for thin Blaschke products.

 \begin{lem}\label{componenthigh}
Let $B$ be a finite Blaschke product and let
 $0<t<1$. Suppose
$\Omega_1$ and  $\Omega_2$  are distinct components of $\{|B|<t\}$.
If the sets $\{B(z): z \in \Omega_1,B'(z)=0\}$
 and  $\{B(z): z \in \Omega_2,B'(z)=0\}$ are disjoint, then for any
$z_1,z_2\in \Omega_1$ and $w_0 \in \Omega_2$ with
$$B(z_1)=B(z_2)=B(w_0)\not\in  \mathcal{E}_B ,$$
we have   $\rho_{w_0,z_1}\sim \rho_{w_0,z_2}.$
\end{lem}
\begin{proof}
Let $B$ be a finite Blaschke product.  Suppose   $0<t<1$, and let
$\Omega_1$ and  $\Omega_2$  be distinct components of $\{|B|<t\}$.

  By  Proposition \ref{sproper}, both $\Omega_1$ and  $\Omega_2$
    are simply connected domains. Assume
$z_1,z_2\in \Omega_1$ and $w_0 \in \Omega_2$ such that
$$B(z_1)=B(z_2)=B(w_0)\not\in  \mathcal{E}_B .$$
There are two cases to distinguish.

\textbf{Case I.} \emph{$z_1$ and $z_2$ are directly $B$-glued in $\Omega_1$.}
In this case, we will prove that $\rho_{w_0,z_1}\sim \rho_{w_0,z_2}.$

 By definition of directly $B$-gluedness, we have
\begin{itemize}
\item[(1)]  $\sigma_i(0)=z_i$ for $i =1,2,$ and $\sigma_1(1)=\sigma_2(1)$;
\item[(2)] $B(\sigma_1(s)) = B(\sigma_2(s))$ for $0 \leq s \leq 1$;
\item[(3)] $\sigma_i(s) \in \Omega_1  \backslash B^{-1}(\mathcal{E}_B)$,
 $i=1,2  $  for  $0\leq  s<  1.$
 \end{itemize}
 By Lemma \ref{proper}, $B|_{\Omega_2}: \Omega_2\to t\mathbb{D}$ is
 a proper map
   onto $t\mathbb{D}$, and $B(\Omega_2)=B(\Omega_1)=t\mathbb{D}$.   Therefore,
 the image of $B|_{\Omega_2}$ covers the image of $B$ on
 the curve $\sigma_1$. By Lemma \ref{component3},
   There exists a curve $\sigma_0$ satisfying $\sigma_0(0)=w_0$ and $\sigma_0(s) \in \Omega_2$ for $0 \leq s \leq 1$, with $B(\sigma_0(s)) = B(\sigma_1(s)) = B(\sigma_2(s))$ for $0 \leq s \leq 1.$
 Note that by the above Condition (2) and the condition $\sigma_1(1)=\sigma_2(1)$ in (1),
 $\sigma_1(1)$ is a critical point.  That is, $B'(\sigma_1(1))=0.$
By assumption, the critical value sets $$\{B(z): z \in \Omega_1,B'(z)=0\}$$
 and  $\{B(z): z \in \Omega_2,B'(z)=0\}$ are disjoint, and
 $B(\sigma_0(1))= B(\sigma_1(1))$. This implies
  $B'(\sigma_0(1))\neq 0.$ Thus  the conditions in
Lemma \ref{gluel} are satisfied to force
$\rho_{z_0,z_1}\sim \rho_{z_0,z_2} $.

\textbf{Case II.} \emph{General case.}

In general,
by Proposition \ref{sproper},
  $\Omega_1$  is simply connected and it follows from
Lemma  \ref{proper}  that the map $B: \Omega_1 \rightarrow t \mathbb{D}$
 is a proper map.
 A further application of Proposition  \ref{sproper} shows that
  after composing two conformal maps,
 $ B|_{\Omega_1}$ becomes a finite Blaschke product.
 Furthermore, by Lemma \ref{45}, $B|_{\Omega_1} $ is gluable.
 In particular, $z_1$ and $z_2$ are   $B$-glued on $\Omega_1$.
This means that there is a finite sequence $z_1=\lambda_1,
\lambda_2,\cdots,\lambda_N=z_2$
in  $\Omega_1$ such that  $\lambda_i$ and $\lambda_{i+1}$
are directly $B$-glued in $\Omega_1$
for each $i~ (1\leq    i\leq     N).$ Then by the conclusion in Case I,
 we have
$\rho_{w_0,\lambda_i} \sim  \rho_{w_0,\lambda_{i+1}}  $
for $1\leq   i\leq     N-1$, which gives
$\rho_{w_0,\lambda_1} \sim \rho_{w_0,\lambda_N}  $.
 That is, $\rho_{w_0,z_1} \sim \rho_{w_0,z_2} $, as desired.

The proof of Lemma \ref{component} is complete.
\end{proof}
As a special case of Lemma \ref{componenthigh}, we get the following consequence.
\begin{cor}\label{component}
Let $B$ be a finite Blaschke product that  is injective on
its critical value set $\mathcal{E}_B$.
For $0<t<1$,  suppose
$\Omega_1$ and  $\Omega_2$ are distinct components
 of $\{|B|<t\}$. Then for any
$z_1,z_2\in \Omega_1$ and $w_0 \in \Omega_2$   such that
$$B(z_1)=B(z_2)=B(w_0)\not\in  \mathcal{E}_B ,$$ we have $\rho_{w_0,z_1}\sim \rho_{w_0,z_2} .$
\end{cor}

\section{Proof of our main result} \label{sect3}

\subsection{Proof of Theorem \ref{submain}}

For the reader's convenience, we restate
 Theorem \ref{submain} as below.

\begin{thm} Suppose $B$ is a finite Blaschke product with
  at least two distinct critical points and
   $B|_{Z(B')}$ is injective.
 Then  $\mathcal{S}_B$ has exactly two components. In this case,
  the multiplication operator $M_B$ on $L_a^2(\mathbb{D})$
  has exactly two minimal reducing subspaces.
 \end{thm}

\subsubsection{Some observations}

 For $0 < s < 1$, let $n_s$ be the number of components
of $\{ z\in \mathbb{D}: |B(z)| < s \}$, and let $N_s$ be the number
 of components of $\{ z\in \mathbb{D}: |B(z)| \leq    s \}$.
 Using Proposition \ref{Jordancurve}, it is direct to see that
$N_s$ also equals the number of components of
$\{z\in \mathbb{D}: |B(z)| = s \}$.

\begin{rem} \label{rem33}
We have the following properties for $n_s$ and
 $N_s $.
\begin{itemize}\item[\emph{(1)}]
  The integer $n_s$ is a decreasing
function of $s$ in $(0,1]$,
and $n_s$ is discontinuous at $s=t$ if and only if
$\{ |B| = t \}$ contains at least one critical value. Moreover,
$n_s$ is left-continuous in $s$ on $(0,1]$.
For a detailed proof of these facts,
we refer the reader to (a) and (i) of \cite[Remark 3.9]{GHLiLuo}.
\item[\emph{(2)}]
 If $n_s$ is continuous at the point $t\in (0,1)$, then
  $N_t= n_t $.

To see this, let $\Omega_{0,t},
\Omega_{1,t}, \cdots, \Omega_{k,t}$ be all
components of $\{ |B| < t\}$.
We first prove that
$\partial \Omega_{i,t}
\cap \partial \Omega_{j,t} = \emptyset$
for all $i \neq j$. Otherwise, by Lemma \ref{component3},
we have for some $i\neq j,$
 $\partial \Omega_{i,t}
\cap \partial \Omega_{j,t}$ contains at least one critical point.
By (1),
$n_s$ is not continuous at some point $s^*\in (0,1)$ if and only if
  $ \{|B|=s^*\}$ contains at least one critical value of $B$.
Setting $s^*=t,$ it follows that
 $n_s$ is discontinuous at $t$, which is a contradiction.

Since $\partial \Omega_{i,t}
\cap \partial \Omega_{j,t} = \emptyset$
for all $i \neq j$, it follows that
 $\overline{\Omega}_{i,t} \cap \overline{\Omega}_{j,t}
 = \emptyset$ for all $i \neq j$. Moreover,
  each $\Omega_{i,t} $
 is simply connected  and bounded by a Jordan curve, and
\[
\{ |B| \leq    t \} = \overline{\{|B|<t\}}  =
 \bigcup_{i=0}^k \overline{\Omega}_{i,t}.
\]
Hence $N_t = n_t$.

\item[\emph{(3)}]
We claim that if $n_s$ is discontinuous at $t$, then
$N_t = \lim\limits_{s \to t^+} n_s < n_t$.

Having established this claim, recall that $N_t$ denotes the
 number of connected components of the set $\{|B| \leq t\}$.
  Consequently, by  (2), at this
  discontinuity point $t$, there exists
  some $\delta>0$ such that $N_s = n_s$
  for all $t < s < t+\delta$. For such $s$,
   every connected component of $\{|B| \leq t\}$ is contained in a nested family of connected components of $\{|B| < s\}$, and the total number of components on both sides is identical.

%

To see $N_t = \lim\limits_{s \to t^+} n_s < n_t$, let $\Omega_{0,t}, \cdots, \Omega_{k,t}$ be all components
of $\{ |B| < t \}$. Then $n_t = k + 1 \geq  2$.
Since $n_s$ is discontinuous at $t$, by  (1), $\{ |B| = t \}$ contains at least one critical value.
At least two of the boundaries $\partial \Omega_{0,t}, \cdots,
\partial \Omega_{k,t}$ share a common point on $\{ |B| = t\}$.
   Let $W = \bigcup_{i=0}^k \overline{\Omega}_{i,t}$,
   and set $N_t = l$, meaning that $W$
   has exactly $l$ connected components
   $V_1, \cdots, V_l$,
   with each $V_j$ containing at least one $\Omega_{i,t}$.
 We have $ l< k + 1 = n_t $.
  Since $n_s$ is a decreasing integer-valued function,
 it has only finitely many discontinuities. Hence there exists $\delta > 0$ such that for any $s' \in (t, t+\delta)$,
  $n_{s'}\equiv N $  for some integer $N$. For such $s'$, let
  $\Omega_{0,s'},\cdots, \Omega_{N-1,s'}$ be all components of
  $\{ |B| < s' \}$.

   Since  every $\Omega_{i,t}$ is contained in some $\Omega_{j,s'}$,
   and every $\Omega_{j,s'}$ contains at least one $\Omega_{i,t}$,
   it is easy to see  $l\geq   N. $  Besides,
    when $s'$ is greater than and sufficiently close to $t$,
     every component $\Omega_{j,s'}$ of $\{ |B| < s' \}$
     contains at most one component of
     $W = \bigcup\limits_{i=0}^k \Omega_{i,s'}$.
Therefore, $N\geq   l$, forcing $N=l.$ Thus, we have
 $ N_t =N= \lim\limits_{s \to t^+} n_s=l<k+1=n_t$.
  This proves the claim.
 \end{itemize}
\end{rem}

\subsubsection{Proof of Theorem \ref{submain}}

Throughout this subsection,  let $B$ denote a fixed finite
 Blaschke product.
  For each $s \in (0,1)$,
 let $\{\Omega_{j,s}\}_{j=1}^{n_s}$ be the finitely many
 components of the domain $\{z\in \mathbb{D}:|B(z)| < s\}$ (note
 that the number $n_s$ of these components
 decreases in $s$). Then
 $$\bigcup_{j=1}^{n_s}\partial \Omega_{j,s}= \{z\in \mathbb{D}:|B(z)| = s\}.$$
For this  fixed  $B$,
  we may assume without loss of generality that one of
  its critical value is zero. Otherwise,
   let $a$ be a critical value of $B$,
   and we replace $B$ by $\varphi_a(B) $, where
   $\varphi_a(z)=\frac{a-z}{1-\overline{a}z}, ~z \in \mathbb{D}$.
   Let $$\Lambda = \{ |B(z)| : B'(z) = 0, z \in \mathbb{D} \}$$ denote
the finite set of absolute critical values. Then $\Lambda$ can be
ordered as $0 = t_0 < t_1 < \cdots < t_N $.
 Consider 
  the intervals
 $$(t_0,t_1], (t_1,t_2],\cdots, (t_{N-1}, t_N],
 (t_N,t_{N+1}]=(t_N,1],$$
 where we define $t_{N+1}=1.$
 \begin{prop}\label{eA} With the above notations, let $0\leq     i\leq    N.$
 For $s'<s''$ and $s',s''\in (t_i,t_{i+1}]$,
  if $\Omega_{0,s'} $ and $\Omega_{0,s''}$
   are two components of $\{|B|<s'\}$ and $\{|B|<s''\}$ respectively,
    and $\Omega_{0,s'} \subseteq \Omega_{0,s''}$, then
    these two components share the same multiplicity. In this case, if
     $\Omega_{0,s'}$ is purely gluable,
    then $\Omega_{0,s''}$ is also purely gluable.
   \end{prop}
   \begin{proof}

We begin with the following  claim.
\vskip2mm
 \noindent \textbf{Claim.} \emph{For $0 < s' < s'' < 1$, suppose $\Omega_{0,s'}$ and
  $\Omega_{0,s''}$ share the same multiplicity and
   $\Omega_{0,s'} \subseteq \Omega_{0,s''}$.
    If $\Omega_{0,s'}$ is purely gluable,
    then $\Omega_{0,s''}$ is also purely gluable.}
    \vskip2mm

    \noindent    This property will be used when both $s'$ and $s''$ lie
 within $(t_i, t_{i+1}]$  for $0\leq i \leq N.$

 To see this, let $\phi = B|{\Omega_{0,s'}}$ and
  $\psi = B|{\Omega_{0,s''}}$. Since $\phi$ is essentially
   a finite Blaschke product,  by Propositions \ref{sproper} and \ref{45} we have that  $\phi$ is purely gluable
   on $\Omega_{0,s'}$. The set $\mathcal{S}_\phi=
 \{(z,w)\in \Omega_{0,s'}^2:
       B(z)=B(w) \not\in \mathcal{E}_B\} $ consists
   of exactly two components. In the language of local inverses,
   if $z_1, \cdots, z_k$ map to the same point under $B$
   and $k$ denotes the multiplicity of $\Omega_{0,s'}$,
    then we have
    $$ \rho_{z_1,z_2} \sim \rho_{z_1,z_i}, \quad 2 \leq   i \leq    k. $$
    This means the local inverse $\rho_{z_1,z_2}$
    can be analytically extended to $\rho_{z_1,z_i}$
    along a curve in $\Omega_{0,s'}$
    (and consequently in $\Omega_{0,s''}$).
    Given that $k$ also matches the multiplicity of
     $\Omega_{0,s''}$, the points $z_1, \cdots, z_k$
     are the only points in $B^{-1}(B(z_1)) \cap \Omega_{0,s''}$.
      We thereby conclude that $$\mathcal{S}_\psi=\{(z,w)\in \Omega_{0,s''}^2:
       B(z)=B(w) \not\in \mathcal{E}_B\} $$
      consists of exactly two components, establishing
      that $B$ remains purely gluable
 on $\Omega_{0,s''}$.   The proof of the claim is complete.

By Remark \ref{rem33}(1), for all $s\in (t_i,t_{i+1}]$, the quantity $n_s$ is equal to a fixed positive integer $K$. By left-continuity, this extends
to $n_s = K$ for all $s\in (t_i, t_{i+1}]$.
For any two values $s', s''$ with $t_i < s' < s'' \leq    t_{i+1}$,
the set $\{|B| < s'\}$ admits exactly $K$
connected components,
$\Omega_{0,s'}, \Omega_{1,s'}, \cdots, \Omega_{K-1,s'}$.
Similarly, the set $\{|B| < s''\}$ also contains exactly
$K$  components, which we label as $\Omega_{0,s''},
\Omega_{1,s''}, \cdots, \Omega_{K-1,s''}$.
Since every  component of $\{|B| < s''\}$ contains
 at least a   component of $\{|B| < s'\}$,
 and every   component of $\{|B| < s'\}$ is
 contained in exactly one  component of $\{|B| < s''\}$,
 there exists a permutation of the indices such that
 the inclusion relation
\begin{equation} \label{one}
\Omega_{i,s'} \subseteq \Omega_{i,s''}, \quad 0
 \leq    i \leq    K-1 \end{equation}
 holds.
Furthermore, the order of the Blaschke product $B$
equals the sum of the multiplicities of the components
 $\{\Omega_{i,s'}\}_{0\leq    i\leq     K-1}$, and also
 coincides with the sum of the multiplicities of the
 components $\{\Omega_{i,s''}\}_{0\leq    i\leq    K-1}$.
This, together with (\ref{one}),  forces the multiplicity of $\Omega_{i,s'}$
to be identical to that of $\Omega_{i,s''}$
for every $0\leq    i\leq     K-1$, and
 in particular for $i=0$.
Combining this result with the previously established
 claim and the fact that $\Omega_{0,s'}$ is purely gluable,
 we immediately conclude that $\Omega_{0,s''}$ is also purely gluable.
\end{proof}
\textbf{Outline of the proof of Theorem \ref{submain}. } 
Let $B$ be a finite Blaschke product with at least two distinct critical points, such that $B|_{Z(B')}$ is injective. We prove that for any $0<s\leq    1$, every component of $\{|B|<s\}$ containing at least two critical points is purely gluable.

We proceed by induction on the intervals $(t_j, t_{j+1}]$:
\begin{enumerate}
    \item \textbf{Base case:} The claim holds for $s\in(t_0,t_1]$.
    \item \textbf{Inductive step:} Assume the claim holds for $s\in(t_k,t_{k+1}]$; we extend it to $s\in(t_{k+1},t_{k+2}]$ in two parts:
    \begin{itemize}
        \item First, we show there exists
        $\varepsilon>0$ such that the claim holds
        for all $s \in(t_{k+1},t_{k+1}+\varepsilon)$ (the core technical step of the proof).
        \item  Second, we extend the result to the full interval $(t_{k+1}, t_{k+2}]$.
 By Remark \ref{rem33}(1), we have $n_{s'} = n_{s''}$ for any $s', s'' \in (t_{k+1}, t_{k+2})$.
Furthermore, due to the left continuity of $n_s$, this equality extends to $s', s'' \in (t_{k+1}, t_{k+2}]$.
Observe that for such $s', s''$ with $s' < s''$,
every component of $\{ |B| < s' \}$ is contained
in a unique component of $\{ |B| < s'' \}$,
and conversely, every component of $\{ |B| < s'' \}$
 contains exactly one component of $\{ |B| < s' \}$.
This implies that the components of $\{ |B| < s \}$
 are nested as $s$ increases over $(t_{k+1}, t_{k+2}]$.
Consequently, by Proposition \ref{eA}, for each $s
 \in (t_{k+1}, t_{k+2}]$, every component of
  $\{ |B| < s \}$ containing at least
  two distinct critical points is purely gluable.
      \end{itemize}
\end{enumerate}

After finitely many steps, the claim holds for $s\in(t_N,t_{N+1}]$. Since $t_{N+1}=1$, the only component of $\{|B|<1\}$ is $\mathbb{D}$, so $\mathbb{D}$ is purely gluable. By definition, $\mathcal{S}_B$ has exactly two components, which completes the proof.
\hfill $\square$

To complete the proof of Theorem \ref{submain}, it therefore suffices to
 establish the inductive claim
  for $s\in(t_{k+1},t_{k+1}+\varepsilon)$,
    rather than verifying  it over  the entire interval.

\vskip2mm

\begin{proof} [\textbf{\emph{Proof of Theorem \ref{submain}.}}]
Let $B$ be a finite Blaschke product with at least two distinct critical points, and
suppose that $B|_{Z(B')}$ is injective.
     By Proposition \ref{thinfin},
     the proof is complete
     once we establish that
      for $0<s\leq  1,$
     each component  $V$  of $\{|B|<s\}$
     that contains at least two distinct critical points
        is purely gluable; that is,
 $$\{(z,w)\in V^2: B(z)=B(w)\not\in  \mathcal{E}_B \}$$
contains exactly two components, the identity component
$\mathcal{S}_0$ and its complement. When this is done,
for $s=1$, we have $\mathbb{D}=\{|B|<1\}$, and then
 $$\{(z,w)\in  \mathbb{D}^2: B(z)=B(w)\not\in  \mathcal{E}_B \}$$
contains exactly two components. That is,
  $\mathcal{S}_B$ has exactly two components.
  By Proposition \ref{thinfin}, $M_B$ on $L_a^2(\mathbb{D})$
   has precisely two minimal reducing subspaces.

We proceed via induction to verify this statement.
  Note that $t_0=0$, $t_{N+1}=1$, and
  $(0,1] = \bigcup\limits_{j=0}^N (t_j, t_{j+1}]$.
  For $s\in (t_0,t_1]=(0,t_1]$, we will shortly prove that each component
of $\{|B|<s\} $  has at most one critical point.
 Under this condition,  we will soon see that
the desired statement immediately holds  for $s\in (t_0,t_1]$.
In fact, let $\lambda_1,\cdots, \lambda_N$
be the zeros of $B.$ There exists a sufficiently small $\delta>0$
such that each component of $\{|B|<\delta\}$
lies in an $\varepsilon$-neighborhood of some $\lambda_i$ for
$ 1\leq    i\leq     N .$
By the discreteness of the critical set of $B$,
each such domain   contains
 at most one possible critical point, whose  only candidate is the zero
 $\lambda_i$. When $\delta$ is replaced by a number
 $s'\in (t_0,t_1],$ note that the components
 are nested as $s'$ increases.

 To finish the proof of this case,
 we claim that each component of
  \(\{|B|<s'\}\) contains at most
  one critical point. The desired
   conclusion then follows. Suppose, to reach a contradiction,
    that there exist some \(s'\in(0,t_1]\) and a component \(\Omega\) containing at least two critical points. Since every critical point \(p\) in \(\Omega\) satisfies
\[
|B(p)| < s' \leq  t_1,
\]
it follows from the definition
 of \(t_i\) that \(B(p)=0\). But as
  the components of \(\{|B|<t'\}\) are nested as \(t'\) increases, each such component can contain at most one critical point. This contradiction completes the proof of this case.

Assume each component  of $\{|B|<s\}$ that contains
 at least two distinct critical points is purely gluable for
$s\in(t_j,t_{j+1}]$ with $0\leq    j\leq     k$. We aim to prove this holds
for $s\in(t_{k+1},t_{k+2}]$. By the Proposition \ref{eA}, it suffices to find
$\varepsilon\in(0,t_{k+2}-t_{k+1})$ such that each component of $\{|B|<s\}$
containing
 at least two distinct critical points is purely
  gluable for any $s\in(t_{k+1},t_{k+1}+\varepsilon)$.

 \textbf{Notation.}  $t$ and $t^+.$
 In the subsequent induction step, we rewrite
 $t = t_j$ for each index  $j$.
 Note that $t_j$ denotes the time $t$ at which at least two
 boundaries $\partial \Omega_{i,t}$ intersect (see Remark \ref{rem33}(1)).

 Let us see what happens to the components of
 $\{|B|\leq   s\}$ when $s=t$ is replaced by
$s=t+\varepsilon.$
Assume $\ell \geq 2$. After relabeling the components $\Omega_{i,t}$, we may suppose that the union $\bigcup\limits_{i=1}^\ell \overline{\Omega}_{i,t}$ is connected, and that no other $\overline{\Omega}_{j,t}$
intersects this union.
 Let $\varepsilon > 0$ be a prescribed parameter,
and define $\Omega$ as the connected component of the set
$\{|B| < t + \varepsilon\}$ that contains $\bigcup\limits_{i=1}^\ell
\overline{\Omega}_{i,t}$. By choosing $\varepsilon$ sufficiently small,
we ensure that $\Omega$ lies in a  small neighborhood
of $\bigcup\limits_{i=1}^\ell \overline{\Omega}_{i,t}$, which contains
no more critical points than the original union $\bigcup\limits_{i=1}^\ell
\overline{\Omega}_{i,t}$. Since the full union  $\bigcup\limits_{j}
\overline{\Omega}_{j,t}=\{|B|\leq    t\}$ consists of only finitely many connected
components, there exists a common constant $\varepsilon > 0$ such
that this property holds for every corresponding component $\Omega$.
\emph{We let $t^+$ denote an arbitrary number in $(t,t+\varepsilon )$.
 In the sequel, $t^+$ consistently keeps this meaning.}

 It suffices to consider every
  component of \(\{|B|<t^+\}\) containing at least two
  distinct critical points.
  We shall prove by induction on $s$ that,
 for a fixed $t=t_k$, such components of \(\{|B|<t^+\}\)
  are purely gluable.

 \textbf{Step I}.  \textbf{Flower model.} Remind $t=t_k.$
   \emph{In this case, some of the boundaries $\partial \Omega_{j,t}$ share a common point $*$, with any two of them intersecting only at $*$.
For all remaining indices $i$, $\partial \Omega_{i,t}$ is disjoint from every such $\partial \Omega_{j,t}$.}
Since the union of the $\partial \Omega_{j,t}$ passing through $*$ forms a petal-like shape, this set of boundaries sharing the common point $*$ is defined as the \emph{flower model}.
By hypothesis, the union of the corresponding closed regions $\overline{\Omega}_{j,t}$ contains at least two distinct critical points.


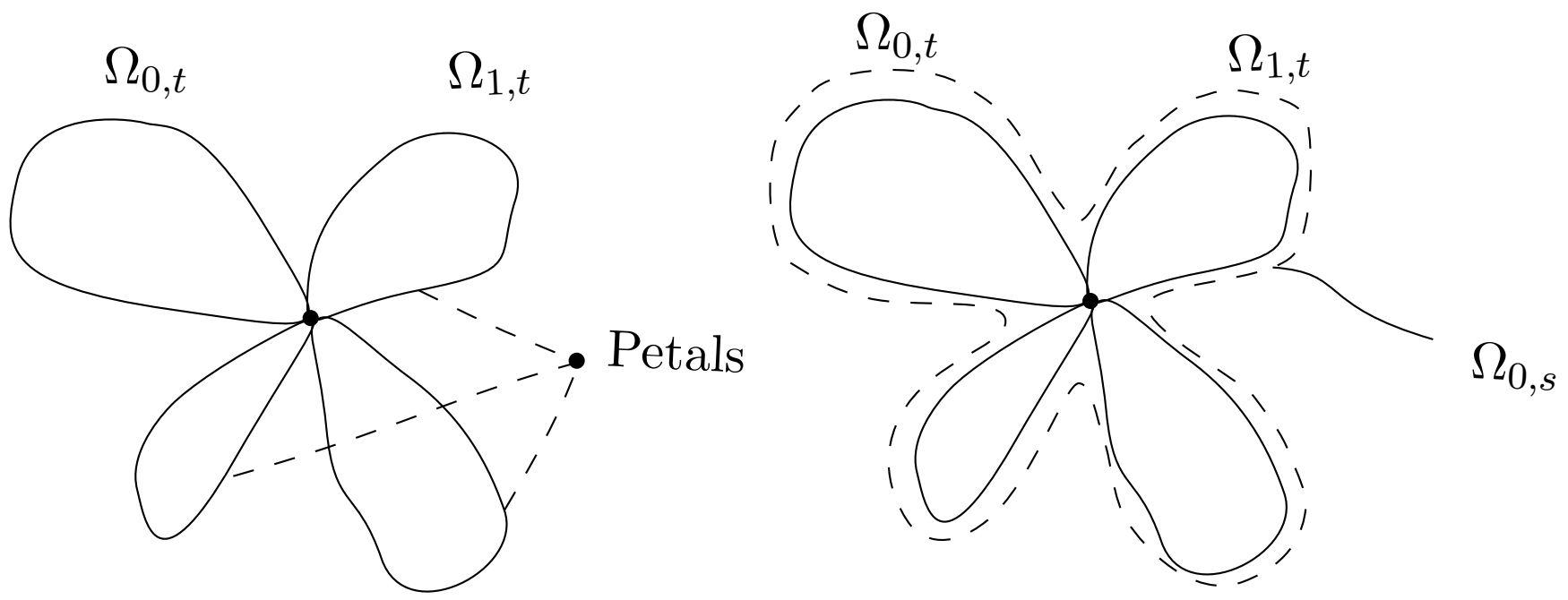
\begin{figure}[htbp]
\centering
\begin{tikzpicture}[x=0.6pt,y=0.6pt,yscale=-1,xscale=1]
\draw   (374.45,177.98) .. controls (382.78,181.95) and (392.4,174.69) .. (416.22,214.13) .. controls (440.05,253.56) and (441.22,251.62) .. (382.67,243.12) .. controls (324.12,234.63) and (323.57,221.8) .. (329.51,197.32) .. controls (335.46,172.84) and (366.12,174.02) .. (374.45,177.98) -- cycle ;
\draw   (469.96,236.49) .. controls (432.71,243.82) and (430.04,255.92) .. (431.16,236.18) .. controls (432.27,216.45) and (442.1,202.73) .. (460.22,188.23) .. controls (478.34,173.73) and (510.61,183.84) .. (503.8,204.96) .. controls (496.99,226.08) and (507.2,229.15) .. (469.96,236.49) -- cycle ;
\draw   (383.25,276.86) .. controls (389.28,270.66) and (403.91,260.57) .. (423.37,250.46) .. controls (442.83,240.35) and (431.94,250.08) .. (404.94,296.82) .. controls (377.95,343.57) and (374.08,317.13) .. (371.38,306.04) .. controls (368.67,294.96) and (377.23,283.06) .. (383.25,276.86) -- cycle ;
\draw  [line width=3] [line join = round][line cap = round] (432.17,246.22) .. controls (432.17,246.22) and (432.17,246.22) .. (432.17,246.22) ;
\draw   (437.73,285.22) .. controls (434.88,256.79) and (429.24,249.14) .. (435.17,246.35) .. controls (441.1,243.56) and (451.18,255.49) .. (466.66,266.94) .. controls (482.14,278.39) and (492.92,292.79) .. (499.99,313.57) .. controls (507.07,334.36) and (465.7,355.85) .. (456.99,330.53) .. controls (448.27,305.21) and (440.57,313.64) .. (437.73,285.22) -- cycle ;
\draw  [dash pattern={on 4.5pt off 4.5pt}] (334.64,172.96) .. controls (341.64,164.63) and (372,162.54) .. (385.23,169.49) .. controls (398.47,176.44) and (404.14,181.97) .. (409.36,191.72) .. controls (414.59,201.47) and (423.37,217.42) .. (428.04,218.12) .. controls (432.71,218.81) and (440.5,196.58) .. (447.5,191.02) .. controls (454.51,185.47) and (462.77,176.79) .. (469.3,175.05) .. controls (475.83,173.3) and (478.14,171.36) .. (487.2,172.96) .. controls (496.26,174.57) and (506.66,176.44) .. (508.22,184.08) .. controls (509.77,191.72) and (510.55,216.73) .. (504.32,228.54) .. controls (498.1,240.35) and (452.17,238.96) .. (452.17,247.99) .. controls (452.17,257.02) and (484.87,273.7) .. (489.54,279.95) .. controls (494.21,286.2) and (499.24,292.59) .. (501.99,300.09) .. controls (504.73,307.6) and (508.22,313.29) .. (507.44,319.55) .. controls (506.66,325.8) and (503.57,337.84) .. (485.64,344.56) .. controls (467.71,351.27) and (450.62,332.05) .. (445.17,323.71) .. controls (439.72,315.38) and (439.72,302.87) .. (436.92,295.09) .. controls (434.12,287.31) and (432.84,275.08) .. (428.04,275.08) .. controls (423.25,275.08) and (411.37,310.62) .. (396.13,323.02) .. controls (380.89,335.42) and (371.68,328.15) .. (368.89,324.41) .. controls (366.1,320.67) and (355.8,304.7) .. (365.77,286.2) .. controls (375.75,267.7) and (400.8,262.58) .. (402.36,254.24) .. controls (403.91,245.91) and (383.44,247.73) .. (364.22,246.6) .. controls (344.99,245.47) and (333.57,236.7) .. (326.86,232.71) .. controls (320.14,228.71) and (319.07,207.7) .. (320.63,195.89) .. controls (322.19,184.08) and (327.63,181.3) .. (334.64,172.96) -- cycle ;
\draw   (101.45,183.98) .. controls (108.78,185.95) and (119.4,180.69) .. (143.22,220.13) .. controls (167.05,259.56) and (168.22,257.62) .. (109.67,249.12) .. controls (51.12,240.63) and (50.57,227.8) .. (56.51,203.32) .. controls (62.46,178.84) and (94.12,182.02) .. (101.45,183.98) -- cycle ;
\draw   (196.96,242.49) .. controls (159.71,249.82) and (157.04,261.92) .. (158.16,242.18) .. controls (159.27,222.45) and (169.1,208.73) .. (187.22,194.23) .. controls (205.34,179.73) and (237.61,189.84) .. (230.8,210.96) .. controls (223.99,232.08) and (234.2,235.15) .. (196.96,242.49) -- cycle ;
\draw   (110.25,282.86) .. controls (116.28,276.66) and (130.91,266.57) .. (150.37,256.46) .. controls (169.83,246.35) and (158.94,256.08) .. (131.94,302.82) .. controls (104.95,349.57) and (101.08,323.13) .. (98.38,312.04) .. controls (95.67,300.96) and (104.23,289.06) .. (110.25,282.86) -- cycle ;
\draw  [line width=3] [line join = round][line cap = round] (159.17,252.22) .. controls (159.17,252.22) and (159.17,252.22) .. (159.17,252.22) ;
\draw   (164.73,291.22) .. controls (161.88,262.79) and (156.24,255.14) .. (162.17,252.35) .. controls (168.1,249.56) and (178.18,261.49) .. (193.66,272.94) .. controls (209.14,284.39) and (219.92,298.79) .. (226.99,319.57) .. controls (234.07,340.36) and (192.7,361.85) .. (183.99,336.53) .. controls (175.27,311.21) and (167.57,319.64) .. (164.73,291.22) -- cycle ;
\draw  [dash pattern={on 4.5pt off 4.5pt}]  (253.11,267.69) .. controls (210.43,278.83) and (189.43,291.83) .. (128.11,308.69) ;
\draw  [dash pattern={on 4.5pt off 4.5pt}]  (196.96,242.49) .. controls (229.33,258.05) and (215.33,252.05) .. (253.11,267.69) ;
\draw  [dash pattern={on 4.5pt off 4.5pt}]  (226.99,319.57) .. controls (243.43,291.83) and (247.33,283.05) .. (253.11,267.69) ;
\draw    (495.96,234.49) .. controls (522.67,236.35) and (510.44,247.51) .. (552.11,259.69) ;
\draw  [line width=3] [line join = round][line cap = round] (252.31,267.16) .. controls (252.31,267.16) and (252.31,267.16) .. (252.31,267.16) ;

\node at (432.17,246.22) [circle, fill=black, inner sep=1.2pt] {}; %
\node at (159.17,252.22) [circle, fill=black, inner sep=1.2pt] {}; %
\node at (252.31,267.16) [circle, fill=black, inner sep=1.2pt] {}; 

\draw (348.14,143.9) node [anchor=north west][inner sep=0.75pt]  [rotate=-0.66] [align=left] {$\Omega_{0,t}$};
\draw (477.95,151.95) node [anchor=north west][inner sep=0.75pt]  [rotate=-357.63] [align=left] {$\Omega_{1,t}$};
\draw (85.14,155.9) node [anchor=north west][inner sep=0.75pt]  [rotate=-0.66] [align=left] {$\Omega_{0,t}$};
\draw (204.95,157.95) node [anchor=north west][inner sep=0.75pt]  [rotate=-357.63] [align=left] {$\Omega_{1,t}$};
\draw (261.81,254.61) node [anchor=north west][inner sep=0.75pt]  [rotate=-2.8] [align=left] {Petals};
\draw (564.14,259.02) node [anchor=north west][inner sep=0.75pt]  [rotate=-3.11] [align=left] {$\Omega_{0,s}$};
\end{tikzpicture}
\caption{Petals associated with $\Omega_{0,t}$}
\label{Fig-7}
\end{figure}
After reordered, assume that  $\partial \Omega_{0,t}$,
$\partial  \Omega_{1,t},\cdots$, $\partial \Omega_{k,t}$ are petals
 of the flower. We will soon show that
 there exists at least one $\Omega_{j,t}$ $(1 \leq    j \leq    k)$
 with multiplicity at least two. To prove this,
 we proceed by contradiction and assume
  conversely that each $\Omega_{j,t}$ ($1 \leq j \leq k$)
   has multiplicity one, meaning that
   $B$ is univalent on each such $\Omega_{j,t}$.
    Consequently, by Corollary \ref{e2} the only critical
  point on $\bigcup\limits_{i=1}^k \overline{\Omega}_{i,t}$ would be $*$.
Let $\delta > 0$ be a prescribed parameter and let $\Omega$
denote the component of $\{|B| < t + \delta\}$ containing
$\bigcup\limits_{i=1}^k \overline{\Omega}_{i,t}$. By selecting $\delta$
sufficiently small, we ensure that $\Omega$ lies in
a small neighborhood $\mathcal{N}$ of
$\bigcup\limits_{i=1}^k \overline{\Omega}_{i,t}$.
Thus either $\Omega$ or this neighborhood
$\mathcal{N}$ contains  no more critical points than  $\bigcup\limits_{i=1}^k \overline{\Omega}_{i,t}$
does. Consequently, $\mathcal{N}$ contains exactly one critical point,  $*$.
However, by the induction hypothesis, $\Omega$ must possess
at least two distinct critical points.
 This contradiction completes the proof.

Without loss of generality, we
     assume that $\partial \Omega_{0,t}$
   is of multiplicity greater than $  2 $.
 There are two situations to distinguish:
  either $\Omega_{0,t}$  contains at least two
   critical points or not.  In either case, we will
  assume $s$ is a bit larger than $t$ so that
  $\Omega_{0,s} $ denote the component containing
   $\Omega_{0,t}$ and the petal associated with it,
  but $\Omega_{0,s} $ has no intersection with other
  $  \Omega_{j,t}$ (see Figure \ref{Fig-7}).
   Following our earlier convention
  we write $t^+$ for this $s$. As follows,
  we now show that $\Omega_{0,t^+} $ is purely gluable.

\textbf{Case I.1}   $\Omega_{0,t}$ 
\emph{contains at least two distinct
critical points.}

 By induction on $s$  with $s=t$,
we conclude that $\Omega_{0,t}$ is purely gluable since it contains at least
two distinct critical points.

 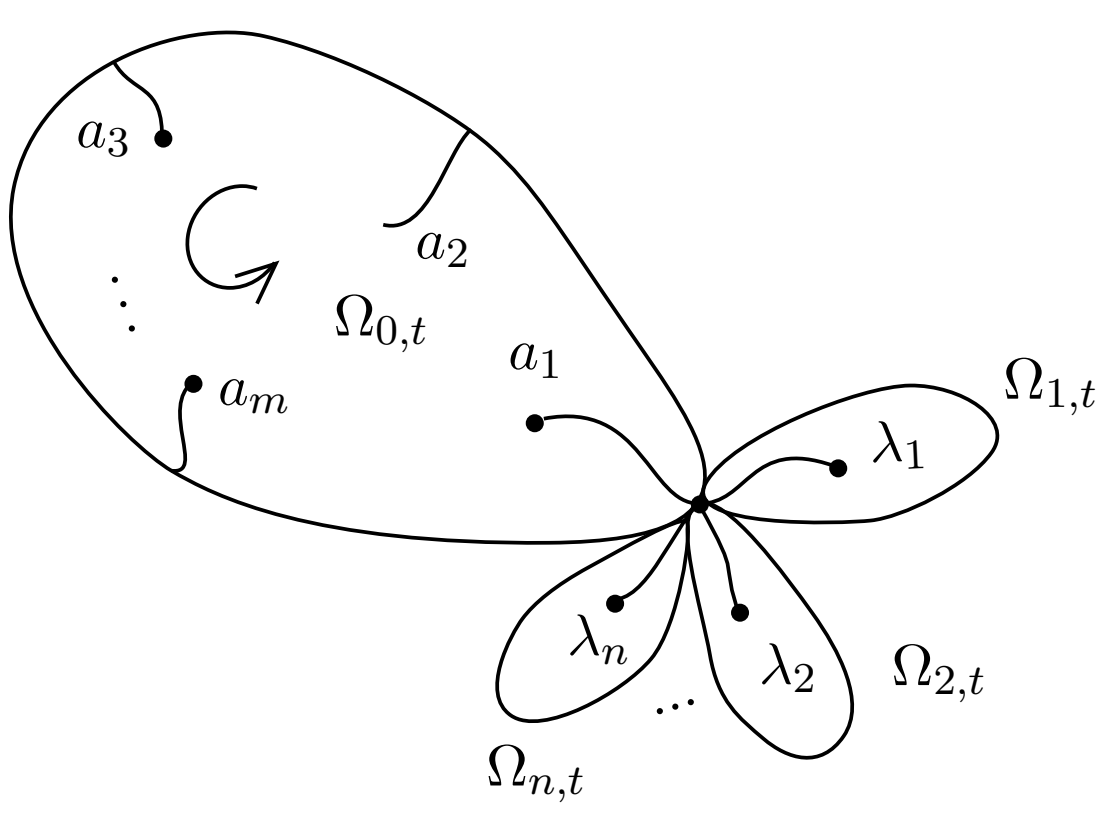
\begin{figure}[htbp]
\centering
\tikzset{every picture/.style={line width=0.70pt}} 

\begin{tikzpicture}[x=0.70pt,y=0.70pt,yscale=-1,xscale=1]

\draw   (242.71,125.41) .. controls (264.31,130.27) and (291.73,144.47) .. (303.73,155.21) .. controls (315.73,165.96) and (321.42,176.23) .. (340.62,204.12) .. controls (359.81,232.02) and (383.81,261.13) .. (319.51,261.64) .. controls (255.21,262.16) and (229.55,249.11) .. (217.92,242.26) .. controls (206.3,235.4) and (169.01,198.88) .. (175.52,166.3) .. controls (182.04,133.73) and (221.12,120.56) .. (242.71,125.41) -- cycle ;
\draw   (405.99,255.66) .. controls (397.75,256.53) and (354.76,257.82) .. (361.64,245.52) .. controls (368.52,233.22) and (402.94,219.52) .. (416.51,219.36) .. controls (430.07,219.19) and (444.89,228.34) .. (438.45,237.37) .. controls (432,246.41) and (414.24,254.8) .. (405.99,255.66) -- cycle ;
\draw   (312.3,282.08) .. controls (315.38,277.92) and (321.21,273.07) .. (330.12,268.22) .. controls (339.04,263.36) and (354.93,254.11) .. (356.49,256.05) .. controls (358.04,257.99) and (354.46,283.81) .. (346.92,292.82) .. controls (339.38,301.83) and (317.78,313.96) .. (309.21,308.07) .. controls (300.64,302.17) and (309.21,286.23) .. (312.3,282.08) -- cycle ;
\draw  [line width=3] [line join = round][line cap = round] (397.21,241.63) .. controls (397.21,241.63) and (397.21,241.63) .. (397.21,241.63) ;
\draw   (362.73,291.84) .. controls (360.57,279.35) and (352.26,253.74) .. (359.69,251.13) .. controls (367.11,248.53) and (383.26,270.81) .. (390.32,280.64) .. controls (397.38,290.47) and (404.85,304.6) .. (398.68,313.61) .. controls (392.51,322.62) and (383.6,319.85) .. (377.09,314.3) .. controls (370.57,308.76) and (364.88,304.34) .. (362.73,291.84) -- cycle ;
\draw  [line width=3] [line join = round][line cap = round] (370.81,280.44) .. controls (370.81,280.44) and (370.81,280.44) .. (370.81,280.44) ;
\draw  [line width=3] [line join = round][line cap = round] (337.22,278.02) .. controls (337.22,278.02) and (337.22,278.02) .. (337.22,278.02) ;
\draw  [line width=3] [line join = round][line cap = round] (315.62,229.5) .. controls (315.62,229.5) and (315.62,229.5) .. (315.62,229.5) ;
\draw  [line width=3] [line join = round][line cap = round] (360.01,251.33) .. controls (360.01,251.33) and (360.01,251.33) .. (360.01,251.33) ;
\draw    (318.06,228.28) .. controls (344.82,223.13) and (346.12,250.94) .. (359.69,251.13) .. controls (373.25,251.32) and (373.25,233.13) .. (395.81,240.89) ;
\draw    (359.69,251.13) .. controls (352.85,257.39) and (345.65,278.01) .. (336.06,276.79) ;
\draw    (370.85,281.64) .. controls (364.85,264.66) and (372.05,273.15) .. (359.69,251.13) ;
\draw    (202.45,132.33) .. controls (207.6,141.38) and (215.89,137.73) .. (215.39,153.73) ;
\draw    (298.5,150.47) .. controls (291.66,156.72) and (287.19,179.13) .. (274.87,176.13) ;
\draw    (223.89,218.23) .. controls (214.39,225.23) and (227.52,243.47) .. (217.92,242.26) ;
\draw  [line width=3] [line join = round][line cap = round] (236.16,179.38) .. controls (236.16,179.38) and (236.16,179.38) .. (236.16,179.38) ;
\draw  [line width=3] [line join = round][line cap = round] (215.76,152.97) .. controls (215.76,152.97) and (215.76,152.97) .. (215.76,152.97) ;
\draw  [line width=3] [line join = round][line cap = round] (275.55,176.46) .. controls (275.55,176.46) and (275.55,176.46) .. (275.55,176.46) ;
\draw  [line width=3] [line join = round][line cap = round] (223.86,218.91) .. controls (223.86,218.91) and (223.86,218.91) .. (223.86,218.91) ;
\draw  [draw opacity=0] (245.68,186.31) .. controls (242.93,190.39) and (238.4,193.12) .. (233.61,193.11) .. controls (226.44,193.1) and (221.39,186.96) .. (222.33,179.4) .. controls (223.26,171.84) and (229.83,165.73) .. (237,165.74) .. controls (238.4,165.75) and (239.71,165.98) .. (240.92,166.42) -- (235.3,179.43) -- cycle ; \draw   (245.68,186.31) .. controls (242.93,190.39) and (238.4,193.12) .. (233.61,193.11) .. controls (226.44,193.1) and (221.39,186.96) .. (222.33,179.4) .. controls (223.26,171.84) and (229.83,165.73) .. (237,165.74) .. controls (238.4,165.75) and (239.71,165.98) .. (240.92,166.42) ;
\draw   (235.05,190.02) -- (246.05,186.56) -- (240.94,197.34) ;

\node at (397.21,241.63) [circle, fill=black, inner sep=1.2pt] {}; 
\node at (370.81,280.44) [circle, fill=black, inner sep=1.2pt] {}; 
\node at (337.22,278.02) [circle, fill=black, inner sep=1.2pt] {}; 
\node at (315.62,229.5) [circle, fill=black, inner sep=1.2pt] {}; 
\node at (360.01,251.33) [circle, fill=black, inner sep=1.2pt] {}; %
\node at (215.76,152.97) [circle, fill=black, inner sep=1.2pt] {}; 
\node at (223.86,218.91) [circle, fill=black, inner sep=1.2pt] {}; 

\draw (281.93,176.89) node [anchor=north west][inner sep=0.75pt]   [align=left] {$a_2$};
\draw (306.92,206.62) node [anchor=north west][inner sep=0.75pt]   [align=left] {$a_1$};
\draw (322.76,279.33) node [anchor=north west][inner sep=0.75pt]   [align=left] {$\lambda_n$};
\draw (374.35,287.09) node [anchor=north west][inner sep=0.75pt]   [align=left] {$\lambda_2$};
\draw (403.65,227.61) node [anchor=north west][inner sep=0.75pt]  [rotate=-357.3] [align=left] {$\lambda_1$};
\draw (228.89,216.23) node [anchor=north west][inner sep=0.75pt]   [align=left] {$a_m$};

\draw (361.84,305.75) node [anchor=north west][inner sep=0.75pt]  [rotate=-163.82] [align=left] {...};

\draw (190.74,146.97) node [anchor=north west][inner sep=0.75pt]   [align=left] {$a_3$};

\draw (259.94,193.64) node [anchor=north west][inner sep=0.75pt]   [align=left] {$\Omega_{0,t}$};

\draw (205.86,185.87) node [anchor=north west][inner sep=0.75pt]  [rotate=-70.82] [align=left] {$\cdots$};

\draw (300.94,314.64) node [anchor=north west][inner sep=0.75pt]   [align=left] {$\Omega_{n,t}$};

\draw (409.94,287.64) node [anchor=north west][inner sep=0.75pt]   [align=left] {$\Omega_{2,t}$};

\draw (439.94,210.64) node [anchor=north west][inner sep=0.75pt]   [align=left] {$\Omega_{1,t}$};
\end{tikzpicture}
\caption{Flower model}
\label{Fig-9}
\end{figure}

Let $w_0$ be a point in  $t\mathbb{D}\backslash  \mathcal{E}_B $.
 Let $a_1,\cdots, a_n$ be all preimages of $w_0$ in $\Omega_{0,t}$.
 For each fixed $\Omega_{j,t}\ (j>0)$, there exists an integer
 $k=k(j)$ so that
 $$B^{-1}(w_0)\cap \Omega_{j,t}  \ \ \ \ (\text{see~ Figure} ~\ref{Fig-9})$$
 contains $k$ elements, $b_1,\cdots,b_k$.
 By Corollary \ref{component},
 $\rho_{a_2,b_1}\sim \rho_{a_2,b_i}$ along a curve in
 $\Omega_{j,t^+}$ for each $i$.
   That is,   \begin{equation} \label{eqbew}
\rho_{a_2,z'}\sim \rho_{a_2,z''},\ \ \ \
   z',z''\in B^{-1}(w_0)\cap \Omega_{j,t}, \end{equation}
   and denote this relation by   $\rho_{a_2,\Omega_{j,t}}$;
    for a local inverse $\rho$ of $B$, we write
   $\rho\sim \rho_{a_2,\Omega_{j,t}}$ if there exists a point
    $z'\in B^{-1}(w_0)\cap \Omega_{j,t}$ such that
    $\rho \sim\rho_{a_2,z'}$.

    Remind that  a curve  $\gamma$ is
    non-self-intersecting  if $\gamma$ is injective, that is,
$\gamma(s')\neq \gamma(s'')$   for
$0\leq     s'<s'' \leq     1 $.
 Let $\gamma$ be a non-self-intersecting  curve
 connecting $w_0$ and $B(*) $,   and
 $$\gamma [0,1)\subseteq  t\mathbb{D}
  \backslash \mathcal{E}_B. $$
  Then by Lemma \ref{component3}, in each
 $\Omega_{j,t}~ (j\geq  0)$ there exists exactly one
 curve $\sigma_j $ with endpoint $*$ ($\sigma_j(1)=*$), so that
 $$
 B(\sigma_j(u))=\gamma(u), \ \ \ \  0\leq    u \leq     1, \ \ \ \
  j\geq  0.
 $$
 Write $\lambda_j =\sigma_j(0).$
 Note that $\sigma_0 \subseteq \Omega_{0,t}$
 and $\sigma_0$ connects $*$ to some point $a_i$.
 Without loss of generality,
 assume $\sigma_0(1)=a_1$.

 We next prove that
$$
\rho_{a_2,a_1} \sim \rho_{a_2,\Omega_{j,t}}, \qquad j\geq1.
$$
To this end, let $\sigma$ be a curve satisfying
$$
B(\sigma(u))=\gamma(u), \qquad \sigma(0)=a_2.
$$
We first claim that $B'(\sigma(1))\neq 0$. Suppose, to the contrary, that $B'(\sigma(1))=0$. By Lemma~\ref{component3}, we also have $B'(*)=0$. Since
$$
B(\sigma(u)) = B(\sigma_j(u)) = \gamma(u), \qquad 0\leq  u\leq  1,
$$
it follows that $B(\sigma(1)) = B(\sigma_j(1))$.
 But $\sigma$ ends at $\sigma(1)$, while
  $\sigma_j$ ends at $*$; hence the
   two critical points $\sigma(1)$ and
   $*$ are mapped to the same value under
    $B$. This contradicts the injectivity
    of $B$ on $Z(B')$, so we conclude
    that $B'(\sigma(1))\neq 0$.

Applying Lemma~\ref{gluel} now gives
\[
\rho_{a_2,a_1} \sim \rho_{a_2,\lambda_j}, \qquad j\geq1.
\]
Since $\lambda_j\in \Omega_{j,t}$, by (\ref{eqbew}) we have
 \(\rho_{a_2,\lambda_j} \sim \rho_{a_2,\Omega_{j,t}}.\)
 Combining this with the preceding equivalence yields
\begin{equation}\label{eq2}
\rho_{a_2,a_1} \sim \rho_{a_2,\Omega_{j,t}}, \qquad j\geq1,
\end{equation}
as desired.

 Finally, we will show that  $\Omega_{0,t^+} $ is
   purely gluable. Since $\Omega_{0,t} $ is
   purely gluable by induction, we have
    $ \rho_{a_2,a_1} \sim \rho_{a_2,a_j}~$ for
   $j\neq 2$. This, together with  \eqref{eq2}, yields
\[
\rho_{a_2,a_1} \sim \rho_{a_2,z'}
\]
for every \(z'\in B^{-1}(w_0)\cap\Omega_{0,t^+}\) with \(z'\ne a_2\), with the convention that \(\rho_{a_2,a_2}=\mathrm{id}\).
Thereby, it follows that on
\(\Omega_{0,t^+}\backslash B^{-1}(\mathcal{E}_B)\),
there are exactly two equivalence classes of local inverses of \(B\). By Proposition~\ref{newequi}, the Riemann surface
\[
\{(z,w)\in \Omega_{0,t^+}: B(z)=B(w)\notin \mathcal{E}_B\}
\]
has exactly two components: the component containing \((a_2,a_1)\) and the identity component \(\mathcal{S}_0\). Hence \(\Omega_{0,t^+}\) is purely gluable, which completes the proof in this case.

  \textbf{Case I.2.}   \textbf{Rotation case}.
 In this case, $\Omega_{0,t}$  contains exactly one critical point $p$, called \emph{the rotation case}.

 Since \(B|_{\Omega_{0,t}}\) is a proper map
  of multiplicity,  say $n$, Proposition~\ref{sproper} yields a conformal map \(F:\Omega\to\mathbb D\) and a finite Blaschke product \(\phi\) of order \(n\) such that
\[
B|_{\Omega_{0,t}} = t(\phi\circ F)|_{\Omega_{0,t}}.
\]
Since \(\Omega_{0,t}\) contains exactly
one critical point \(p\), \(\phi\)
also has exactly one critical point. There exist $a \in \mathbb{D}$ and a holomorphic automorphism $\varphi$ of $\mathbb{D}$ such that
$
\phi = \varphi(\phi_a^n).
$
  From the above factorization and the properties of
  \(\phi\), we obtain
  $n$ pairwise disjoint non-self-intersecting curves $\gamma_1, \cdots, \gamma_n$
and $n$ distinct points $a_1, \cdots, a_n$ in $\Omega_{0,t}\backslash
B^{-1}(\mathcal{E}_B)$ such that
$$
B(\gamma_1(u)) = \cdots = B(\gamma_n(u)),
  \ \ \  0\leq      u \leq    1
$$
with $\gamma_i(0) = a_i$ and $\gamma_i(1) = a_{i+1}$ for $1 \leq     i \leq     n$
(where $a_{n+1} \equiv a_1$).
Applying  Lemma \ref{continuation} yields the following
  equivalence
\begin{equation*} 
\rho_{\gamma_i(0), \gamma_j(0)} \sim \rho_{\gamma_i(1), \gamma_j(1)},\ \ \ 1\leq  i\leq     j\leq     n
\end{equation*}
For $i=1 $ and $j=2,\cdots,n$, we have
\begin{equation} \label{equij2}
\rho_{a_1,a_j}\sim \rho_{a_2,a_{j+1}}, \ \ \ 2\leq     j \leq     n.
\end{equation}
Here the case $j=n$ yields $\rho_{a_1,a_n}\sim \rho_{a_2,a_1}.$

As follows, we will prove that $ (a_2,a_i)$ (where $i\neq 2$) lie in a same
component of $\{(z,w)\in \Omega_{0,t^+}: B(z)=B(w)\not\in \mathcal{E}_B\}$.
For this, let $\lambda_1\in \Omega_{1,t}$ be such that
$B(\lambda_1)=B(a_1)$
and $ \lambda_1 \stackrel{B}{\simeq} a_1$, with
$a_1,\cdots,a_n \in \Omega_{0,t}$ (See Figure \ref{Fig-9}).  
Since  $B$ is injective on critical points,
 by   Corollary \ref{component}  we obtain
$$\rho_{\lambda_1, a_1}\sim  \cdots \sim \rho_{\lambda_1, a_n}. $$
By Remark \ref{remxyyx}(2), we have $\rho_{a_1,\lambda_1}  \sim \cdots \sim \rho_{a_n,\lambda_1} $.
In particular,
\begin{equation} \label{sett}
\rho_{a_2,\lambda_1}  \sim \cdots \sim \rho_{a_n,\lambda_1} .
\end{equation}

Next we prove that $\rho_{a_2,\lambda_1}\sim \rho_{a_2,a_1} $.
Since $ \lambda_1 \stackrel{B}{\simeq} a_1$,
$\lambda_1$ and $a_1$ are connected by   $B$-equal curves
ending at $p;$ that is, there exist curves
$\wp_1,\wp_2$ whose ranges
 lie in $\Omega_{0,t^+}$, satisfying
  $\wp_1(0)=\lambda_1,\wp_2(0)=a_1,$
 $\wp_1(1)= \wp_2(1)=p$ and $B(\wp_1(u))=B(\wp_2(u)),0\leq   u\leq    1$
 (in fact, for any $t'>t,$
 the ranges are contained in
  the closure of the component of $\{|B|<t' \}$ which contains
   $\Omega_{0,t }$, as well as  the nearby petals around
   $\Omega_{0,t }$).  Then by Lemma \ref{component3}   $p$ is a critical point.
  Again by Lemma \ref{component3},
  there exists another $B$-equal curve $\wp_0$ beginning   
   from $a_2$ and ending at a distinct point $v \neq p$;
   that is, $\wp_1(0)=v,\wp_2(0)=p,$ and
   $$B(\wp_0(u))=B(\wp_1(u))=B(\wp_2(u)), \ \ \ 0\leq     u\leq   1.$$
By the assumption of the flower model, $p$ serves as the only
 critical point on $\partial \Omega_{j,t}$.  Consequently,
 at the endpoint $v=\gamma_1(0) $, we must have
 $B'(v)\neq 0.$ Applying  Lemma  \ref{gluel},
we obtian  $\rho_{a_2,\lambda_1}\sim \rho_{a_2,a_1} $, as desired.

 By similar reasoning,
 $\rho_{a_j,\lambda_1}\sim  \rho_{a_j,a_1} $ for $j\geq    2$.
  This, combined with the  identity (\ref{sett}),
  shows that
$$\rho_{a_2,a_1} \sim \cdots \sim \rho_{a_n,a_1}.$$
Therefore,  by Proposition \ref{inversecomponent}, we obtain
 $\rho_{a_1,a_2}\sim \cdots \sim \rho_{a_1,a_n}.$
 Combining this relation with \eqref{equij2},
 we obtain that
\begin{equation}\label{rotation}
\rho_{a_2,a_3} \sim \rho_{a_2,a_4} \sim \cdots \sim \rho_{a_2,a_1}.
\end{equation}
 As in Case  I.1, we conclude that
\[
\rho_{a_2,a_1} \sim \rho_{a_2,\Omega_{j,t}}, \qquad j\geq1.
\]
Combining this with \eqref{rotation}, we obtain
\[
\rho_{a_2,a_1} \sim \rho_{a_2,z'}
\]
for every \(z'\in B^{-1}(w_0)\cap\Omega_{0,t^+}\)
with \(z'\ne a_2\) where  $\rho_{a_2,a_2}=id  $.
Hence,  on \(\Omega_{0,t^+} \backslash B^{-1}(\mathcal{E}_B)\),
  there are exactly two equivalence classes of local
   inverses of \(B\). By Proposition~\ref{newequi},
   the Riemann surface
$$\{(z,w)\in \Omega_{0,t^+}: B(z)=B(w)\not\in \mathcal{E}_B\}$$
admits exactly   two components: the one containing \((a_2,a_1)\)
 and the identity component \(\mathcal{S}_0\).

\textbf{Step II}. We proceed to address the general case,
 hereafter referred to as  \textbf{the Cactuses Case.}
The name will become clear as we see the shape of the branching
diagram shown later.

As $s$ grows bigger,   boundaries of distinct  components of $\{|B|<s\}$
intersect at one or more common points.
 Theorem~\ref{component2} guarantees
  the existence of a finite-vertex tree structure.
  After relabeling, let $\Omega_1,\cdots,\Omega_n$ be the components of $\{|B|<t\}$ whose boundaries intersect those of other components, chosen such that their union is maximal with respect to inclusion, subject to the condition that $\bigcup\limits_{i=1}^n \partial \Omega_i$ is connected.
We now treat $\Omega_1,\cdots,\Omega_N$ as abstract vertices.
  By Theorem \ref{component2}, any nonempty intersection
 $\partial  \Omega_i\cap \partial  \Omega_j $
   $(i\neq j)$,
   contains
 exactly one point $P$.  Then  we construct edges by
 connecting $\Omega_1$ and $\Omega_2 $
 for each such pair, and
 mark the intersection  point  $P $ as a vertex.
  By Theorem \ref{component2}, this procedure yields a connected
  graph without cycles. We denote the graph vertices as
  $\{\Omega_1,\Omega_2,\cdots, \Omega_n\} $ and
   $\{P_1,\cdots,P_m \}$, where the latter denote a subset of
   the critical set of $B$ (See Figure \ref{202610}).

 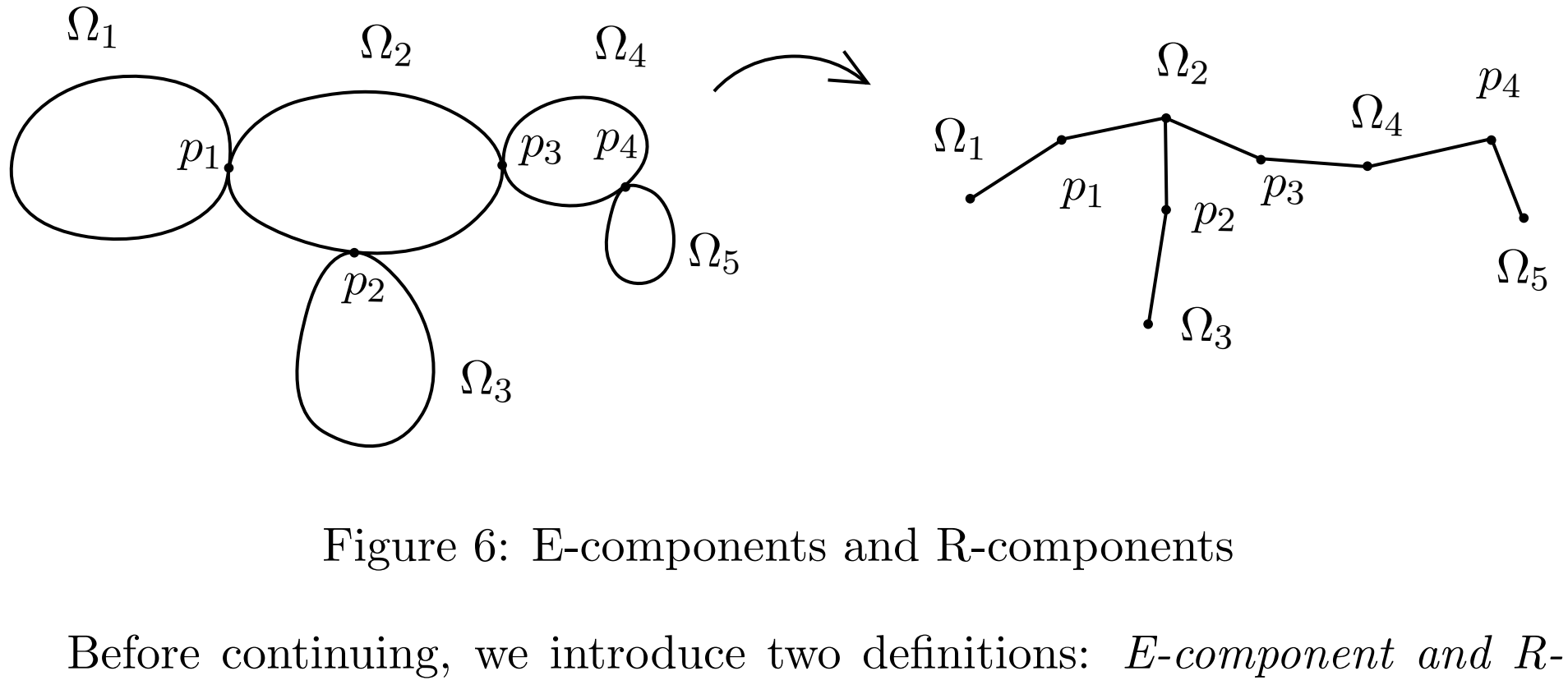
\begin{figure}[htbp]
\centering
\tikzset{every picture/.style={line width=0.75pt}} 
\begin{tikzpicture}[x=0.75pt,y=0.75pt,yscale=-1,xscale=1]
\draw (251.65,247.85) .. controls (235.68,261.77) and (206.34,258.98) .. (186.93,247.88) .. controls (167.53,236.77) and (180.24,215.31) .. (199.97,210.43) .. controls (219.69,205.55) and (235.84,208.26) .. (248.26,215.51) .. controls (260.69,222.75) and (267.62,233.93) .. (251.65,247.85) -- cycle ;
\draw [line width=3] [line join = round][line cap = round] (177.44,231.21) .. controls (177.44,231.21) and (177.44,231.21) .. (177.44,231.21) ;
\fill[black] (177.44,231.21) circle (1.1pt);
\draw [line width=3] [line join = round][line cap = round] (298.56,237.03) .. controls (298.56,237.03) and (298.56,237.03) .. (298.56,237.03) ;
\fill[black] (298.56,237.03) circle (1.1pt);
\draw (176.88,234.8) .. controls (174.22,247.69) and (153.62,255.94) .. (132.51,251.93) .. controls (111.39,247.91) and (108.7,235.31) .. (112.4,223.81) .. controls (116.1,212.32) and (132.53,200.59) .. (155.33,203.71) .. controls (178.12,206.84) and (179.53,221.91) .. (176.88,234.8) -- cycle ;
\draw (295.48,263.25) .. controls (291.81,258.36) and (292.1,252.62) .. (294.5,243.91) .. controls (296.9,235.2) and (300.63,235.59) .. (305.98,238.12) .. controls (311.33,240.65) and (315.23,250.71) .. (312.34,259.36) .. controls (309.45,268.02) and (299.15,268.15) .. (295.48,263.25) -- cycle ;
\draw (302.55,232.08) .. controls (296.42,240.4) and (286.31,245.16) .. (273.34,241.51) .. controls (260.38,237.87) and (259.65,229.33) .. (262.86,221.81) .. controls (266.07,214.3) and (279.02,207.43) .. (291.38,210.33) .. controls (303.74,213.23) and (308.68,223.77) .. (302.55,232.08) -- cycle ;
\draw (237.16,305.21) .. controls (230.87,317.16) and (219.46,319.25) .. (206.76,312.06) .. controls (194.05,304.87) and (198.6,284.94) .. (201.34,274.84) .. controls (204.08,264.74) and (211.36,249.31) .. (225.17,261.45) .. controls (238.99,273.6) and (243.44,293.27) .. (237.16,305.21) -- cycle ;
\draw [line width=3] [line join = round][line cap = round] (260.82,230.38) .. controls (260.82,230.38) and (260.82,230.38) .. (260.82,230.38) ;
\fill[black] (260.9,230.38) circle (1.1pt);
\draw [line width=3] [line join = round][line cap = round] (215.75,257.12) .. controls (215.75,257.12) and (215.75,257.12) .. (215.75,257.12) ;
\fill[black] (215.75,257.12) circle (1.1pt);
\draw [draw opacity=0] (325.61,207.8) .. controls (331.43,201.26) and (340.25,197.1) .. (350.12,197.1) .. controls (358.74,197.1) and (366.55,200.28) .. (372.27,205.43) -- (350.12,226.47) -- cycle ; \draw (325.61,207.8) .. controls (331.43,201.26) and (340.25,197.1) .. (350.12,197.1) .. controls (358.74,197.1) and (366.55,200.28) .. (372.27,205.43) ;
\draw (365.2,195.54) -- (372.45,205.25) -- (360.37,204.35) ;
\draw [line width=3] [line join = round][line cap = round] (431.71,222.66) .. controls (431.71,222.66) and (431.71,222.66) .. (431.71,222.66) ;
\fill[black] (431.71,222.66) circle (1.1pt);
\draw [line width=3] [line join = round][line cap = round] (403.8,240.6) .. controls (403.8,240.6) and (403.8,240.6) .. (403.8,240.6) ;
\fill[black] (403.8,240.6) circle (1.1pt);
\draw [line width=3] [line join = round][line cap = round] (492.58,229.15) .. controls (492.58,229.15) and (492.58,229.15) .. (492.58,229.15) ;
\fill[black] (492.58,228.7) circle (1.1pt);
\draw [line width=3] [line join = round][line cap = round] (463.64,216.89) .. controls (463.64,216.89) and (463.64,216.89) .. (463.64,216.89) ;
\fill[black] (463.64,216) circle (1.1pt);
\draw [line width=3] [line join = round][line cap = round] (562.94,222.71) .. controls (562.94,222.71) and (562.94,222.71) .. (562.94,222.71) ;
\fill[black] (562.94,222.71) circle (1.1pt);
\draw [line width=3] [line join = round][line cap = round] (458.14,278.93) .. controls (458.14,278.93) and (458.14,278.93) .. (458.14,278.93) ;
\fill[black] (458.14,278.93) circle (1.1pt);
\draw [line width=3] [line join = round][line cap = round] (525.11,230.67) .. controls (525.11,230.67) and (525.11,230.67) .. (525.11,230.67) ;
\fill[black] (525.11,230.67) circle (1.1pt);
\draw (463.38,215.91) -- (463.77,243.86) ;
\draw [line width=3] [line join = round][line cap = round] (463.64,244.17) .. controls (463.64,244.17) and (463.64,244.17) .. (463.64,244.17) ;
\fill[black] (463.64,244) circle (1.1pt);
\draw [line width=3] [line join = round][line cap = round] (572.86,246.52) .. controls (572.86,246.52) and (572.86,246.52) .. (572.86,246.52) ;
\fill[black] (572.86,246.52) circle (1.1pt);
\draw (562.99,222.31) -- (572.83,247.51) ;
\draw (403.67,240.81) -- (432.27,222.36) ;
\draw (432.27,222.36) -- (463.38,215.91) ;
\draw (463.77,243.86) -- (458.27,278.86) ;
\draw (463.38,215.91) -- (492.27,228.36) ;
\draw (492.27,228.36) -- (525.27,230.86) ;
\draw (525.27,230.86) -- (562.99,222.31) ;

\draw (430.35,233.01) node [anchor=north west][inner sep=0.75pt] [rotate=-1.24] [align=left] {$p_1$};
\draw (459.17,191.52) node [anchor=north west][inner sep=0.75pt] [rotate=-1.24] [align=left] {$\Omega_2$};
\draw (391.04,214.6) node [anchor=north west][inner sep=0.75pt] [rotate=-1.24] [align=left] {$\Omega_1$};
\draw (126.37,180.54) node [anchor=north west][inner sep=0.75pt] [rotate=-1.24] [align=left] {$\Omega_1$};
\draw (160.81,221.18) node [anchor=north west][inner sep=0.75pt] [rotate=-1.24] [align=left] {$p_1$};
\draw (215.98,186.04) node [anchor=north west][inner sep=0.75pt] [rotate=-1.24] [align=left] {$\Omega_2$};
\draw (264.96,219.78) node [anchor=north west][inner sep=0.75pt] [rotate=-1.24] [align=left] {$p_3$};
\draw (288,218) node [anchor=north west][inner sep=0.75pt] [rotate=-1.24] [align=left] {$p_4$};
\draw (287.56,186.4) node [anchor=north west][inner sep=0.75pt] [rotate=-1.24] [align=left] {$\Omega_4$};
\draw (210.92,261.92) node [anchor=north west][inner sep=0.75pt] [rotate=-1.24] [align=left] {$p_2$};
\draw (246.56,288.59) node [anchor=north west][inner sep=0.75pt] [rotate=-1.24] [align=left] {$\Omega_3$};
\draw (316.19,249.66) node [anchor=north west][inner sep=0.75pt] [rotate=-1.24] [align=left] {$\Omega_5$};
\draw (470.48,240.35) node [anchor=north west][inner sep=0.75pt] [rotate=-1.24] [align=left] {$p_2$};
\draw (491.64,232.02) node [anchor=north west][inner sep=0.75pt] [rotate=-1.24] [align=left] {$p_3$};
\draw (518.48,207.35) node [anchor=north west][inner sep=0.75pt] [rotate=-1.24] [align=left] {$\Omega_4$};
\draw (557.48,199.35) node [anchor=north west][inner sep=0.75pt] [rotate=-1.24] [align=left] {$p_4$};
\draw (563,254.68) node [anchor=north west][inner sep=0.75pt] [rotate=-1.24] [align=left] {$\Omega_5$};
\draw (466.48,272.35) node [anchor=north west][inner sep=0.75pt] [rotate=-1.24] [align=left] {$\Omega_3$};
\end{tikzpicture}
\caption{E-components and R-components}
\label{202610}
\end{figure}

Before continuing, we introduce two definitions:
  \emph{E-component and R-component}.
Given components $\Omega_1,\cdots,\Omega_n$ of $\{|B|<t\}$,
if one component $\Omega_k$ is
 such that there exist two distinct components $\Omega_i$ and
 $\Omega_j$ satisfying
 $$\partial \Omega_i \cap \partial \Omega_k \neq \emptyset,\ \ \ \
 \partial \Omega_j \cap \partial \Omega_k \neq \emptyset$$
 and $\partial \Omega_i \cap \partial \Omega_k \neq
 \partial \Omega_j \cap \partial \Omega_k $, then $\Omega_k$ is called
 an \emph{E-component} (edge-component).
 If a component $\Omega_k$ is
 such that there is a  component $\Omega_i$   satisfying
  $\partial \Omega_i \cap \partial \Omega_k \neq \emptyset$
and  $\Omega_k$ is not an E-component, then
$\Omega_k$ is called an \emph{R-component}
(root component).

The resulting graph \(\mathcal{G}\) contains at
least two \(R\)-components; here, a component
 \(\Omega_k\) is an \(R\)-component exactly
  when it acts like  a root
   when regarded as a vertex in  $\mathcal{G}$.
   This usage is not to be confused with the
   standard notion of a root vertex in graph theory.

\textbf{Case II.1.} \emph{There exists one $R$-component of multiplicity no less than two.}

  Let $\Omega_{0,t}$ be such an $R$-component. There are two situations under consideration:
either $\Omega_{0,t}$ contains at least two distinct critical points (and by induction it is purely
gluable) or $\Omega_{0,t}$ contains exactly one critical point.

\emph{We first consider the case where \(\Omega_{0,t}\) contains at least two distinct critical points;
 the remaining case is similar}.
  As done in Case I.2, we can pick $w_0 \in t\mathbb{D}\backslash \mathcal{E}_B$
and let $a_1,\cdots, a_n$ be all preimages of $w_0$ in $\Omega_{0,t}$.
Since $\Omega_{0,t}$ is an $R$-component, there exist
$\Omega_{1,t}^*$, $\cdots, \Omega_{N,t}^*$ ($N\geq 1$) such that
$$\partial \Omega_{j,t}^*\cap \partial \Omega_{0,t}$$
share a common critical point $*$ for $1\leq j \leq N.$
These are petals of the flower with critical point $*$.
Let $\gamma$ be a non-self-intersecting curve
connecting $w_0$ and $B(*)$, and
$$\gamma([0,1))\subseteq t\mathbb{D} \backslash \mathcal{E}_B.$$
Then by Lemma \ref{component3}, in each $\Omega_{j,t}^*$ ($1\leq j\leq N$) there exists exactly one
curve $\sigma_j$ with endpoint $*$ ($\sigma_j(1)=*$), and
$$B(\sigma_j(u))=\gamma(u), \qquad 0\leq u \leq 1.$$
Write $\lambda_j =\sigma_j(0).$
In $\Omega_{0,t}$, the preimage
$B^{-1}(\gamma)\cap \Omega_{0,t}$ of $\gamma$ contains
a curve connecting $a_1$ and $*$.
Again by Lemma \ref{component3},
there exists a curve $\sigma_0$ with $\sigma_0(0)=a_1$
and $\sigma_0(1)=*$, such that
$$B(\sigma_0(u))=\gamma(u), \quad 0 \leq u\leq 1.$$
Therefore,
$$B(\sigma_0(u))= B(\sigma_j(u))=\gamma(u), \quad 0\leq u\leq 1,$$
meaning $a_1$ and $\lambda_j$ are both connected
to the point $*$ via $B$-equal curves.
Applying Lemma \ref{gluel}, we obtain
\begin{equation} \label{eqqn1}  \rho_{a_2,a_1}
\sim \rho_{a_2,\lambda_j},\ \ \  1\leq     j\leq    N.
\end{equation} Then by Corollary \ref{component},
   $$\rho_{a_2,a_1}\sim \rho_{a_2,\lambda_j}\sim
   \rho_{a_2,\Omega_{j,t}},\ \ \  1\leq   j\leq     N;$$
remind that  for a set $G$,   $\rho_{a_2, G }$ means with $w_0=B(a_2)$,
$$\rho_{a_2,z'}\sim \rho_{a_2,z''},\ \
   z', z''\in B^{-1}(w_0)\cap G.  $$
Furthermore,
for certain domains \(\Omega\) (to be specified shortly), we will prove that
\[
\rho_{a_2,a_1} \sim \rho_{a_2,z}, \qquad z \in B^{-1}(w_0)\cap \Omega.
\]
More precisely, these domains lie in a chain of components associated with some \(\Omega_{j,t}^*\), as illustrated below. For instance, there exist components
\[
\Omega_{2,t}, \cdots, \Omega_{k,t}
\]
with \(\Omega_{k,t}=\Omega\) such that
 $$\partial \Omega_{1,t}^*\cap \partial \Omega_{2,t}\neq \emptyset,\
 \partial \Omega_{2,t} \cap \partial \Omega_{3,t}\neq \emptyset,\ \cdots,\
 \partial \Omega_{k-1,t} \cap \partial \Omega_{k,t}\neq \emptyset . $$
 The sequence $(\Omega_{1,t}^*,\Omega_{2,t},\cdots,\Omega_{k,t}) $
  is then called a \emph{chain} (see Figure \ref{Fig-11}).

\begin{figure}[htbp]
\centering
\tikzset{every picture/.style={line width=0.75pt}} 

\begin{tikzpicture}[x=0.75pt,y=0.75pt,yscale=-1,xscale=1]

\draw   (225.96,169.46) .. controls (234.25,173.31) and (242.9,183.28) .. (248.55,192.11) .. controls (254.2,200.94) and (254.23,203.16) .. (244.25,202.48) .. controls (234.28,201.8) and (210.63,191.92) .. (205.27,187.24) .. controls (199.91,182.57) and (196.16,170.7) .. (201.82,165.51) .. controls (207.47,160.32) and (217.66,165.6) .. (225.96,169.46) -- cycle ;
\draw   (275.05,173.03) .. controls (269.42,181.18) and (256.12,202.38) .. (252.87,200.24) .. controls (249.62,198.11) and (249.57,171.05) .. (250.67,164.18) .. controls (251.77,157.31) and (261.28,137.88) .. (269.36,135.45) .. controls (277.44,133.03) and (282.07,139.11) .. (282.78,142.74) .. controls (283.49,146.37) and (283.93,151.33) .. (282.88,154.39) .. controls (281.84,157.44) and (280.69,164.88) .. (275.05,173.03) -- cycle ;
\draw  [line width=3] [line join = round][line cap = round] (253.26,200.66) .. controls (253.26,200.66) and (253.26,200.66) .. (253.26,200.66) ;
\draw   (175.8,275.14) .. controls (175.8,235.1) and (208.26,202.64) .. (248.3,202.64) .. controls (288.34,202.64) and (320.8,235.1) .. (320.8,275.14) .. controls (320.8,315.19) and (288.34,347.64) .. (248.3,347.64) .. controls (208.26,347.64) and (175.8,315.19) .. (175.8,275.14) -- cycle ;
\draw   (343.58,152.66) .. controls (363.51,149.97) and (377.32,159.28) .. (376.22,168.14) .. controls (375.11,177) and (354.65,183.37) .. (348.32,184.32) .. controls (341.99,185.27) and (330.69,186.81) .. (321.82,178.19) .. controls (312.95,169.57) and (323.64,155.35) .. (343.58,152.66) -- cycle ;
\draw   (396.87,150.07) .. controls (416.8,147.38) and (420.34,154.81) .. (419.24,163.67) .. controls (418.13,172.53) and (409.65,176.37) .. (403.32,177.32) .. controls (396.99,178.27) and (388.74,181.69) .. (379.87,173.07) .. controls (371,164.45) and (376.93,152.75) .. (396.87,150.07) -- cycle ;
\draw   (419.24,157.15) .. controls (419.24,146.94) and (427.51,138.67) .. (437.72,138.67) .. controls (447.93,138.67) and (456.2,146.94) .. (456.2,157.15) .. controls (456.2,167.35) and (447.93,175.63) .. (437.72,175.63) .. controls (427.51,175.63) and (419.24,167.35) .. (419.24,157.15) -- cycle ;
\draw   (456.2,157.15) .. controls (456.2,146.94) and (464.47,138.67) .. (474.68,138.67) .. controls (484.89,138.67) and (493.16,146.94) .. (493.16,157.15) .. controls (493.16,167.35) and (484.89,175.63) .. (474.68,175.63) .. controls (464.47,175.63) and (456.2,167.35) .. (456.2,157.15) -- cycle ;
\draw   (288.94,205.46) .. controls (279.03,205.19) and (254.83,203.07) .. (254.02,201.76) .. controls (253.22,200.44) and (276.63,184.53) .. (287.55,180.44) .. controls (298.47,176.36) and (315.37,172.74) .. (321.82,178.19) .. controls (328.27,183.63) and (321.17,192.79) .. (318.52,195.37) .. controls (315.87,197.94) and (311.96,201.02) .. (308.83,201.82) .. controls (305.7,202.62) and (298.84,205.72) .. (288.94,205.46) -- cycle ;
\draw    (175.92,269.99) .. controls (202.36,258.68) and (191.9,281.26) .. (214.52,270.25) ;
\draw  [line width=3] [line join = round][line cap = round] (213.34,270.48) .. controls (213.34,270.48) and (213.34,270.48) .. (213.34,270.48) ;

\draw    (251.97,200.31) .. controls (264.38,226.26) and (241.38,216.76) .. (253.34,238.89) ;
\draw  [line width=3] [line join = round][line cap = round] (253.06,237.33) .. controls (253.06,237.33) and (253.06,237.33) .. (253.06,237.33) ;

\node at (253.26,200.66) [circle, fill=black, inner sep=1.2pt] {}; %
\node at (214.26,270.66) [circle, fill=black, inner sep=1.2pt] {}; %
\node at (176.26,270.26) [circle, fill=black, inner sep=1.2pt] {}; %
\draw (240.89,119) node [anchor=north west][inner sep=0.75pt]  [rotate=-358.49] [align=left] {$\Omega^*_{2,t}$};
\draw (240.89,205) node [anchor=north west][inner sep=0.75pt]  [rotate=-358.49] [align=left] {$p$};
\draw (248.89,245) node [anchor=north west][inner sep=0.75pt]  [rotate=-358.49] [align=left] {$a_1$};
\draw (218.89,270) node [anchor=north west][inner sep=0.75pt]  [rotate=-358.49] [align=left] {$a_2$};
\draw (160,270) node [anchor=north west][inner sep=0.75pt]  [rotate=-358.49] [align=left] {$v$};
\draw (287.19,156.93) node [anchor=north west][inner sep=0.75pt]  [rotate=-2.08] [align=left] {$\Omega^*_{1,t}$};

\draw (431.65,119.91) node [anchor=north west][inner sep=0.75pt]  [rotate=-359.39] [align=left] {$\Omega_{k,t}$};
\draw (338.65,133.91) node [anchor=north west][inner sep=0.75pt]  [rotate=-359.39] [align=left] {$\Omega_{2,t}$};
\draw (314.65,316.91) node [anchor=north west][inner sep=0.75pt]  [rotate=-359.39] [align=left] {$\Omega_{0,t}$};

\end{tikzpicture}
\caption{Chain}
\label{Fig-11}
\end{figure}

 By Corollary \ref{component},
 we clearly have the relationships:
 $$\rho_{a_2,\Omega_{1,t}^*}  ,  \rho_{a_2,\Omega_{2,t} } ,
  \cdots,  \rho_{a_2, \Omega_{k,t} }.$$
      We will verify that
\[
\rho_{a_2,\Omega_{i,t}} \sim \rho_{a_2,\Omega_{i+1,t}},
\quad  i\geq2,
\]
and also
\[
\rho_{a_2,\Omega_{1,t}^*} \sim \rho_{a_2,\Omega_{2,t}}.
\]
For this, Theorem~\ref{component2}(1) implies that
\(\partial \Omega_{i,t} \cap \partial \Omega_{i+1,t}\) contains exactly one critical point, denoted by \(Q\). Let \(w_0 = B(a_1)\), and let \(\wp\) be a non-self-intersecting curve connecting \(w_0\) to \(B(Q)\) such that
\[
\wp(u) \in t  \mathbb{D} \backslash  \mathcal{E}_B, \qquad 0\leq  u<1.
\]
By Lemma~\ref{component3}, there are two curves,
lying in \(\Omega_{i,t}\) and \(\Omega_{i+1,t}\)
respectively, connecting a point in \(B^{-1}(w_0)\)
to \(Q\). Let \(\lambda\) and \(\mu\) be
the endpoints of these two curves lying in
 \(B^{-1}(w_0)\); then $\lambda \overset{B}{ \simeq } \mu$.
   Applying Lemma~\ref{gluel} immediately yields
\[
\rho_{a_2,\lambda} \sim \rho_{a_2,\mu}.
\]
Combined with the equivalences
\[
\rho_{a_2,\lambda} \sim \rho_{a_2,\Omega_{i,t}}
\quad\text{and}\quad
\rho_{a_2,\mu} \sim \rho_{a_2,\Omega_{i+1,t}},
\]
we conclude that
\[
\rho_{a_2,\Omega_{i,t}} \sim \rho_{a_2,\Omega_{i+1,t}}.
\]
The argument for \(\rho_{a_2,\Omega_{1,t}^*} \sim \rho_{a_2,\Omega_{2,t}}\) is similar.

Combining the above identities with \eqref{eqqn1}, we obtain
\begin{equation}\label{add}
\rho_{a_2,a_1} \sim \rho_{a_2,\lambda_1} \sim
\rho_{a_2,\Omega_{1,t}^*} \sim \rho_{a_2,\Omega_{2,t}} \sim
\cdots \sim \rho_{a_2,\Omega_{k,t}}.
\end{equation}
Since $\Omega_{0,t}$ is purely gluable, we have
$\rho_{a_2,a_1} \sim \rho_{a_2,a_j}$ for all $j\neq 2$.
This, together with \eqref{add}, shows that
\[
\rho_{a_2,a_1} \sim \rho_{a_2,z'}
\]
for every \(z'\in B^{-1}(w_0)\cap\Omega_{0,t^+}\) with \(z'\ne a_2\).
Hence, on $\Omega_{0,t^+} \backslash B^{-1}(\mathcal{E}_B)$,
there are exactly two equivalence classes of local inverses of $B$.
By Proposition~\ref{newequi}, the Riemann surface
\[
\{(z,w)\in \Omega_{0,t^+}: B(z)=B(w)\notin \mathcal{E}_B\}
\]
has exactly two components.
That is, $\Omega_{0,t^+}$ is purely gluable.

We now address the case where
  $\Omega_{0,t}$ contains exactly one critical point.
 By the arguments developed above, we immediately derive
   (\ref{add}). Following the reasoning in Rotation Case of Case II,
we further establish the relation (\ref{rotation}).
As a result, the conclusion that $\Omega_{0,t^+} $ is purely gluable
 still holds.

\textbf{Case II.2.} \emph{All  $R$-components are of multiplicity one.}

\begin{figure}[htbp]
\centering
\tikzset{every picture/.style={line width=0.75pt}} 

\begin{tikzpicture}[x=0.65pt,y=0.65pt,yscale=-1,xscale=1]

\draw   (345.9,225.41) .. controls (352.32,242.7) and (348.42,264.71) .. (344.89,279.81) .. controls (341.35,294.91) and (339.36,298.88) .. (328.74,324.8) .. controls (318.12,350.73) and (308.84,380.56) .. (279.37,335.8) .. controls (249.89,291.04) and (252.24,270.78) .. (252.37,265.11) .. controls (252.51,259.45) and (258.29,198.06) .. (287.5,188.07) .. controls (316.71,178.09) and (339.48,208.13) .. (345.9,225.41) -- cycle ;
\draw   (326.77,388.8) .. controls (324.87,384.76) and (293.4,355.86) .. (313,358.41) .. controls (332.6,360.97) and (353.91,381.11) .. (358.7,390.9) .. controls (363.5,400.69) and (355.63,408.83) .. (349.33,407.68) .. controls (343.02,406.53) and (328.68,392.85) .. (326.77,388.8) -- cycle ;
\draw   (268.52,342.76) .. controls (276.04,345.12) and (277.14,344.97) .. (281.99,347.73) .. controls (286.84,350.49) and (300.83,357.08) .. (300.19,359.04) .. controls (299.56,360.99) and (278.26,367.67) .. (266.23,367.56) .. controls (254.2,367.46) and (233.78,359.6) .. (241.67,350.12) .. controls (249.57,340.64) and (261,340.4) .. (268.52,342.76) -- cycle ;
\draw  [line width=3] [line join = round][line cap = round] (328.61,381.07) .. controls (328.61,381.07) and (328.61,381.07) .. (328.61,381.07) ;
\draw   (280.03,377.35) .. controls (287.7,370.25) and (299.88,355.04) .. (305.05,359.07) .. controls (310.21,363.1) and (304.86,384.95) .. (301.34,392.53) .. controls (297.83,400.12) and (292.74,403.66) .. (288.62,408.57) .. controls (284.5,413.47) and (272.4,417.85) .. (266.83,413.71) .. controls (261.26,409.58) and (262.09,402.3) .. (264.91,396.33) .. controls (267.72,390.36) and (272.36,384.45) .. (280.03,377.35) -- cycle ;
\draw  [line width=3] [line join = round][line cap = round] (289.8,379.99) .. controls (289.8,379.99) and (289.8,379.99) .. (289.8,379.99) ;
\draw  [line width=3] [line join = round][line cap = round] (276.26,355.39) .. controls (276.26,355.39) and (276.26,355.39) .. (276.26,355.39) ;
\draw  [line width=3] [line join = round][line cap = round] (300.06,318.55) .. controls (300.06,318.55) and (300.06,318.55) .. (300.06,318.55) ;
\draw  [line width=3] [line join = round][line cap = round] (305.06,359.39) .. controls (305.06,359.39) and (305.06,359.39) .. (305.06,359.39) ;
\draw    (302.01,319.71) .. controls (317.7,336.12) and (299.04,349.5) .. (305.05,359.07) .. controls (311.06,368.65) and (323.65,360.5) .. (328.49,379.76) ;
\draw    (305.05,359.07) .. controls (297.62,357.09) and (280.09,361.29) .. (276.58,354.03) ;
\draw    (288.98,380.55) .. controls (294.84,372.91) and (295.8,373.92) .. (298.37,370.55) .. controls (299.77,368.72) and (301.65,365.59) .. (305.05,359.07) ;
\draw   (273.92,149.97) .. controls (283.13,165.95) and (309.49,187.31) .. (289.74,186.41) .. controls (270,185.52) and (256.17,175.84) .. (241.48,157.88) .. controls (226.78,139.92) and (237.9,132.39) .. (247.08,132.52) .. controls (256.26,132.64) and (264.71,134) .. (273.92,149.97) -- cycle ;
\draw   (332.14,193.61) .. controls (324.87,192.26) and (326.18,191.12) .. (320.53,189.84) .. controls (314.87,188.57) and (300.18,186.64) .. (300.65,184.64) .. controls (301.11,182.63) and (323.21,175.23) .. (335.21,175.52) .. controls (347.2,175.81) and (366.47,178.91) .. (364.6,187.2) .. controls (362.74,195.49) and (339.41,194.96) .. (332.14,193.61) -- cycle ;
\draw  [line width=3] [line join = round][line cap = round] (272.62,165.13) .. controls (272.62,165.13) and (272.62,165.13) .. (272.62,165.13) ;
\draw   (324.18,162.59) .. controls (317.13,170.31) and (300.22,186.89) .. (297.95,185.07) .. controls (295.68,183.25) and (298.73,166.98) .. (304.14,152.39) .. controls (309.55,137.8) and (319.29,124.54) .. (324.81,122.26) .. controls (330.33,119.97) and (338.23,118.88) .. (339.81,127.54) .. controls (341.39,136.2) and (331.23,154.86) .. (324.18,162.59) -- cycle ;
\draw  [line width=3] [line join = round][line cap = round] (311.39,162.94) .. controls (311.39,162.94) and (311.39,162.94) .. (311.39,162.94) ;
\draw  [line width=3] [line join = round][line cap = round] (326.94,186.31) .. controls (326.94,186.31) and (326.94,186.31) .. (326.94,186.31) ;
\draw  [line width=3] [line join = round][line cap = round] (297.91,184.75) .. controls (297.91,184.75) and (297.91,184.75) .. (297.91,184.75) ;
\draw  [dash pattern={on 1.5pt off 1.5pt on 1.5pt off 1.5pt}]  (300,317.92) .. controls (298.1,310.06) and (307.6,279.52) .. (309.99,266.94) .. controls (312.37,254.36) and (307.27,228.91) .. (303.59,216.14) .. controls (299.9,203.37) and (299.75,187.46) .. (297.95,185.07) .. controls (296.15,182.67) and (279.29,185.21) .. (272.85,166.42) ;
\draw  [dash pattern={on 1.5pt off 1.5pt on 1.5pt off 1.5pt}]  (297.95,185.07) .. controls (305.51,186.41) and (322.63,180.76) .. (326.74,187.7) ;
\draw  [dash pattern={on 1.5pt off 1.5pt on 1.5pt off 1.5pt}]  (312.15,162.31) .. controls (304.14,174.83) and (306.01,165.81) .. (297.95,185.07) ;

\node at (300.06,318.55) [circle, fill=black, inner sep=1pt] {}; 
\node at (272.62,165.13) [circle, fill=black, inner sep=1pt] {}; 
\node at (328.61,381.07) [circle, fill=black, inner sep=1pt] {}; 
\node at (289.8,379.99) [circle, fill=black, inner sep=1pt] {}; 
\node at (276.26,355.39) [circle, fill=black, inner sep=1pt] {}; 
\node at (326.94,186.31) [circle, fill=black, inner sep=1pt] {}; 
\node at (311.39,162.94) [circle, fill=black, inner sep=1pt] {}; 
\node at (298.3,184.75) [circle, fill=black, inner sep=1.2pt] {}; 
\node at (305.06,359.39) [circle, fill=black, inner sep=1pt] {}; %

\draw (285,299) node [anchor=north west][inner sep=0.75pt]  [rotate=-359.77] [align=left] {$b_1$};
\draw (308.96,385.66) node [anchor=north west][inner sep=0.75pt]  [rotate=-21.15] [align=left] {$...$};
\draw (282.28,148.39) node [anchor=north west][inner sep=0.75pt]  [rotate=-327.94] [align=left] {$...$};
\draw (287.39,247.59) node [anchor=north west][inner sep=0.75pt]  [rotate=-1.91] [align=left] {$\Omega_3$};
\draw (329.55,177.45) node [anchor=north west][inner sep=0.75pt]  [rotate=-359.77] [align=left] {$a_1$};
\draw (309.44,139.96) node [anchor=north west][inner sep=0.75pt]  [rotate=-359.77] [align=left] {$a_2$};
\draw (253.19,149.62) node [anchor=north west][inner sep=0.75pt]  [rotate=-359.77] [align=left] {$a_n$};
\draw (258.55,343.8) node [anchor=north west][inner sep=0.75pt]  [rotate=-359.77] [align=left] {$b_2$};
\draw (256.05,184.74) node [anchor=north west][inner sep=0.75pt]  [rotate=-359.77] [align=left] {$P$};
\draw (324.58,343.38) node [anchor=north west][inner sep=0.75pt]  [rotate=-359.77] [align=left] {$Q$};

\end{tikzpicture}

\caption{Case 2-1}
\label{Fig-10}
\end{figure}
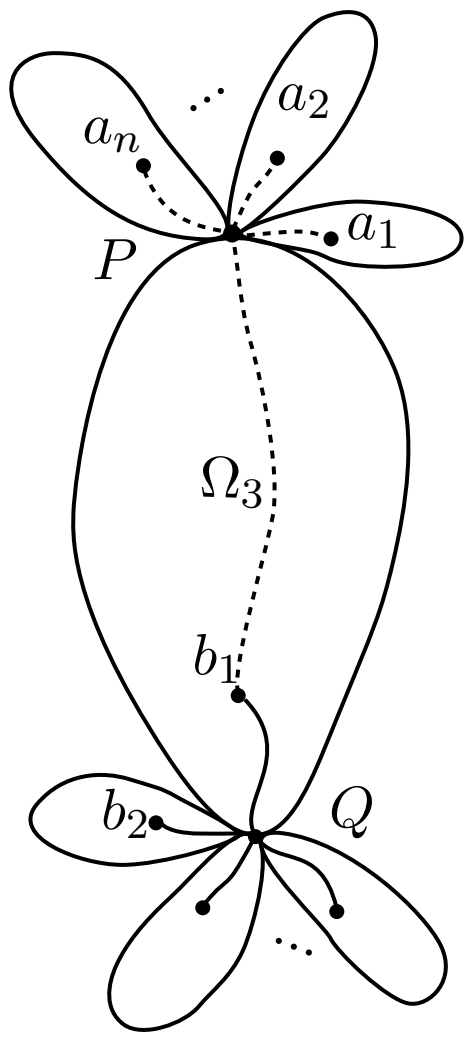

To illustrate the idea, we begin with a simplified example (see Figure \ref{Fig-10})
 before addressing the
 general case.

 \textbf{A special case.}  We assume,
  in this simplified case,
  that all components   are of multiplicity one,
   and that there exists only
 one $E$-component, denoted by $\Omega_3$, whose boundary $\partial \Omega_3$
 contains exactly two critical points $P$ and $Q$.
 There are $R$-components $\Omega_1^{1},\cdots, \Omega_n^{1} $ and
 $R$-components  $\Omega_2^{2},\cdots, \Omega_m^{2} $
   such that
 $$P\in \partial \Omega_i^{1}\cap \partial \Omega_3,\ \  \ 1\leq     i \leq    n$$
 and $$Q\in
 \partial \Omega_j^{2}\cap \partial \Omega_3,
  \ \ \  2\leq     j \leq   m.$$
   This labeling $\Omega_2^{2},\cdots, \Omega_m^{2} $
    is adopted for notational convenience in the sequel.
  By Theorem~\ref{component2}, \(P\) and \(Q\)
  are the only possible intersection points of
  the boundaries of any two such domains
  (among those components \(\Omega_j^i\) for \(i=1,2\), and \(\Omega_3\)).

Let
$w_0$ be a point in $t\mathbb{D} \backslash \mathcal{E}_B $,
and define:
\begin{itemize}
    \item $a_i $ as the unique point in $\Omega_i^1 \cap B^{-1}(w_0) $
     for $1 \leq    i \leq    n $,
    \item $ b_j $ as the unique point in $\Omega_j^2 \cap B^{-1}(w_0) $
     for $ 2 \leq     j \leq     m $.
\end{itemize}
In the special case, $ b_1 $ is the only point in
$ \Omega_3 \cap B^{-1}(w_0) $ (as illustrated in Figure \ref{Fig-10}).

We will soon obtain that
\begin{equation} \label{t1}
\rho_{b_2,b_1} \sim \rho_{b_2,a_i}, \quad 1\leq i\leq n.
\end{equation}
To this end, fix $1\leq i\leq n$. The relation $b_1 \overset{B}{\simeq} a_i$ implies the existence of curves $\wp_1$ and $\wp_2$ such that
$$
B(\wp_1(u))=B(\wp_2(u))\notin \mathcal{E}_B, \qquad 0\leq u<1,
$$
where $\wp_1$ connects $b_1$ to $P$, and $\wp_2$ connects $a_i$ to $P$. By Lemma~\ref{component3}, there exists a curve $\wp_0$ such that
$$
B(\wp_0(u))=B(\wp_1(u))=B(\wp_2(u)), \qquad 0\leq u\leq 1,
$$
with $\wp_0(0)=b_2$ and $\wp_0(u)\in \Omega_2^2$ for $0\leq u<1$. Let $v=\wp_0(1)$.

 We claim that \(B'(v)\neq 0\),
 which, together with Lemma~\ref{gluel}, gives
  \(\rho_{b_2,b_1}\sim \rho_{b_2,a_i}\), as desired.
To prove the claim  \(B'(v)\neq 0\), suppose for contradiction that \(v\) is a critical point of \(B\). Since \(B\) is injective on its critical set \(Z(B')\), we must have \(v=Q\). But \(B(\wp_0(1))=B(\wp_1(1))\) implies \(B(Q)=B(v)=B(P)\), contradicting the assumption that
 \(P\) and \(Q\) have distinct \(B\)-values.
This  finishes the proof of (\ref{t1}).
 By symmetry, we may exchange
 the positions   of $a_i$ and $b_j$ to obtain
 $$\rho_{a_1, b_2}\sim \cdots \sim \rho_{a_1, b_m} .$$
 By Remark \ref{remxyyx}(2),
\begin{equation} \label{abe}
\rho_{b_2,a_1}\sim \cdots \sim \rho_{b_m,a_1} .
\end{equation}
By (\ref{t1}) $\rho_{b_2,b_1}\sim \rho_{b_2,a_1} ,$
and similarly, $  \rho_{b_i,b_1}\sim \rho_{b_i,a_1}, i\geq   2. $
This, combined with (\ref{abe}), gives
$$\rho_{b_2,b_1}\sim \cdots\sim \rho_{b_m,b_1}.$$
By Remark \ref{remxyyx}(2), we have
$$\rho_{b_1,b_2}\sim \cdots\sim \rho_{b_1,b_m}.$$
 By applying the arguments developed in  Case I.2,
  we obtain an analogue of (\ref{rotation}),
$\rho_{b_2,b_3}\sim \cdots\sim \rho_{b_2,b_m}\sim \rho_{b_2,b_1}.$
This, together with (\ref{t1}),
imply that
$$ \rho_{b_2,b_i}\sim \rho_{b_2,b_1}\sim \rho_{b_2,a_j},\ \ \  i\geq   1, \  i\neq 2;\  \ j\geq   1. $$
That is,
\[
\rho_{b_2,b_1} \sim \rho_{b_2,z'}
\]
for every \(z'\in B^{-1}(w_0)\cap\Omega_{0,t^+}\)
with \(z'\ne b_2\) where  $\rho_{b_2,b_2}=id  $.
Hence,  on \(\Omega_{0,t^+}
\backslash B^{-1}(\mathcal{E}_B)\),
  there are exactly two equivalence classes of local
   inverses of \(B\). By Proposition~\ref{newequi},
   the Riemann surface
$ \mathcal{S}_{B,\Omega_{0,t^+}}  $
admits exactly   two components,
   the identity component $\mathcal{S}_0$
 and the  component $\mathcal{S} $
  of $\mathcal{S}_{B,\Omega_{0,t^+}}  $
 containing the point  $ (b_2,b_1)$.

\textbf{ The general case.}
   Let $\Omega_1 $ be an $R$-component whose boundary contains at least one critical point, say  $*$.  There are two situations
 to distinguish:
   $\Omega_1 $ is the unique $R$-component whose boundary contains $*$ or not.

   We first treat the former case.
   After relabeling, let \(\Omega_2',\cdots,\Omega_N'\)
   denote the maximal family of components of
    \(\{|B|<t\}\), other than \(\Omega_1\),
    such that their boundaries intersect those
    of other components and the union
    \(\bigcup\limits_{i=2}^N   \overline{\Omega_i'}\backslash \{*\} \)
    is connected. This union is
    like a cactus pup growing from a single node $\{*\}$.
   Recall that in the penultimate paragraph of Case  II.1,
    we discussed a chain of components.
   By the discussions there, by letting $b_1$ be the only point in
    $B^{-1}(w_0)\cap \Omega_1$ (which follows from the assumption that $\Omega_1 $
    is of multiplicity one),  we have
   $$\rho_{b_1,\Omega_2'}\sim  \cdots\sim \rho_{b_1,\Omega_N'},$$
 since each $\Omega_j' \ (2\leq    j \leq    N)$ is contained in some chain.
 Here,  for each $i$,  $\rho_{b_1,\Omega_i'}$ means
   $$\rho_{b_1,z'}\sim \rho_{b_1,z''} , \ \ \ \
   z',z''\in B^{-1}(w_0)\cap \Omega_i'. $$
In this way, we will have that
$\{(z,w)\in \Omega_{0,t^+}: B(z)=B(w)\not\in \mathcal{E}_B\}$ has exactly two components.

 It remains to consider the case where \(\Omega_1\) is not the only \(R\)-component whose boundary contains \(*\).

Below, we write $\Omega_i$ for $\Omega_{i,t}$ to simplify
 notation. After relabeling, let $\Omega_1,\cdots,\Omega_m$
  be the $R$-components of $\{|B|<t\}$ of multiplicity
  one whose boundaries contain $*$, and let
   $\Omega_{m+1},\cdots,\Omega_N$ be the remaining
    components of $\{|B|<t\}$ whose boundaries also
     contain $*$ (these are $E$-components).
      No other component $\Omega_{i}$ has $*$ on its boundary. Let
$$
W_0 = \bigcup_{i=1}^N \overline{\Omega}_i,
$$
and let $W$ be the component of $W_0 \backslash \{q\}$ containing $\overline{\Omega}_1 \backslash \{q\}$. By Theorem~\ref{component2}, the boundaries $\partial \Omega_i$ for $1\leq i\leq m$ intersect only at $*$:
$$
\partial \Omega_i \cap \partial \Omega_j = \{*\}, \qquad 1\leq i < j \leq m.
$$
Moreover, each $\partial \Omega_k$ ($1\leq k\leq m$) intersects $\partial \Omega_j$ only at $*$ for all $j\geq m+1$. Let $\gamma$ be a non-self-intersecting curve connecting $w_0$ to $B(*)$ such that
$$
\gamma(u) \in t\mathbb{D} \backslash \mathcal{E}_B, \qquad 0\leq u<1.
$$
By Lemma~\ref{component3}, for each $1\leq j\leq N$, there is exactly one curve $\sigma_j$ in $\Omega_j$ satisfying
$$
B(\sigma_j(u))=\gamma(u), \quad 0\leq u<1, \qquad \sigma_j(1)=*.
$$
Write $b_j= \sigma_j(0)$ for each $j$.

 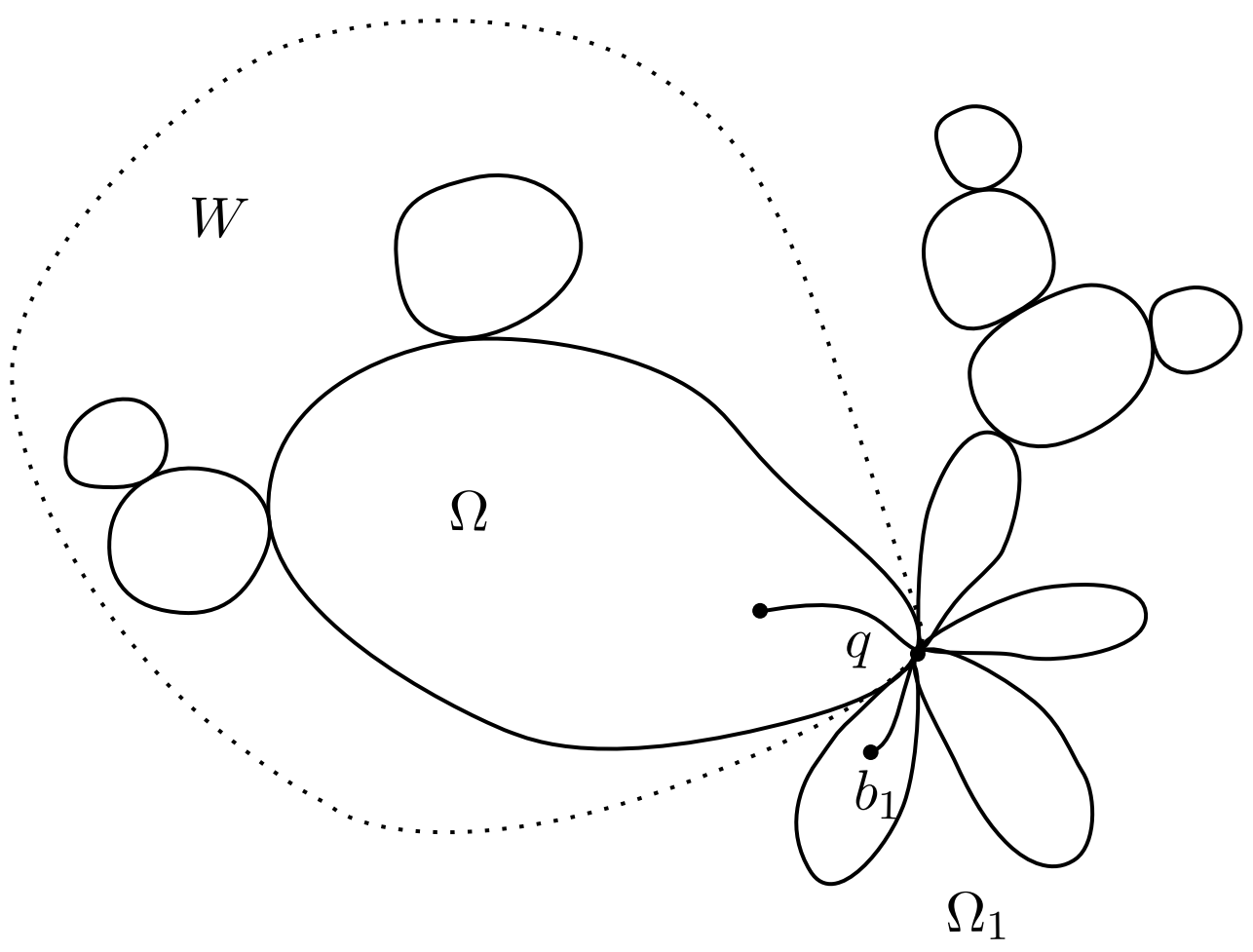
\begin{figure}
\centering
\tikzset{every picture/.style={line width=0.75pt}} 
\begin{tikzpicture}[x=0.75pt,y=0.75pt,yscale=-1,xscale=1]
\draw (278.07,475.33) .. controls (296.5,474.64) and (319.1,479.51) .. (332.34,487.58) .. controls (345.58,495.66) and (344.2,501.6) .. (365.51,519.79) .. controls (386.82,537.98) and (413.85,559.04) .. (362.12,573.05) .. controls (310.39,587.05) and (291.87,578.48) .. (286.54,576.58) .. controls (281.2,574.68) and (224.69,550) .. (224.34,519.13) .. controls (223.99,488.26) and (259.64,476.03) .. (278.07,475.33) -- cycle ;
\draw (428.42,502.26) .. controls (415.44,505.92) and (405.18,495.47) .. (404.97,484.32) .. controls (404.76,473.16) and (426.49,462.89) .. (434.11,461.66) .. controls (441.73,460.43) and (451.13,465.23) .. (452.11,476.96) .. controls (453.09,488.68) and (441.4,498.61) .. (428.42,502.26) -- cycle ;
\draw (365.38,584.63) .. controls (369.97,578.22) and (370.17,577.13) .. (374.31,573.39) .. controls (378.44,569.64) and (389.08,558.42) .. (390.74,559.63) .. controls (392.4,560.85) and (392.08,583.16) .. (388.22,594.56) .. controls (384.36,605.95) and (370.51,622.9) .. (363.97,612.43) .. controls (357.44,601.97) and (360.79,591.03) .. (365.38,584.63) -- cycle ;
\draw (401.79,585.41) .. controls (397.45,575.9) and (386.81,559.58) .. (392.25,555.93) .. controls (397.7,552.28) and (416.77,564.2) .. (422.88,569.91) .. controls (428.98,575.62) and (430.76,581.56) .. (434.13,587.01) .. controls (437.5,592.45) and (437.87,605.32) .. (432.2,609.32) .. controls (426.54,613.32) and (419.88,610.25) .. (415.09,605.7) .. controls (410.3,601.16) and (406.14,594.91) .. (401.79,585.41) -- cycle ;
\draw [line width=3] [line join = round][line cap = round] (379.76,582.12) .. controls (379.76,582.12) and (379.76,582.12) .. (379.76,582.12) ;
\draw [line width=3] [line join = round][line cap = round] (392.56,556.02) .. controls (392.56,556.02) and (392.56,556.02) .. (392.56,556.02) ;
\draw (392.25,555.93) .. controls (387.54,560.14) and (387.21,580.53) .. (379.22,581.59) ;
\draw (414.35,527.76) .. controls (412.98,530.77) and (414.11,530.4) .. (405.4,538.68) .. controls (396.68,546.97) and (394.43,555.6) .. (392.67,554.53) .. controls (390.91,553.45) and (391.65,528.03) .. (394.34,519.32) .. controls (397.02,510.62) and (404.34,495.1) .. (413.2,500.51) .. controls (422.05,505.91) and (415.72,524.75) .. (414.35,527.76) -- cycle ;
\draw (424.59,557.74) .. controls (418.05,557.68) and (418.72,556.63) .. (412.93,556.35) .. controls (407.14,556.07) and (392.56,556.58) .. (392.67,554.53) .. controls (392.79,552.47) and (412.63,541.17) .. (424.5,539.39) .. controls (436.36,537.61) and (450.77,538.44) .. (450.36,546.93) .. controls (449.95,555.42) and (431.13,557.8) .. (424.59,557.74) -- cycle ;
\draw (406.85,436.9) .. controls (401.81,436.48) and (399.37,433.06) .. (397.12,426.61) .. controls (394.88,420.17) and (397.87,418.01) .. (403.07,416.03) .. controls (408.27,414.04) and (415.88,417.22) .. (417.73,423.94) .. controls (419.58,430.67) and (411.89,437.31) .. (406.85,436.9) -- cycle ;
\draw (412.94,470.62) .. controls (400.93,476.75) and (396.03,468.28) .. (393.63,457.39) .. controls (391.23,446.49) and (397.67,440.53) .. (404.9,437.83) .. controls (412.13,435.13) and (422.29,437.99) .. (425.55,449.3) .. controls (428.81,460.6) and (424.95,464.49) .. (412.94,470.62) -- cycle ;
\draw (459.37,483.87) .. controls (454.01,482.74) and (452.3,478.94) .. (451.65,472.13) .. controls (451.01,465.32) and (454.85,463.57) .. (461.05,462.3) .. controls (467.25,461.04) and (474.63,465.32) .. (474.78,472.36) .. controls (474.92,479.4) and (464.73,485.01) .. (459.37,483.87) -- cycle ;
\draw (196.79,507.94) .. controls (193.49,512.95) and (188.46,513.96) .. (180.32,513.41) .. controls (172.19,512.85) and (171.58,508.81) .. (172.4,502.56) .. controls (173.23,496.32) and (180.92,489.88) .. (189.14,490.97) .. controls (197.35,492.05) and (200.08,502.93) .. (196.79,507.94) -- cycle ;
\draw (223.15,531.33) .. controls (218.28,542.34) and (211.08,547.74) .. (198.11,545.37) .. controls (185.14,542.99) and (182.44,533.7) .. (183.65,525.04) .. controls (184.86,516.38) and (193.66,507.18) .. (207.56,508.84) .. controls (221.46,510.5) and (228.02,520.32) .. (223.15,531.33) -- cycle ;
\draw (272.1,475.03) .. controls (261.17,473.35) and (257.98,466.28) .. (257.25,453.38) .. controls (256.52,440.48) and (264.59,436.81) .. (277.48,433.85) .. controls (290.36,430.89) and (305.19,438.38) .. (304.86,451.76) .. controls (304.53,465.14) and (283.04,476.72) .. (272.1,475.03) -- cycle ;
\draw [dash pattern={on 0.84pt off 2.51pt}] (232.26,398.63) .. controls (259.68,390.11) and (299.02,391.51) .. (321.72,405.71) .. controls (344.42,419.91) and (354.82,434.12) .. (371.17,484.09) .. controls (387.51,534.07) and (388.35,535.58) .. (392.69,550.08) .. controls (397.03,564.57) and (335.55,591.2) .. (302.13,598.59) .. controls (268.71,605.98) and (247.87,601.11) .. (241.35,596.38) .. controls (234.82,591.64) and (236.67,595.74) .. (210.64,575.37) .. controls (184.61,554.99) and (159.17,513.86) .. (158.3,485.18) .. controls (157.44,456.5) and (204.85,407.16) .. (232.26,398.63) -- cycle ;
\draw (392.25,555.93) .. controls (382.25,551.93) and (383.25,538.93) .. (349.25,545.93) ;
\draw [line width=3] [line join = round][line cap = round] (349.76,546.12) .. controls (349.76,546.12) and (349.76,546.12) .. (349.76,546.12) ;

\draw (371.41,549.35) node [anchor=north west][inner sep=0.75pt] [rotate=-1.24] [align=left] {$q$};
\draw (373.76,584.78) node [anchor=north west][inner sep=0.75pt] [rotate=-1.28] [align=left] {$b_1$};
\draw (397.13,616.53) node [anchor=north west][inner sep=0.75pt] [rotate=-359.34] [align=left] {$\Omega_1$};
\draw (202.41,437.35) node [anchor=north west][inner sep=0.75pt] [rotate=-1.24] [align=left] {$W$};
\draw (269.41,513.01) node [anchor=north west][inner sep=0.75pt] [rotate=-1.24] [align=left] {$\Omega$};

\fill (351, 545.35) circle (1.5pt);
\fill (379.5, 581.78) circle (1.5pt);
\fill (391.6, 556.45) circle (1.5pt);
\end{tikzpicture}
 \caption{Structure of $W_0$}
 \label{Fig-13}
\end{figure}

Consider the components
of $W  \backslash \{*\}.$
For $m+1\leq    k \leq     N$, let $W_{k}$ be the component of $W  \backslash \{*\} $
containing $\Omega_{k}$.
Since $\Omega_{m+1},\cdots,
\Omega_N$ are $E$-components, each such
 $W_k\cap B^{-1}(w_0)$ contains at least two points.
   We make the following claims.

  \textbf{Claim 1.} \emph{For $m+1\leq     k \leq    N$ and
  for $z',z''\in B^{-1}(w_0)\cap W_k,$ we
  have    $$\rho_{b_1,z'}\sim \rho_{b_1,z''}.$$ Then we write}
 \begin{equation} \label{last0}
\rho_{b_1,b_k}\sim  \rho_{b_1,W_k},\ \ \ \  m+1\leq     k \leq     N.
 \end{equation}

 To verify this, by Corollary~\ref{component}
 or Lemma~\ref{componenthigh}, we have the relation
$
\rho_{b_1,\Omega_{i,t}}$
 for all \(\Omega_{i,t}  \subseteq W_k\).
Fix one such component \(\Omega_{j,t}\) in \(W_k\), and denote it by \(\Omega_{k,t}^*\). For every component \(\Omega = \Omega_{i,t} \subseteq W_k\), there exists a chain of components
\[
\Omega_{1,t}' := \Omega_{k,t}^*,\quad
\Omega_{2,t}',\quad \cdots,\quad
\Omega_{m,t}' = \Omega
\]
such that
\[
\partial \Omega_{1,t}' \cap \partial \Omega_{2,t}' \neq \varnothing,\quad
\partial \Omega_{2,t}' \cap \partial
\Omega_{3,t}' \neq \varnothing,\quad \cdots,\quad
\partial \Omega_{m-1,t}' \cap \partial
\Omega_{m,t}' \neq \varnothing.
\]
To prove Claim 1, it suffices to show that
\[
\rho_{b_1,\Omega_{i,t}'} \sim \rho_{b_1,\Omega_{i+1,t}'}, \qquad 1\leq  i\leq  m-1.
\]
To see this, note that
by Theorem~\ref{component2}(1),
\(\partial \Omega_{i,t}' \cap \partial \Omega_{i+1,t}'\) contains exactly one critical point, denoted by \(Q'\). Let \(w_0 = B(b_1)\), and let \(\gamma\) be a non-self-intersecting curve connecting \(w_0\) to \(B(Q')\) such that
\[
\gamma(u) \in t\mathbb D  \backslash
 \mathcal{E}_B, \qquad 0\leq  u<1.
\]
By Lemma~\ref{component3}, there are
 two curves, lying in \(\Omega_{i,t}'\)
 and \(\Omega_{i+1,t}'\) respectively,
  connecting a point in \(B^{-1}(w_0)\) to \(Q'\).
   Let \(\lambda\) and \(\mu\) be the endpoints
    of these two curves lying in \(B^{-1}(w_0)\);
     then \(\lambda \overset{B}{\simeq} \mu\). Applying Lemma~\ref{gluel} immediately yields
\[
\rho_{b_1,\lambda} \sim  \rho_{b_1,\mu}.
\]
Combined with the equivalences
$
\rho_{b_1,\lambda} \sim \rho_{b_1,\Omega_{i,t}'}$
 and $
\rho_{b_1,\mu} \sim \rho_{b_1,\Omega_{i+1,t}'},$
 we conclude that
\[
\rho_{b_1,\Omega_{i,t}'} \sim \rho_{b_1,\Omega_{i+1,t}'}.
\]
Therefore \(\rho_{b_1,\Omega_{k,t}^*} \sim \rho_{b_1,\Omega}\) for every \(\Omega\subseteq W_k\). This proves Claim 1.

Let $\lambda $ be a point in $B^{-1}(w_0)\cap W_N$
such that $\lambda\neq b_N$.
The same method
as in the proof of Claim 1 yields the following.

\textbf{Claim 2.} For  $1\leq    i \leq      N-1,$ and
for $z',z''\in B^{-1}(w_0)\cap W_N,$ we have
  $$\rho_{b_i,z'}\sim \rho_{b_i,z''}.$$
Consequently, $ \rho_{b_i,\lambda}\sim \rho_{b_i,b_N}$.

\textbf{Claim 3.}  We have the equivalence chain:
$ \rho_{b_1,\lambda}\sim \rho_{b_2,\lambda}\sim \cdots\sim \rho_{b_N,\lambda}.$

To prove Claim 3 holds, observe $w_0=B(b_1)=\cdots=B(b_N).$
By Lemma \ref{component3}, there exist curves $\gamma_0,\gamma_1,\cdots,\gamma_N$
 such that $\gamma_i(0)=b_i,\gamma_i(1)=* \
 $ for $1\leq    i\leq     N)$,
 $$B(\gamma_0(u))=B(\gamma_1(u))=\cdots=B(\gamma_N(u)), \ \ \
 0\leq    u\leq     1,$$
 and $\gamma_0(0)=\lambda$.
 Let \(\lambda_0(1)=v\). We
  will soon see that \(B'(v)\neq 0\).
   Clearly \(v\neq *\) (otherwise,
   since both \(\gamma_0\) and \(\gamma_N\) lie in \(W_N\),
    they would share the same endpoints and hence
     belong to the same component of \(\{|B|<t\}\).
     But this contradicts Lemma~\ref{component3},
      as \(\lambda_0(0)=\lambda \neq b_N = \lambda_N(0)\)).
       Since \(B(v)=B(*)\) and \(*\) is a critical point, the injectivity of \(B|_{Z(B')}\) forces \(B'(v)\neq 0\), as desired.

By applying Lemma  \ref{gluel},  we get
 $$ \rho_{\lambda, b_1}\sim \rho_{\lambda, b_2}\sim
 \cdots\sim \rho_{\lambda, b_N}.$$
 Remark \ref{remxyyx}(2) then yields
$ \rho_{b_1,\lambda}\sim \rho_{b_2,\lambda}
\sim \cdots\sim \rho_{b_N,\lambda},$
 finishing the proof of Claim 3.

We will use the Claims 1 and 2 to prove the following.

\textbf{Claim 4.}  We have $\rho_{b_1,b_2}\sim \rho_{b_1,b_i}$ for $2 \leq     i \leq     N$.

For this, we first show that  Claim 4 is equivalent
to prove  $$\rho_{b_N,b_1}\sim \rho_{b_N,b_i}, \ \ \ \ 1\leq     i \leq      N-1,$$
 reminding that $b_i\in \Omega_i$ for each  $i$.
  In fact, by Theorem~\ref{e1}(2),
  the points \(b_1,\cdots,b_N\) are precisely those around \(q\) satisfying
\[
B(b_1)=\cdots=B(b_N)=w_0,
\]
and \(\mathrm{order}_{B'}(q)=N-1\). By Theorem~\ref{bott}, after a suitable permutation, we may arrange these points in anticlockwise order; denote them by \(b_1'=b_1,\cdots,b_N'\). Moreover, there exist pairwise
 disjoint non-self-intersecting curves \(\sigma_1,\cdots,\sigma_N\) in \(\Omega_{0,t}\backslash B^{-1}(\mathcal{E}_B)\) such that
\[
B(\sigma_1(u)) = \cdots = B(\sigma_N(u)), \qquad 0\leq  u\leq  1.
\]
with $ \sigma_i(0) = b_i'$ and $\sigma_i(1) = b_{i+1}'$ for $1 \leq     i \leq    N$
(where $b_{N+1}' \equiv b_1'$). Then
 by  Lemma \ref{continuation}, we obtain the following
  equivalence
\begin{equation*} 
\rho_{\sigma_i(0), \sigma_j(0)} \sim \rho_{\sigma_i(1), \sigma_j(1)}, \ \ \ \ 1 \leq  i\leq     j\leq      N.
\end{equation*}
Using this, we obtain the following equivalent conditions:
\begin{itemize}
    \item $\rho_{b_1,b_2}\sim \rho_{b_1,b_i}$ for $2 \leq     i \leq    N$; that is,
     $\rho_{b_1',b_2}\sim \rho_{b_1',b_i}$ for $2 \leq     i \leq    N$;
    \item for each $j~(1\leq     j\leq    N)$ and $m,n\neq j$,
$\rho_{b_j',b_m'}\sim \rho_{b_j',b_n'}$. 
\end{itemize}
Consequently, Claim 4 is equivalent
to    $$\rho_{b_N,b_1}\sim \rho_{b_N,b_i},\ \ \  1\leq     i \leq      N-1.$$

To proceed,  by Claims 2 and 3,
    we immediately have
 $$  \rho_{b_1,b_N}\sim \rho_{b_2,b_N}\sim\cdots\sim \rho_{b_{N-1},b_N}.$$
By Remark \ref{remxyyx}(2),  this gives
 $$  \rho_{b_N, b_1}\sim \rho_{b_N ,b_2}\sim
 \cdots\sim \rho_{b_N, b_{N-1}}.$$
Then $\rho_{b_1,b_2}\sim \rho_{b_1,b_i}$ for $2 \leq     i \leq     N$.
The proof of Claim 4 is complete.

So far, we have completed the proofs
of Claims 1--4. By Claim 4, we have
 $\rho_{b_1,b_2}\sim \rho_{b_1,b_i} $
 for $2 \leq     i \leq    N$.
By (\ref{last0}) in Claim 1,
$$\rho_{b_1,b_2}\sim \cdots\sim \rho_{b_1,b_m} \sim \rho_{b_1,b_k}\sim \rho_{b_1,W_k}$$ for $m+1\leq     k\leq    N.$
Hence, restricted on $\Omega_{0,t^+} $, $B$
admits exactly two equivalence classes of local inverses,
the class containing $ \rho_{b_1,b_2} $, and  the identity map.
   Hence
  $\Omega_{0,t^+} $ is purely gluable.

The proof of Theorem \ref{submain}
  is finished.
\end{proof}

\subsection{Proof of Theorem \ref{main}}

We start with a lemma that follows
 directly from the definition of a thin Blaschke product.

\begin{lem}\label{claimlem} \emph{For  a thin Blaschke product  $\phi$,
  $\{z\in \mathbb{D}:\phi'(z)=0, \phi(z)=a\}$
  is a finite set for every $a$ in $\mathbb{D}$.}
\end{lem}
\begin{proof}
  Suppose, to the contrary,  that there exists
   a point $a$ in $\mathbb{D}$ so that
   the set $\{z\in \mathbb{D}:\phi'(z)=0, \phi(z)=a\}$
  is infinite. Then there exists
  a sequence of distinct points $\{p_n\}_{n=1}^\infty$  such that
   for $n\geq  1 $,
$$\phi'(p_n)=0 \ \ \ \  \mathrm{and}  \ \ \ \  \phi(p_n)=a. $$
 The points  $  p_1,p_1, p_2, p_2, \cdots$
 form a subsequence of zeros  of $\varphi_a \circ \phi$.
Consequently, $\varphi_a \circ \phi$ cannot
be a thin Blaschke product.
      On the other hand,
       since \(\phi\) is a thin Blaschke product, Proposition~\ref{thin2} implies that \(\varphi_a\circ\phi\) is also thin, a contradiction.
      Thus the set is finite for every \(a\in\mathbb D\).
\end{proof}

  We now proceed to the proof of Theorem \ref{main}.

\begin{proof}[\textbf{\emph{Proof of Theorem \ref{main}.}}]
 Suppose $B$ is a thin Blaschke product
 such that $B$ is essentially injective on its critical set.
 That is, there exists a finite subset $F$ of
  $Z(B') $ such that $B$ is injective on $Z(B')\backslash F $.
 To prove that the multiplication operator $M_B$ is irreducible,
by Proposition \ref{thinfin} it suffices to show that $\mathcal{S}_B$  has exactly two components.

 \textbf{Case I.} \emph{First assume that $B$ is injective on the set
 $Z(B')$ of its critical points.
We will show that  $\mathcal{S}_B$  has exactly two components.}

By Proposition \ref{infcritical},
a thin Blaschke product $B$ always has infinitely many critical points. Let
$z_0  \in  \mathbb{D} \backslash B^{-1}( \mathcal{E}_B)$ and
 all points in $B^{-1}(B(z_0))$
are listed as
$$z_0,z_1,\cdots,z_n,\cdots.$$
For a fixed $n,$ there exists an $r\in (0,1)$ such that
 $z_0,\cdots,z_n$ are contained in
$r\mathbb{D}$, and moreover, we require that
$r\mathbb{D}$ contains at least two critical points
 of $B$.   By Proposition \ref{atomz},
 for every $s$ in $(0,1)$, each component $\Omega$ of $\{|B|<s\}$
 is mapped onto $s\mathbb{D}$
 by the proper restriction $B|_{\Omega}.$
 Let $$r'=\max \{|B(z)|: |z|\leq     r\},$$
 and let $\Omega_0$ be the component of $\{|B|<r'\}$
 containing  $z_0,z_1,\cdots,z_n.$
 It is routine to check that $\Omega_0$ is simply connected.
   By Proposition~\ref{sproper}, the restriction \(B|_{\Omega_0}:\Omega_0\to r'\mathbb D\) is, up to composition with conformal maps, a finite Blaschke product. More precisely, there exist a conformal map \(F_1:\Omega_0\to\mathbb D\), a finite Blaschke product \(\varphi\), and a scaling map \(F_2(z)=r'z\) such that
\[
B|_{\Omega_0}=F_2\circ \varphi\circ F_1.
\]
Moreover, the map
\[
\sigma \mapsto F_1\circ \sigma \circ F_1^{-1}
\]
gives a one-to-one correspondence between the local inverses of \(\varphi\) and those of \(B|_{\Omega_0}\).
   By Theorem \ref{submain}, we have
   $$ \rho_{z_0,z_1 }\sim\cdots\sim \rho_{z_0,z_n}.$$
   By the arbitrariness of $n$, we have
   $ \rho_{z_0,z_1}\sim  \rho_{z_0,z_j},\ j\geq   1.$
   Thus $B$
admits exactly two equivalence classes of local inverses,
the class containing $ \rho_{z_0,z_1} $, and  the identity map.
    By Remark \ref{remxyyx}(1), $\mathcal{S}_B$
   has exactly two components.

  \textbf{Case II.}\emph{ In general, we assume
   $B$ is essentially injective on its critical set. That is,
    there exists a finite set $F$ such that  $B|_{Z(B')\backslash F}$ is injective.}

 Let $p_1, p_2, \cdots, p_n, \cdots$
 be all critical points of $B$. Since $B$ is
 essentially injective on its critical set $Z(B'),$
 there exists a positive integer $N_0$ such that
 \begin{equation} \label{edd1}
  B(p_n) \neq B(p_m),
  \quad m>n \geq   N_0 .
  \end{equation}
By Lemma \ref{claimlem}, the set
 $\{n: B(p_n)=B(p_i)\}$ is  finite for each $i $
 with $1\leq   i\leq     N_0 $.
 For all sufficiently large $n$, we thus have
  \begin{equation}\label{edd2}  B(p_n) \neq B(p_i),
  \quad i=1, \cdots, N_0.   \end{equation}
   Combing (\ref{edd1})  and (\ref{edd2}), we may choose an integer
     $N > N_0$
   such that
   \begin{equation} \label{edd3}
   B(p_n) \neq B(p_k), \quad n \neq k, \quad n \geq   N.
   \end{equation}

 Fix  an  arbitrary natural number $n$.
  We can select a closed
    disk
    $\overline{D}\subseteq \mathbb{D}$
    that   contains all critical points
     $p_1, \cdots, p_N$ and
     the finite zero sequence $z_0,\cdots,z_n.$
       Define
$$ t_0 = \sup \,\{|B(z)| : z \in \overline{D}\}. $$
Let $\Omega_{1, t_0}$ denote the component
 of  $\{|B| < t_0\}$ that contains $\overline{D}$;
 in particular, every point $p_1, p_2,\cdots, p_N$
  lies in $\Omega_{1, t_0}$. Choose a critical point $Q^* $
   that lies outside $\Omega_{1, t_0}$.
We may then pick a value  $t^*$ with $t_0 < t^* < 1$ such that there exists
 a component $\Omega$ of $\{|B| < t^*\}$
 containing both $\Omega_{1, t_0}$ and  $Q^*$. 
 Observe that by Lemma \ref{atomz}, the restriction
$B|_{\Omega}$ is essentially  a finite
 Blaschke product, with $\Omega$ playing the role
 of the unit disk $\mathbb{D}$.
 Techniques for finite Blaschke
  products remain applicable for $B|_{\Omega}$. The
 following arguments, based on this restricted map,
 will  prove that $\Omega$  is purely gluable.

Consider the family of the components $\Omega_{j, t_0}$ of
$\{z\in \Omega: |B(z)| < t_0\}$
indexed by $j$.
By slightly increasing \(t_0\) while keeping \(t_0<t^*\), we may assume that \(\partial \Omega_{1,t_0}\) intersects at least one other boundary \(\partial \Omega_{k,t_0}\) of a component lying in \(\Omega\). Such a \(t_0\) exists; otherwise, for every \(s<t^*\), the number \(n_s\) of components of \(\{z\in\Omega: |B(z)|<s\}\) would satisfy \(n_s\geq2\). However, by Proposition \ref{sproper}  and
Lemma \ref{entirecompo}, we have
 \(n_s=1\) for \(s\) sufficiently close to \(t^*\), which yields a contradiction. Hence, such a  \(t_0\) must exist.


We then take the maximal sub-collection of these boundaries   $\partial \Omega_{j, t_0}$ whose union is connected and contains  $\partial \Omega_{1, t_0}$.
  By the finite tree structure of these boundaries $\partial \Omega_{j, t_0}$,
 it follows that there exist at least two
    $R$-components.
  Let $\Omega_{0, t_0}$ be such a $R$-component
  distinct from $\Omega_{1, t_0}$ and remind that
  $\Omega_{0, t_0} \subseteq \Omega $ (see Figure \ref{root}).
  We select $\Omega_{0, t_0}  $ as the $R$-component
   of maximal multiplicity among all $R$-components distinct from
 $\Omega_{1, t_0} $. Formally,
  if $\Omega_{i_1, t_0} ,\cdots,\Omega_{i_K, t_0} $
 denote all such $R$-components, then
 $\Omega_{0, t_0}  $ is chosen as any one of them
 that attains the maximal multiplicity.
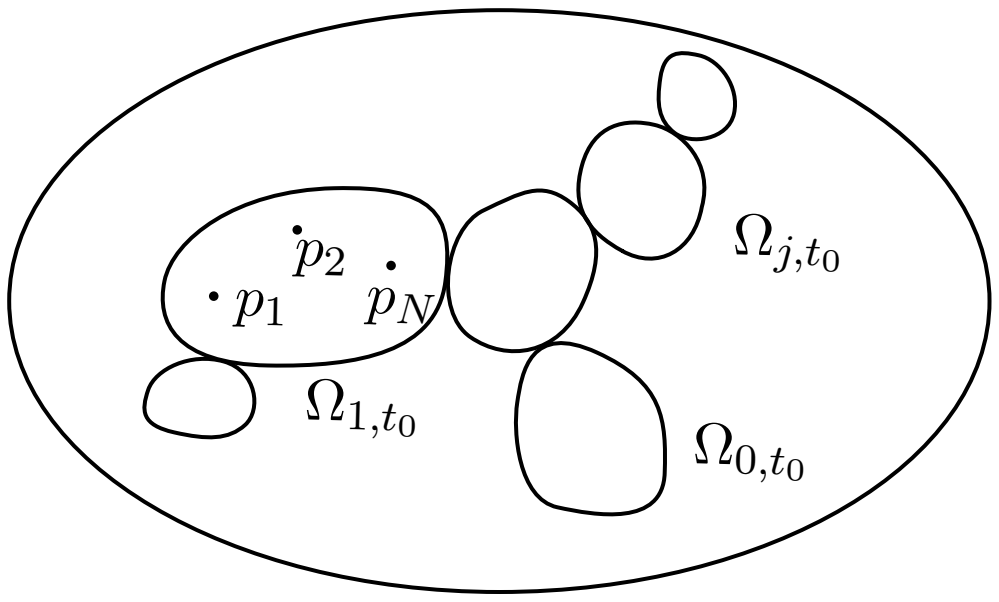
\begin{figure}[htbp]
\centering
\tikzset{every picture/.style={line width=0.75pt}} 
\begin{tikzpicture}[x=0.75pt,y=0.75pt,yscale=-1,xscale=1]
\draw   (333.15,185) .. controls (327.09,198.28) and (314.1,199.94) .. (304.91,193.62) .. controls (295.71,187.29) and (298.57,166.48) .. (308.16,161.88) .. controls (317.74,157.28) and (322.67,154.5) .. (329.57,159.53) .. controls (336.47,164.57) and (339.22,171.72) .. (333.15,185) -- cycle ;
\draw  [line width=3] [line join = round][line cap = round] (240.53,183.36) .. controls (240.53,183.36) and (240.53,183.36) .. (240.53,183.36) ;
 \draw  [line width=3] [line join = round][line cap = round] (261.49,166.73) .. controls (261.49,166.73) and (261.49,166.73) .. (261.49,166.73) ;
\draw   (248.4,215.68) .. controls (244.25,219.43) and (238.91,219.43) .. (230.62,217.55) .. controls (222.32,215.68) and (222.3,212.1) .. (224.09,206.86) .. controls (225.88,201.63) and (234.77,197.41) .. (243.07,199.76) .. controls (251.37,202.1) and (252.55,211.94) .. (248.4,215.68) -- cycle ;
\draw   (298.64,179.82) .. controls (296.31,197.31) and (277.21,201.2) .. (254.07,200.76) .. controls (230.93,200.31) and (226.05,189.92) .. (228.17,179.56) .. controls (230.28,169.2) and (245.93,156.84) .. (270.72,156.28) .. controls (295.52,155.73) and (300.97,162.32) .. (298.64,179.82) -- cycle ;
\draw   (356.16,142.39) .. controls (352.12,139.29) and (351.6,135.1) .. (352.61,128.51) .. controls (353.62,121.92) and (357.09,121.76) .. (362.35,122.95) .. controls (367.6,124.14) and (372.55,130.95) .. (371.08,137.56) .. controls (369.62,144.17) and (360.2,145.5) .. (356.16,142.39) -- cycle ;
\draw   (343.19,172.09) .. controls (334.17,167.21) and (330.06,160.89) .. (332.94,150.51) .. controls (335.82,140.13) and (343.89,138.69) .. (351.16,140.39) .. controls (358,142) and (366,150) .. (363,161) .. controls (361,172) and (352,177) .. (343,172) -- cycle ; 
\draw   (189.2,184.6) .. controls (189.2,144.28) and (244.27,111.6) .. (312.2,111.6) .. controls (380.13,111.6) and (435.2,144.28) .. (435.2,184.6) .. controls (435.2,224.92) and (380.13,257.6) .. (312.2,257.6) .. controls (244.27,257.6) and (189.2,224.92) .. (189.2,184.6) -- cycle ;
\draw   (353.69,227.74) .. controls (353.26,241.22) and (336.48,238.47) .. (326.75,236.34) .. controls (317.03,234.21) and (315.56,217.69) .. (316.69,210.06) .. controls (317.82,202.42) and (320,189.88) .. (336.71,197.52) .. controls (353.42,205.15) and (354.12,214.27) .. (353.69,227.74) -- cycle ;
\draw  [line width=3] [line join = round][line cap = round] (285.03,175.68) .. controls (285.03,175.68) and (285.03,175.68) .. (285.03,175.68) ;
\draw (369.41,161.35) node [anchor=north west][inner sep=0.75pt]  [rotate=-1.24] [align=left] {$\Omega_{j,t_0}$};
\draw (316.76,265.78) node [anchor=north west][inner sep=0.75pt]  [rotate=-1.28] [align=left] {$\Omega$};
\draw (359.45,213.76) node [anchor=north west][inner sep=0.75pt]  [rotate=-1.52] [align=left] {$\Omega_{0,t_0}$};
\draw (244.41,180.35) node [anchor=north west][inner sep=0.75pt]  [rotate=-1.24] [align=left] {$p_1$};
\draw (277.45,179.76) node [anchor=north west][inner sep=0.75pt]  [rotate=-1.52] [align=left] {$p_N$};
\draw (262,202.76) node [anchor=north west][inner sep=0.75pt]  [rotate=-1.52] [align=left] {$\Omega_{1,t_0}$};
\draw (259.45,167.76) node [anchor=north west][inner sep=0.75pt]  [rotate=-1.52] [align=left] {$p_2$};
\fill[black] (240.53,183.36) circle (0.9pt); 
\fill[black] (261.49,166.73) circle (0.9pt); 
\fill[black] (285.03,175.68) circle (0.9pt); 

\end{tikzpicture}
\caption{An $R$-component $\Omega_{0,t_0}$
 distinct from $\Omega_{1,t_0}$}
\label{root}
\end{figure}
We now state a key claim.

\noindent \textbf{Claim.}
\emph{The domain $\Omega_{0, t_0^+}  $  is purely gluable.
}

Recall that the notation $ t_0^+ $ refers to
 a number $s>t_0$ taken sufficiently close
to $t_0$. The claim can be restated formally as follows:
  there exists some $\varepsilon>0,$
  such that for every $s$ with $t_0<s<t_0+\varepsilon,$
the component $\Omega_{0,s}$ of $\{|B|<s\}$ that contains
 $\Omega_{0,t_0}$ is
  purely gluable.

Remind that  $\Omega_{0, t_0^+}  $ contains
 the tree structure associated with
 $\overline{\Omega}_{0, t_0 }  $ (or, with $\overline{\Omega}_{1, t_0 }$)  .
We treat some distinct  cases
for $\Omega_{0, t_0}  $,  corresponding respectively
 to Step I and Step II in the proof of Theorem \ref{submain},
and verify the claim in each case.

 \textbf{Step I.} Flower model.

  \textbf{Step II.} The Cactuses Case.

There are exactly three cases under consideration:
 $\Omega_{0, t_0}$ has at least two   critical points;
 $\Omega_{0, t_0}$ has exactly one   critical point; and
 $\Omega_{0, t_0}$ has no critical point.
  These cases are respectively
  addressed in Steps I  and II
  of the proof of Theorem \ref{main}.

  We now point out the  difference between
   the current general setting
   and the earlier special case (Case I of Theorem \ref{main}).
The key observation is that under Condition (\ref{edd3}),
 for $\lambda, \mu \in \Omega_{1, t}$
and $\omega \in \Omega_{0, t}$ with $B(\omega) =
B(\lambda) = B(\mu)$, Lemma  \ref{componenthigh} yields
$\rho_{\omega, \lambda} \sim  \rho_{\omega, \mu}$.
 In other words,
 we only need to substitute  Corollary \ref{component}
 with Lemma  \ref{componenthigh}. Then
the arguments presented in Steps I and II of Theorem \ref{submain}
 proceeds accordingly,  ensuring the desired conclusion.
 That is, $\Omega_{0,t_0^+}$ is
  purely gluable. The proof of the claim is compete.

  For any other $\Omega_{j,t_0}$ with closure
  lying outside the component of  $\{|B|\leq   t_0\}$
  containing $\overline{\Omega}_{1,t_0}$,
    the arguments in Steps I and II  of Theorem \ref{submain}
     apply without modification.
     Therefore, for such $\Omega_{j,t_0}$, the corresponding domain
     $\Omega_{j,t_0^+}$ either  contains at most
     one critical point or  is purely gluable.
    In conclusion, the induction step at $s=t_0^+$ holds for
    all $\Omega_{j,s}$ lying in $\Omega.$

Following the same procedure as in Steps I and II
 in the proof of Theorem \ref{submain}, we obtain a sequence
     $t_0<t_1<t_2<\cdots<t_d< t^*$,
  where for $i\geq1$ each $t_i $
  denotes  the absolute value of the critical value of
  $B|_{\Omega}$ that is strictly larger than $t_0$.
  This yields a finite collection of stages:
  $$t_0 < s \leq     t_1, \ \ \cdots,
 \ \ t_{d-1} < s \leq    t_d, \ \  t_d <s\leq    t^*.$$
  At each stage,
   the arguments in Steps I and II of Theorem \ref{submain}
     remain applicable with no modification.
     This holds because,
     under condition (\ref{edd3}),
     the problem reduces to the case where
 $B$ is injective on the remaining critical set
 as $\Omega_{0,t_0^+}$ has   absorbed
     the critical points $ p_1,\cdots,p_N $.
  At  $ s=t^* $,
    all components merge into a single
     component $  \Omega$,
    and the final simply connected domain
    $\Omega$ is purely gluable.
    In particular, since $z_0,\cdots,z_n$ lie in $ \Omega $,
 $\rho_{z_0,z_1}\sim \rho_{z_0,z_i},~ 1\leq    i \leq     n. $
By the arbitrariness of $n,$
$\rho_{z_0,z_1}\sim \rho_{z_0,z_i}, i\geq   1$.
This fact in turn implies that the Riemann surface $\mathcal{S}_B$
contains exactly two components.

 The proof of Theorem \ref{main} is now complete.
 \end{proof}

\section{Primeness and essential injectiveness}\label{sect4}
In this section, we show that   the
   essential injectivity of a thin Blaschke
    product implies its primeness, and also
    establishes the implications
  $(3) \Rightarrow (2) \Rightarrow (1)$ in
     Conjecture~\ref{Conj1}.

Remind that in the introduction, we mention
 the equivalence of the following conditions (I) and (II)
 for a thin Blaschke product $B$.

  (I) The mapping $B$ is not essentially injective on its critical set.

(II) There exist two  sequences $\{z_n\}_{n=1}^\infty$ and
  $\{w_n\}_{n=1}^\infty $ in $\mathbb{D}$ such that
 \begin{itemize}
\item[(1)] $B(z_n) = B(w_n)  $, and $z_n \neq w_n$ for all $n$;
\item[(2)] $B'(z_n) = B'(w_n) = 0;$
\item[(3)] $B(z_n)\neq B(z_m)$ if $m\neq n$.
\end{itemize}
\begin{prop} \label{EIII}
The equivalence between \emph{(I)} and \emph{(II)} holds for a thin Blaschke product $B$.
\end{prop} \begin{proof}
 First, we prove that (II) implies (I). By (II),
  there exist two sequences \(\{z_n\}_{n=1}^\infty\) and \(\{w_n\}_{n=1}^\infty\)
   such that
\[
B'(z_n)=B'(w_n)=0 \quad\text{and}\quad B(z_n)=B(w_n).
\]
  Assume,
for contradiction, that \(B\)
 is essentially injective on its critical set.
 Then there exists a finite set \(\mathcal{F}_0\) such that \(B\)
  is injective on \(Z(B')\backslash  \mathcal{F}_0\).
  In particular, \(B\) is injective on
\[
\big(\{z_n:n\geq 1\}\cup\{w_n:n\geq 1\}\big) \backslash  \mathcal{F}_0,
\]
which is impossible. This contradiction
 proves that (II) implies (I).

 Next we show that  (I) implies (II).

  By (I), for any finite subset
  $\mathcal{F}$, $B$ is not injective on
   $Z(B')\backslash \mathcal{F}$. Also recall,
   by Proposition \ref{infcritical}, that $Z(B')\backslash \mathcal{F}$ is infinite.
    Specifically,
   when $\mathcal{F} = \emptyset$, there exist
   distinct points $z_1, w_1 \in \mathbb{D}$ satisfying
    $$B(z_1) = B(w_1) \quad \text{and} \quad B'(z_1) = B'(w_1) = 0.$$
    By induction,   we  have constructed sequences ${z_n}$ and ${w_n}$
     $(1 \leq    n \leq     N)$ in $\mathbb{D}$ with the following
      properties for $1 \leq     m, n \leq     N$:
      \begin{itemize}
      \item[(1)] $B(z_n) = B(w_n)$ with $z_n \neq w_n$;
      \item[(2)] $B'(z_n) = B'(w_n) = 0$;
      \item[(3)] $B(z_n) \neq B(z_m)$ whenever $m \neq n$.
      \end{itemize}
By Lemma \ref{claimlem},
 for each $n \geq   1$, the set
      $$\{z \in \mathbb{D} : B'(z) = 0, B(z) = B(z_n)\}$$
       is finite. Let
       $$\mathcal{F} = \{z \in \mathbb{D} : B'(z) = 0, B(z) = B(z_n)
        \text{ for some } 1 \leq     n \leq     N \},$$
     and $\mathcal{F}$ is a finite set.    Since $B$ is  not injective on
          $Z(B')\backslash \mathcal{F}$,
          there exist critical points $z^*, w^* $ in $Z(B')
          \backslash \mathcal{F}$ satisfying
          $$B(z^*) = B(w^*) \quad \text{and} \quad
          B'(z^*) = B'(w^*) = 0.$$ Setting $z_{N+1} = z^*$
          and $w_{N+1} = w^*$
          extends the sequences to satisfy all conditions
         of (II) for
          $n = 1, \cdots, N+1$, which
 completes the inductive step. This concludes
 the proof.
\end{proof}
The essential injectivity of a
thin Blaschke product implies its primeness.
\begin{prop} \label{onlyif}
Let $B$ be a thin Blaschke product.
If $B$ is essentially injective on
its critical points, then $B$ is prime.
\end{prop}
\begin{proof}
 Suppose $B$ is essentially injective on its critical set.
To establish that $B$ is prime, assume for contradiction
  that $B$ is not prime. Let $B=\eta \circ  u$,
  where $\eta$ and $u$ are inner functions and neither
  $\eta$ nor $u$ is a M\"obius map.  Then by \cite[Theorem 4.1]{GLMR94}, $u$ is a finite Blaschke
   product of order $N \geq   2$.
Since $B(z)=\eta \circ u(z),$
\begin{equation} \label{et2}
B'(z)=\eta'(u(z)) \cdot u'(z).
\end{equation}
Since $u'$ has only finitely many critical points
 and, by Theorem \ref{infcritical}, $B$ has
  infinitely many critical points,
the   function $\eta$ must have infinitely many
 critical points, say,
$\lambda_1, \cdots, \lambda_n, \cdots$.
 By (\ref{et2}),
for all but finitely many points
in $  \mathcal{E}_u $, each $\lambda_j$
has exactly $N$ preimages in $u^{-1}(\lambda_j)$.  
 Without loss of generality, for  $j\geq   1$  we can pick
 $z_j^{(1)}, z_j^{(2)} \in u^{-1}(\lambda_j) $  such that
  $ z_j^{(1)}\neq  z_j^{(2)}$
  and
\begin{eqnarray*}
  & B^{\prime} (z_j^{(1)})=B^{\prime} (z_j^{(2)})=0,\\ 
& B (z_j^{(1)})=B (z_j^{(2)})=\eta (\lambda_j),\ \ \ \  j\geq   1.
\end{eqnarray*}
Since $B$ is a thin Blaschke product,  using
 Lemma \ref{claimlem},  we obtain that
 $ \{ \eta (\lambda_j ): j \geq    1 \}$ is
 an infinite sequence.
  Passing to a subsequence of $\{\lambda_j\}_{j=1}^\infty$, we observe that
  $B$ is not essentially injective on
its critical points. This  contradicts our assumption,
 thus completing  the proof.
\end{proof}

To conclude this subsection, we clarify the implications in  Conjecture \ref{Conj1}  stated below.
\begin{prop}\label{eqq}
In Conjecture \ref{Conj1}, the implications \((3)\Rightarrow(2)\Rightarrow(1)\) hold.
\end{prop}
\begin{proof}
Let $B$ be a thin Blaschke product.
  We will  prove
the   implication chain
  $(3) \Rightarrow (2) \Rightarrow (1) $.

 $(3) \Rightarrow (2)$.  By Theorem 5.3.1 in
\cite{GH4}, for a thin Blaschke
product $B$, $\mathcal{V}^*(B)=\mathbb{C}I $  if and only if
 the Riemann surface
$\mathcal{S}_B$ has a component $\mathcal{S}$
of finite multiplicity, other than
  the identity component $\mathcal{S}_0$.
Given that $\mathcal{S}_B$ has exactly
two components, the   component $\mathcal{S}$
 distinct from $\mathcal{S}_0$  
  must
be of infinite multiplicity. Consequently,
$\mathcal{V}^*(B)=\mathbb{C}I.$

 $(2) \Rightarrow (1)$.
Suppose  $\mathcal{V}^*(B)$ is trivial. We will prove  that $B$ is  prime.
Assume  conversely that $B$ is not prime.
By  \cite[Theorem 4.1]{GLMR94},   if $B=\eta \circ B_0$
where $\eta$ and $B_0$ are inner functions distinct
 from M$\mathrm{\ddot{o}}$bius maps, then $B_0$ must be a finite
  Blaschke product. Moreover, $B_0$   is not a M\"obius map
  by non-primeness of $B$. Clearly,
   we have $\mathcal{V}^*(B) \supseteq \mathcal{V}^* (B_0 )$.
  By \cite{HSXY}, $V^*(B_0) \neq \mathbb{C} I$.
  Thus $\mathcal{V}^*(B) \neq \mathbb{C} I$, which is a contradiction.
\end{proof}

\section{An example tightening Conjecture \ref{Conj1} }\label{sect5}

The following example constructs thin Blaschke products
with prescribed properties, and illustrates how
 to extract subsequences from thin sequences.

\begin{exam} \label{as}
Let \(\{z_k\}_{k=1}^\infty\) be a  Blaschke sequence.
As follows, we will extract a subsequence  $\{\lambda_k\}_{k=1}^\infty$ of
\(\{z_k\}_{k=1}^\infty\)  which
  forms a thin Blaschke
  Blaschke product  $B$ such that for $n\geq   1$,
  the multiplication operator $M_{B^n}$ is
irreducible, and the associated Riemann surface  $\mathcal{S}_{B^n}$
 has precisely $n+1$ components.
\end{exam}
\noindent To define \(\mathcal{S}_{B^n}\),
note that \(B^n\) is not a thin Blaschke product
for \(n\ge 2\). However, by routine computations,
we see that \((B^n)^{-1}(\mathcal{E}_{B^n})\)
is still a discrete subset of \(\mathbb{D}\).
 Hence \(\mathcal{S}_{B^n}\) is well-defined.

This example shows that Conjecture~\ref{Conj1},
which concerns a single thin Blaschke product,
does not extend to products
 of finitely many thin Blaschke products.

To demonstrate Example \ref{as},
we need a lemma from \cite[Lemma 5.1]{GH3}.
For $a\in \mathbb{D}$, let $$\psi_a(z)=
 \frac{|a|}{a} \frac{a-z}{1-\overline{a}z},\ \ \ \  z\in \mathbb{D}.$$
\begin{lem} \label{51gh}
    Suppose that $B$ is a finite Blaschke product
with order greater than $2$. Then for any $r \in (0,1)$,
there exists $s \in (0,1)$ such that for
 all $\lambda \in \mathbb{D}$ with $|\lambda| > s$,
 the function $B\psi_\lambda$ has a critical point
 $w^*$ satisfying $|B(w^*)\psi_\lambda(w^*)| > r$.
 \end{lem}

\noindent  Specifically, for $\lambda$ with $|\lambda|$
 sufficiently close to $1$,
 $\psi_\lambda$ is uniformly close to
 a unimodular constant on each compact subset of $\mathbb{D}$.
  By an application of Rouch$\acute{\mathrm{e}}$'s theorem,
 the  critical points of $B  \psi_\lambda$  other than $w^*$
  lie in a small neighborhood of critical points of $B$. The ``new"
  critical point \(w^*\) can attain values
   \(\lambda^*\) with \(|\lambda^*|\)
   arbitrarily close to \(1\).

\noindent \textbf{The construction of $B$.}
The construction adopts an approach very similar to
the technique employed in the proof
of \cite[Theorem 5.3]{GH3}. For clarity, we present the key idea rather than the detailed steps.

  Let \(\{z_k\}_{k=1}^\infty\) be an arbitrary Blaschke sequence.
  Without loss of generality, we may assume that \(\{z_k\}_{k=1}^\infty\) is a thin Blaschke sequence;
   otherwise, we extract a subsequence to replace  \(\{z_k\}_{k=1}^\infty\).
    We begin by picking $\lambda_1$
   and $\lambda_2$ from the sequence \(\{z_k\}_{k=1}^\infty\)
   so that the order-two Blaschke product
     \[ B_1 =  \psi_{\lambda_1}\psi_{\lambda_2}  \]
 has exactly one    critical point distinct from zero.
 The key step involves constructing a sequence of finite
 Blaschke products $\{B_n\}_{n=1}^\infty$ through the recursive relation:
\[ B_{k+1} = B_k \cdot \psi_{\lambda_k} \]
where  $\{\lambda_k\}_{k=1}^\infty$ is a subsequence of the thin sequence $\{z_k\}_{k=1}^\infty$.
 Hence the limit of  $\{B_k\}_{k=1}^\infty$ exists
 as a thin Blaschke product.
The sequence $\{B_m\}_{m=1}^\infty$ is
accompanied by an increasing sequence
$\{K_m\}_{m=1}^\infty$ of closed domains in
$\mathbb{D}$:
\[ K_1 \subseteq K_2 \subseteq K_3 \subseteq
\cdots \quad \text{with} \quad \bigcup_{m=1}^\infty K_m = \mathbb{D}, \]
and for each $m,$ $K_m$ is contained in the interior of $K_{m+1} $;
 it is required that the interior $K_{m+1}$ is required to contain
 the first $m+1$ critical points of $B_{m+1},B_{m+2},\cdots ,$
  and of  the limit $B.$
For each $m,$ the limit product $B$ can be written as
$$B=B_m \widetilde{B}_m,$$
for an infinite Blaschke product $\widetilde{B}_m$.
 When \(\{|\lambda_k|\}_{k=1}^\infty\) tends rapidly to \(1\),
  the quantity \(\|\widetilde{B}_m-1\|_{K_m,\infty}\) becomes sufficiently small, and
\[
\|B-B_m\|_{K_m,\infty}\leq  \|\widetilde{B}_m-1\|_{K_m,\infty}.
\]
 This error $\|B-B_m \|_{K_m,\infty} $ can be controlled so that
 on $K_m$, $B$ has exactly the same number of critical points
 of $B_m,$ counting multiplicity. By induction, we can ensure
 that $B_m$   admits only simple critical points;
 thus passing to the limit, $B$ retains the same property.

 By rearranging the first $m$ critical points
  of \(B_i\) \((i>m)\) in \(K_m\),
  we can ensure that they lie in a sufficiently small
  neighborhood of the critical points of \(B_m\).
  Let $p_{m, 1}, \cdots, p_{m, m}$  be all
  critical points of $B_m$, and let
   $p_{i, 1}, \cdots, p_{i, m}$ be the only critical points of $B_i$
   on $K_m.$ By Lemma \ref{51gh},
   we  carefully choose $\{\lambda_k\}_{k=1}^\infty$
so that
$$
|B_m(p_{m, 1})|<
|B_m(p_{m, 2})|<\cdots <\mid B_m (p_{m, m}
)|  .
$$
Then with fixed $m$, we can find constants
$c_i'$ and $c_i''$ satisfying
\begin{align*}
&|B_m(p_{m, 1})|<c_1'<c_1''<
|B_m(p_{m, 2})|\\
&\ \ \ <c_2'<c_2''<\cdots  <c_{m-1}'<c_{m-1}''\\
&\ \ \ <| B_m (p_{m, m}
)|.
\end{align*}
The above inequalities holds if each $p_{m, i}~(1 \leq     i \leq   m)$
is replaced with a sufficiently small   closed disk centered at $p_{m, i}$
and entirely contained in  $K_m$.
Furthermore, for $i>m,$ the above inequality is maintained:
\begin{align*}
&|B_i(p_{i, 1})|<c_1'<c_1''<
|B_i(p_{i, 2})|\\
& \ \ \ <c_2'<c_2''<\cdots  <c_{m-1}'<c_{m-1}''\\
&\ \ \ <| B_i (p_{i, m})|,
\end{align*}
where each $\ p_{i, j} $ lie in the above closed disk centered at
$ p_{m, j} $ for $1\leq     j \leq     m$.
By taking a limit, we have
$$
|B(p_1)|<|B(p_2)|<\cdots <|B(p_m)|,
$$
where $p_1,p_2,\cdots, p_m$ denotes the only critical points of $B$ on $K_m.$
In such a way, the critical points of $B$
can be listed as $p_1,p_2,\cdots   $ satisfying
\begin{equation}
\label{gekzero}
0<|B(p_1)|<|B(p_2)|<\cdots  .
\end{equation}
Consequently, $B$ is injective on its critical set $Z(B')$.

\vskip2mm

\noindent  \textbf{Properties of $B$.} Next we will show that $B$
 satisfies the required property by first analyzing
 $\mathcal{S}_{B^n}$ and then establishing
 its connection to the irreducibility
  of $M_{B^n}$.

  It is clear that $B$ is injective on critical points.
Let $\rho$ be a local inverse of $B^n$ defined on a subdomain of
 $\mathbb{D}\backslash \mathcal{E}$ for some discrete set $\mathcal{E}$.
   To illustrate, we assume
$n=2$ without loss of generality. Now for a local inverse of
$B^2,$ consider $$B^2(\rho(z))=B^2(z),$$
which splits into $B (\rho(z))=B(z) ,$ and $B (\rho(z))=-B(z).$
 By Theorem \ref{main}, the local inverses of \(B\)
(the case \(B(\rho(z))=B(z)\)) fall into
 two equivalence classes, corresponding to the two components
  of \(\mathcal S_B\): the identity component
 \(\mathcal S_0\) and its complement \(\mathcal S_1\).

 Now consider
 $$B (\rho(z))=-B(z) .$$
In this case, we can establish an analogue of Lemma \ref{gluel}.
 \begin{lem} \label{neew}
Let $B$ be a thin Blaschke product.
Given $z_0,z_1$ and $z_2$ in $\mathbb{D} \backslash B^{-1} (\mathcal{E}_B)$
with $- B(z_0)=B(z_1)=B(z_2)$, suppose there are three curves $\sigma_0,\sigma_1,\sigma_2$
such that
\begin{itemize}
\item[\emph{(1)}]  $\sigma_i(0)=z_i$ for $i=0,1,2,$ and $\sigma_1(1)=\sigma_2(1)$;
\item[\emph{(2)}]  $-B(\sigma_0(t))=B(\sigma_1(t))=B(\sigma_2(t))$ for $0\leq    t\leq     1;$
\item[\emph{(3)}] $ \sigma_i(t) \in  \mathbb{D} \backslash B^{-1} (\mathcal{E}_B) $,
 $i=0,1,2  $ for $0\leq    t< 1;$
\item[\emph{(4)}]   $B'(\sigma_0(1))\neq 0.$
\end{itemize}
Then as local inverses of $B^2,$ $\rho_{z_0,z_1} \sim \rho_{z_0,z_2}$.
 \end{lem}
 \noindent We now explain the existence of
  \(\sigma_0,\sigma_1,\sigma_2\) in our setting.
  Here \(z_1\) and \(z_2\)
  are  directly \(B\)-glued in a component \(\Omega\)
  of \(\{|B|<t^*\}\), that is,
$
z_1 \overset{B}{\backsimeq} z_2.
$
 Then there exist \(\sigma_1,\sigma_2\) satisfying
the conditions in Lemma~\ref{neew}. In particular,
\[
\sigma_i(u)\in\mathbb D\backslash B^{-1}(\mathcal E_B), \ \ \ \  0\leq  u <1;
\ \ i=1,2 .
\]
Since \(B|_{\Omega}\) is a holomorphic proper map by Proposition
\ref{atomz},
we can find \(\sigma_0\) such that
\[
B(\sigma_0(u))=B(\sigma_1(u)),
\ \ \ \  0\leq  u\leq  1.
\]
With a small perturbation if necessary, we may require
\[
\sigma_0(u)\in\mathbb D\backslash B^{-1}(\mathcal E_B), \ \ \ \ 0\leq  u<1.
\]

Let $B$ be the thin Blaschke product constructed as above, and
we will soon see that the condition
$B'(\sigma_0(1))\neq 0 $ itself holds. To reach a contradiction,
 assume $B'(\sigma_0(1))= 0 $. Either $\sigma_1(1)$ and $\sigma_0(1)$
are equal or not. If $\sigma_1(1)=\sigma_0(1)$,
 the identity
\[
-B(\sigma_0(u))=B(\sigma_1(u))=B(\sigma_2(u)),
\qquad  0\le u\le 1.
\]
at \(1\) gives \(-B(\sigma_0(1)) = B(\sigma_0(1))=0,\)
a contradiction to
\eqref{gekzero}.
If \(\sigma_1(1)\) and \(\sigma_0(1)\) are distinct critical points,
the identity
\[
-B(\sigma_0(u))=B(\sigma_1(u))=B(\sigma_2(u)),
\qquad  0\le u\le 1.
\]
at \(1\) gives \(|B(\sigma_0(1))|=|B(\sigma_1(1))|\).
Consequently, this leads to a contradiction with (\ref{gekzero}).
For each $0<t<1$, $$\{|B^2|<t^2\}=\{|B |<t \}.$$
Let $\Omega$ be a component of $\{|B |<t \}$ and
 $w_0\in \Omega$. Let
$ w_1,w_2,\cdots,w_N $ be all points in $\Omega$
such that
   $$ -B(w_0)=B(w_1)=\cdots=B(w_N).$$
Note that $B|_{\Omega}$   behaves as a finite Blashke product
 (see Proposition \ref{atomz}).
 By  gluable property of  finite Blaschke product  combined with Lemma \ref{neew},
    we obtain the equivalence relations:
     $$ \rho_{w_0, w_1} \sim \rho_{w_0,w_2} \sim \cdots\sim
       \rho_{w_0,w_N}.$$
As the parameter $t$ increases, the domain $\Omega$ expands,
yielding   an infinite sequence $\{w_k:k\geq   1\}$
satisfying  $B (w_k)= -B(w_0) .$
Consequently, we derive the extended equivalence:
$$  \rho_{w_0, w_1}\sim   \rho_{w_0, w_i}, \ \ \  \ i\geq  1.$$
In other words, the set
\[
\{(z,w): B(z)=-B(w),\ z\in\mathbb D \backslash \mathcal E\}
\]
is connected as a component of \(\mathcal S_{B^2}\),
 denoted by \(\mathcal S_2\). Moreover, there exist points \(\{w_i'\}_{i=1}^\infty\subset\mathbb D\) such that
\[
-B(w_0)=-B(w_1')=\cdots=-B(w_N')=\cdots,
\]
and no other points satisfy this relation. That is,
$$
B(w_0)=B(w_1')=\cdots=B(w_N')=\cdots.
$$
By the construction of \(B\), \(B\) is essentially injective on its critical point set. Applying Theorem~\ref{main}, we obtain
\[
\rho_{w_0,w_1'}\sim\rho_{w_0,w_i'},\qquad i\geq 1,
\]
which are local inverses of \(B\) (and hence of \(B^2\)).
 This corresponds to the only two components of \(\mathcal S_B\)
  mentioned before,  the
 identity component \(\mathcal S_0\) and \(\mathcal S_1\).

It remains to discuss irreducibility of $ M_{B^2} ,$ which is equivalent
to $$\mathcal{V}^*(M_{B^2}) =\mathbb{C}I,$$ where $\mathcal{V}^*(M_{B^2})=
\{M_{B^2},M_{B^2}^*\}'$.
For a thin Blaschke product $\phi,$
each operator $T\in  \mathcal{V}^*(M_{\phi})$
admits the representation
 $$ Tf(z) =\sum_{k=0}^\infty c_k f(\rho_k(z)) \rho_k'(z) $$
 on $\mathbb{D}$ minus a discrete subset. Here,   the
coefficients sequence
 $\{c_k\}_{k=1}^\infty\in l^2 $ are uniquely determined,
  $\rho_k$ are local inverses of $\phi.$
 When $\phi$ is replaced by  $B^2,$ the same reasoning applies.
  Following the approach in \cite{GH3} or \cite[Theorem 5.3.1]{GH4}  $\mathcal{V}^*(M_{B^2})\neq \mathbb{C}I$
 if and only if there is a local inverse $\rho$ of finite multiplicity,
  if and only if $\mathcal{S}_{B^2}$ has a non-identity component
  of finite multiplicity.
  From our earlier analysis,
   $\mathcal{V}^*(M_{B^2})= \mathbb{C}I$.
   Equivalently, $M_{B^2} $ is
   irreducible. The same argument extends to $M_{B^n}.$

The construction of Example \ref{as} is finished. $\hfill \square$

\section{A prime thin Blaschke product failing
 essential injectivity} \label{6cexample} 
 This section will present a thin
 Blaschke product \(B\) with the following properties:
\begin{enumerate}
  \item[(1)]  All zeros of \(B\) tend to \(1\).
  With further analysis, this implies that \(B\) is prime.
  \item [(2)] \(B\) has infinitely many pairs of distinct critical points with the same image, and hence \(B\) is not essentially injective on its critical set.
\end{enumerate}

This construction yields the following result.
 \begin{thm} \label{plusexample}
 There exists a prime thin Blaschke product $B$
 that is not essentially injective on its critical set.
 \end{thm}

 Throughout this section, the notations \(\eta\),
 \(\lambda\), $r$, $ s$ and $P$, 
  are used with meanings independent of their usage in earlier sections.
\subsection{Some preparations}\label{orpre}

\textbf{1. The Cayley transform}

Define
\[
\eta(z)=\frac{1+z}{1-z}.
\]
This is a biholomorphic map from \( \mathbb{D}\)
onto the right half-plane
\[
 \mathbb{H}=\{\eta\in\mathbb C: \mathrm{Re}\,\eta>0\}.
\]
Its inverse is
\[
z=h(\eta)=\frac{\eta-1}{\eta+1}.
\]
For \(A>1\), set
\[
a=h(A)=\frac{A-1}{A+1}\in(0,1).
\]
Then \(A\to\infty\) if and only if \(h(A)\to1\).
Note that $a= \overline{a}$.
  The corresponding Blaschke factor then satisfies
\begin{equation}
\frac{a -z}{1-az}
=
\frac{A-\eta}{A+\eta}.
 \end{equation}
Thus zeros on the interval \((0,1)\) may be parametrized in the right half-plane by positive real numbers \(A\).

To construct a thin Blaschke product on $\mathbb{D}$
with desired properties, it suffices to construct one
on $\mathbb{H}.$\vspace{2mm}\\
\textbf{2. Unique critical point between consecutive positive real zeros}

Assume that \(\{A_j\}_{j=1}^\infty\) is a Blaschke sequence
 in the right half-plane \(\mathbb{H}\), and the positive constants
\(A_j\) will be constructed later.
Let
\[
P(\eta)=\prod_{j=1}^\infty \frac{A_j-\eta}{A_j+\eta} \ \ \ \  (\eta\in \mathbb{H}),
\] where $$ 0<A_1<A_2<\cdots.$$
For \(\eta\ne A_j\) and $\eta\in \mathbb{H}$,
\begin{equation}
\frac{P'(\eta)}{P(\eta)}
=-\sum_{j=1}^{\infty}\frac{2A_j}{A_j^2-\eta^2}.
\label{A14}
\end{equation}
In fact, for a finite product \(P\), it is routine
to verify that (\ref{A14}) holds. In general, since \(P\) and \(P'\)
can be approximated by sequences of finite Blaschke
 products, \eqref{A14}
 holds for every \(\eta\ne A_j\) as well.

Denote the right-hand side of (\ref{A14}) by \(L(\eta)\), and
\begin{equation} \label{B14}
L(\eta)=-\sum_{j=1}^\infty \frac{2A_j}{A_j^2-\eta^2}=
\frac{P'(\eta)}{P(\eta)} .
 \end{equation} Then
\begin{equation}
L'(\eta)
=-\sum_{j=1}^\infty \frac{4A_j\eta}{(A_j^2-\eta^2)^2}<0,
\qquad \eta>0.
\label{A15}
\end{equation}
Hence \(L\) is strictly decreasing on each
 interval \((A_j,A_{j+1})\). Moreover,
\[
L(\eta)\to +\infty
\quad(\eta\to A_j^+);
\qquad
L(\eta)\to -\infty\quad(\eta\to A_{j+1}^-).
\]
By the intermediate value theorem,
there is a unique point
\(c_j\in(A_j,A_{j+1})\) such that \(L(c_j)=0\).
Since \(P(c_j)\ne0\),  (\ref{B14})
 gives \(P'(c_j)=0\). Thus \emph{there is
 exactly one real critical point of $P$ between
 any two consecutive positive real zeros $A_j$ and $A_{j+1}$}.

Furthermore, by (\ref{A15}), at such a critical point,
\[
P''(c_j)=P(c_j)L'(c_j)\ne0.
\]
Hence these critical points are simple.
\vspace{2mm}\\
\noindent \textbf{3. A symmetric four-zero block}

 A four-zero block is a set of four zeros corresponding to
  a Blaschke product of order four.

Fix \(\lambda>1\) and consider
\begin{equation}
\phi_{\lambda}(s,\eta)
=
\frac{1-\eta}{1+\eta}
\frac{\lambda-\eta}{\lambda+\eta}
\frac{\lambda^2-\eta}{\lambda^2+\eta}
\frac{\lambda^3s-\eta}{\lambda^3s+\eta},
\label{A16}
\end{equation}
where \(s>0\). When \(s=1\), the four zeros are
\[
1,\quad \lambda,\quad \lambda^2,\quad \lambda^3.
\]
For \(j=0,1,2,3\),
\[
\frac{\lambda^j-(\lambda^3/\eta)}
     {\lambda^j+(\lambda^3/\eta)}
=-\frac{\lambda^{3-j}-\eta}
        {\lambda^{3-j}+\eta},
\]
and then by (\ref{A16}) we obtain
 \begin{equation}
\phi_{\lambda}\!\left(1,\frac{\lambda^3}{\eta}\right)
=
\phi_{\lambda}(1,\eta).
\label{A17}
\end{equation}

Let \(x_{\lambda}\) denote the unique critical point
 in \((1,\lambda)\). By (\ref{A17}),
\[
y_{\lambda}=\frac{\lambda^3}{x_{\lambda}}\in (\lambda^2,\lambda^3)
\]
  is also a critical point, and satisfies
\begin{equation}
\phi_{\lambda}(1,x_{\lambda})
=
\phi_{\lambda}(1,y_{\lambda}).
\label{A18}
\end{equation}
Thus a single four-zero block already produces two distinct critical points with the same critical value.\vspace{2mm}\\
\textbf{4. The role of the perturbation parameter \(s\)}

If infinitely many four-zero blocks are multiplied together,
the remaining blocks would perturb
  both the critical points and
 the critical values of any fixed block,
 thereby destroying the exact identity \eqref{A18} in general.
  In the following, the parameter \(s>0\), chosen sufficiently close to \(1\),
 is introduced to absorb this perturbation.

To be precise,
for $s>\frac{1}{\lambda}$, let $c_1(s)$ and $c_3(s)$ be
 the unique real critical  points of  $\phi_{\lambda}(s,\cdot)$ in the intervals
\[
(1,\lambda) \quad\text{and}\quad (\lambda^2,\lambda^3s),
\]
respectively; the critical point \(c_2(s)\) in \((\lambda,\lambda^2)\) is not needed here.
  The function $\phi_{\lambda}(s,\cdot)$ is zero at each interval endpoint and negative in its
interior. Hence these critical values are precisely the minima:
\[
 m_1(s)=\min_{\eta\in [1,\lambda]}\phi_{\lambda}(s,\eta),
 \qquad
 m_3(s)=\min_{\eta\in [\lambda^2,\lambda^3s]}
 \phi_{\lambda}(s,\eta).
\]
   Define
\[
D_{\lambda}(s)=\phi_{\lambda}(s,c_1(s))-
\phi_{\lambda}(s,c_3(s)),
\]
and we have
 $$D_\lambda(s)=m_1(s)-m_3(s).$$ In particular,
$D_\lambda(1)=0$ by~\eqref{A18}.

\begin{prop}\label{global-sign}
For every $\lambda>1$, we have that
\begin{equation*}
 \begin{aligned}
 \begin{cases}
 D_\lambda(s)<0, &\text{if }\ \frac{1}{\lambda}<s<1,\\
 D_\lambda(s)>0, &\text{if }\ s>1.
 \end{cases}
 \end{aligned}
\end{equation*}
\end{prop}
\begin{proof}
Rewrite $a=\lambda^3$. For $0<t<\min\{a,as\}$ define
\[
 H_s(t)=\frac{(as-t)(a+t)}{(as+t)(a-t)}.
\]
Only the fourth factor
$\dfrac{\lambda^3s-\eta}{\lambda^3s+\eta}$ of~\eqref{A16}
 changes with $s$, and we have
\begin{equation*}
 \phi_{\lambda}(s,t)=H_s(t)\phi_{\lambda,1}(t).
\end{equation*}
We first determine the monotonicity of $H_s$. Subtracting one gives
\[
 H_s(t)-1=
 \frac{(as-t)(a+t)-(as+t)(a-t)}{(as+t)(a-t)}
 =\frac{2at(s-1)}{(as+t)(a-t)}.
\]
By direct computations,
\[
 H_s'(t)=\frac{dH_s(t)}{dt}=\frac{2a(s-1)(a^2s+t^2)}{(as+t)^2(a-t)^2}.
\]
Thus, throughout the  domain $(0,\min\{a,as\})$, $H_s$ is strictly increasing if $s>1$ and strictly decreasing if
$s<1$.

Suppose first that $s>1$. Choose $x\in(1,\lambda)$ with
$\phi_{\lambda}(s,x)=m_1(s)$ and set $y=\frac{a}{x}$. Then
\[
 \lambda^2<y<a<as  \quad \mathrm{and} \quad 0<x<y<a.
\]
In particular, $y$ lies in the third zero interval. Monotonicity gives
$H_s(x)<H_s(y)$, while the equality (\ref{A17}) gives
$\phi_{\lambda}(1,x)=\phi_{\lambda}(1,y)<0$.
Hence, by $\phi_{\lambda}(s,t)=H_s(t)\phi_{\lambda}(1,t) $
we obtain
\[
 \phi_{\lambda}(s,y)<\phi_{\lambda}(s,x).
\]
Therefore
\[
 m_3(s)=  \min_{\eta\in [\lambda^2,as]}
 \phi_{\lambda}(s,\eta)  \leq  \phi_{\lambda}(s,y)
 <\phi_{\lambda}(s,x)=m_1(s).
\]
Hence $D_\lambda(s)>0$.

Now let $\frac{1}{\lambda}<s<1$. Choosing  $y\in(\lambda^2,as)$ with
$\phi_{\lambda}(s,y)=m_3(s)$ and put $x=\frac{a}{y}$, we have
\[
 1<\frac{1}{s}<x<\lambda,\qquad 0<x<y<as<a.
\]
Thus both $x$ and $y$ lie in the domain $(0,\min\{a,as\})$
  where $H_s$ is strictly decreasing. Consequently
$H_s(x)>H_s(y)$. Multiplying again by the common negative value
$\phi_{\lambda}(1,x)=\phi_{\lambda}(1,y)$ gives
\[
 \phi_{\lambda}(s,x)<\phi_{\lambda}(s,y).
\]
It follows that
\[
 m_1(s)=\min_{\eta\in [1,\lambda]}\phi_{\lambda}(s,\eta)
 \leq  \phi_{\lambda}(s,x)
 <\phi_{\lambda}(s,y)=m_3(s),
\]
which immediately gives $D_\lambda(s)<0$.
\end{proof}

For the construction, choose an increasing sequence  $\lambda_k \geq 2$ with $\lambda_k\to\infty$. Fix, once and for all,
\[
 \boxed{\varepsilon_0=\frac{1}{4},\qquad I=\Big[\frac{3}{4}, \frac{5}{4}\Big].}
\]
Since $\frac{3}{4}>\frac{1}{\lambda_k}$ for every $k$, Proposition~\ref{global-sign}
gives
\begin{equation}
 D_{\lambda_k}(1-\varepsilon_0)<0,
 \qquad D_{\lambda_k}(1+\varepsilon_0)>0
 \qquad(k\geq 1). \label{eq:20}
\end{equation}
This is one common $\varepsilon_0$.
  Let \(\{\lambda_k\}_{k=1}^\infty\) be any  increasing sequence
   tending to infinity with \(\lambda_1>2\). Then
    \(0<\varepsilon_0<1-\frac{1}{\lambda_k}\) holds for all \(k\).
 \subsection{Patching together infinitely many four-zero blocks}
 In what follows, we will construct a sequence
of four-zero blocks to build an infinite Blaschke
 product satisfying our requirements.

Let the four zeros in the $k$th block $J_k$ be
\begin{equation}
 R_k,\quad R_k\lambda_k,\quad R_k\lambda_k^2,
 \quad R_k\lambda_k^3s_k,
  \label{eq:A21}
\end{equation} with $s_k\in I.$
Set
\[
 \Lambda_k=\lambda_k^3(1+\varepsilon_0)=\tfrac54\lambda_k^3,
 \quad \mathrm{and} \quad J_k=[R_k,\Lambda_kR_k].
\]  The interval $J_k$ contains every possible
zero of block $k$ and is
independent of the actual value of $s_k$.
Later, the scales \(R_k\) will be chosen to increase
sufficiently rapidly so that distinct blocks
 are widely separated.
Remind that $$
\phi_{\lambda}(s,\eta)
=
\frac{1-\eta}{1+\eta}
\frac{\lambda-\eta}{\lambda+\eta}
\frac{\lambda^2-\eta}{\lambda^2+\eta}
\frac{\lambda^3s-\eta}{\lambda^3s+\eta},
$$ and for each
$ \mathbf s=(s_1,s_2,\cdots)\in I^{\infty}$,  define
\begin{equation}
P_{\mathbf s}(\eta)
=
\prod_{k=1}^{\infty}
\phi_{\lambda_k}\Big(s_k,\frac{\eta}{R_k}\Big),
\ \ \ \ \eta\in \mathbb{H}.
 \label{eq:22}
\end{equation}
We first fix the error term appearing in \eqref{eq:20}.
 Note that the numbers
\begin{equation}
 \mu_k=\frac12\min\{-D_{\lambda_k}(1-\varepsilon_0),
                         D_{\lambda_k}(1+\varepsilon_0)\}>0
\label{eqmu}
\end{equation}
are well defined. Put
\[
 \beta_k=\min\Big\{\frac{\mu_k}{8},\,2^{-k-4}\Big\}>0.
\]
The  errors $\mu_k$ and $\beta_k$ need not be bounded away from zero.

The two intervals associated with block \(k\) are
\[
J_{k,1}=[R_k,R_k\lambda_k],\qquad
J_{k,3} =[R_k\lambda_k^2,R_k\lambda_k^3s_k],
\]
both contained in \(J_k\).  The functions
$\phi_{\lambda_k}(s_k,\frac{\eta}{R_k})$ and
$P_{\mathbf s} $ vanish at the endpoints and
 are negative inside,
 since the other blocks have positive product there.
 By Subsection~\ref{orpre},
 $P_{\mathbf s} $ has exactly one critical point in either open interval,
 and its value is the interval minimum.
  Let \(c_1^{(k)}(\mathbf s)\) and \(c_3^{(k)}(\mathbf s)\)
    denote the critical points of
     $P_{\mathbf s}$ in $J_{k,1}$ and
$J_{k,3}$, respectively.
  Let \(\Phi_k(\mathbf s)\) denote the difference between the corresponding two critical values of \(P_{\mathbf s}\):
\[
\Phi_k(\mathbf s)=P_{\mathbf s}(c_1^{(k)}(\mathbf s))-P_{\mathbf s}(c_3^{(k)}(\mathbf s))
=\min_{\eta\in J_{k,1}}P_{\mathbf s}(\eta)
-\min_{\eta\in J_{k,3} }P_{\mathbf s}(\eta).
\]
We have the following key result.
\begin{thm}\label{key1} With a fixed increasing
 sequence \(\{\lambda_k\}_{k=1}^\infty\)
 satisfying \(\lambda_1>2\), there exists
  a carefully chosen infinite sequence \(\{R_k\}_{k=1}^\infty\), associated with \(\{\lambda_k\}_{k=1}^\infty\)
as in \eqref{eq:A21}, such that
\begin{equation}
 \begin{aligned}
 \Phi_k(\mathbf s)&\leq   -\mu_k<0&&\text{when }s_k=1-\varepsilon_0,\\
 \Phi_k(\mathbf s)&\geq\mu_k>0&&\text{when }s_k=1+\varepsilon_0,
 \end{aligned}
 \qquad\varepsilon_0= \frac{1}{4},
  \label{A23}
\end{equation}
where $\mathbf s=(s_1,s_2,\cdots)\in I^\infty. $
\end{thm}
\noindent
 Note that \(\mu_k\), defined in (\ref{eqmu}),
  depends only on \(\lambda_k\), not on \(R_k\).
 Later, we will see that,
 by the construction,
the infinite sequence  \eqref{eq:A21} yields a thin Blaschke
 sequences.

To prove Theorem~\ref{key1}, we begin with two facts.

\noindent \textbf{Fact 1.} If
$0\leq  t_j\leq  1$ holds for each $j,$  then
 \begin{equation}
 0\leq  1-\prod_{j=1}^mt_j
\leq  \sum_{j=1}^m(1-t_j).
 \label{fact1}
\end{equation} This also holds for $m=\infty.$

\noindent  \textbf{Fact 2.} Let both $f$ and $g$ be
continuous real-valued functions on a compact interval $J$.
 Then
\begin{equation}
 \Big|\min_J f-\min_J g\Big|\leq  \sup_J|f-g|.
  \label{eq:min-est}
\end{equation}

For completeness, note that (\ref{fact1}) follows from the identity
\[
 1-t_1\cdots t_m=(1-t_1)+t_1(1-t_2)+\cdots+
 t_1\cdots t_{m-1}(1-t_m),
\]
and every coefficient of $1-t_j$ is at most one.
To verify Fact 2, let \(z_0\) be a point in \(J\) satisfying \(g(z_0)=\min_J g\). Then
\[
\min\limits_J f-\min\limits_J g \leq  f(z_0)-g(z_0) \leq  \sup\limits_J|f-g|.
\]
Similarly, \(\min\limits_J g-\min\limits_J f \leq  \sup\limits_J|f-g|\).
We therefore obtain \eqref{eq:min-est}.

 Before proving Theorem \ref{key1}, we first introduce
 the construction of the sequence \(\{R_k\}_{k=1}^\infty\)
 mentioned above. Given that
\(\{\lambda_j\}_{j=1}^\infty\) is a fixed increasing sequence
 with \(\lambda_1>2\), we will choose \(R_k\) so that
\[
R_j\ll R_j\lambda_j\ll R_{j+1}.
\]
Precisely,
take $R_1>2$. Having chosen $R_1,\cdots,R_{n-1}$, choose $R_n$ so that
\begin{align}
 &R_n>2^n,\qquad R_n>2\Lambda_{n-1}R_{n-1},
  \label{22a}\\
 &\frac8{R_n}\sum_{j<n}\Lambda_jR_j<\beta_n,
 \label{22b}\\
 &\frac{8\Lambda_jR_j}{R_n}<2^{-n}\beta_j
 \qquad(1\leq  j<n),
 \label{22c}
\end{align}
where $\beta_n=\min\big\{\frac{\mu_n}{8},\,2^{-n-4}\big\}>0.$
At stage $n$ there are only finitely many requirements, and every
quantity other than $R_n$ is already fixed.
This choice is maintained throughout the
 remainder of Section \ref{6cexample}.

\begin{proof}[\textbf{\emph{Proof of Theorem \ref{key1}}.}]
Since the proof of this theorem is relatively lengthy, we divide it into two steps.

 \textbf{Step 1.} \emph{We first prove}
\begin{equation}
 \sup_{\mathbf s\in I^{\infty}}\sup_{\eta\in J_k}
 \Big|P_{\mathbf s}(\eta)-\phi_{\lambda_k}\Big(s_k,\frac{\eta}{R_k}\Big)\Big|<2\beta_k.
 \label{22h}
\end{equation}

Fix $k$ and $\eta\in J_k$, We will estimate \(1-\phi_{\lambda_k}(s_k,\frac{\eta}{R_k})\)
according to whether \(j<k\) or \(j>k\). If $j<k$, all four zeros $A$ of block $j$
are less than $\eta$. Its four minus signs cancel, and thus
\[
 \phi_{\lambda_j}(s_j,\frac{\eta}{R_j}) =\prod_{A\text{ in block }j}\frac{\eta-A}{\eta+A}.
\]
Each factor belongs to $(0,1)$ and
\[
 1-\frac{\eta-A}{\eta+A}=\frac{2A}{\eta+A}
 \leq  \frac{2\Lambda_jR_j}{R_k}.
\]
Hence, by~\eqref{fact1},
\begin{equation}
 0\leq  1- \phi_{\lambda_j}\Big(s_j,\frac{\eta}{R_j}\Big)\leq  \frac{8\Lambda_jR_j}{R_k}
 \qquad(j<k). \label{eq:earlier}
\end{equation}
If $j>k$, all its zeros satisfy $A>\eta$ and
\[
 1-\frac{A-\eta}{A+\eta}=\frac{2\eta}{A+\eta}
 \leq  \frac{2\Lambda_kR_k}{R_j}.
\]
Thus
\begin{equation}
 0\leq  1- \phi_{\lambda_j}\Big(s_j,\frac{\eta}{R_j}\Big) \leq  \frac{8\Lambda_kR_k}{R_j}
 \qquad(j>k). \label{eq:later}
\end{equation}
All right-hand sides are independent of
the values of the parameters $s_j$.

Write $$Q_{k,\mathbf s}(\eta)=\prod_{{\substack{j\geq 1\\ j \neq k}}}
\phi_{\lambda_j}\Big(s_j,\frac{\eta}{R_j}\Big),$$  and
$$P_{\mathbf s}(\eta)=\phi_{\lambda_k}
\Big(s_k,\frac{\eta}{R_k}\Big)Q_{k,\mathbf s} .$$  On $J_k$,
 every factor in
$Q_{k,\mathbf s}$ is positive and at most one. Applying
\eqref{fact1} to finite products and passing to the
limit, then using~\eqref{eq:earlier}--\eqref{eq:later} and
\eqref{22b}--\eqref{22c}, gives
\begin{align*}
 0\leq  1-Q_{k,\mathbf s}(\eta)
 &\leq  \frac8{R_k}\sum_{j<k}\Lambda_jR_j
       +\sum_{j>k}\frac{8\Lambda_kR_k}{R_j}\\
 &<\beta_k+\beta_k\sum_{j=k+1}^{\infty}2^{-j}\\
 &=\beta_k(1+2^{-k})\leq    \frac{3}{2}\beta_k.
\end{align*}
Since $|\phi_{\lambda_k}(s_k,\frac{\eta}{R_k})| \leq  1$ on the positive
real axis, and
  \begin{equation*}
\Big|P_{\mathbf s}(\eta)-\phi_{\lambda_k}\Big(s_k,\frac{\eta}{R_k}\Big)\Big|
=\Big| \phi_{\lambda_k}\Big(s_k,\frac{\eta}{R_k}\Big)\Big|~
 |1-Q_{k,\mathbf s}(\eta)|
  \end{equation*}
  we immediately obtain (\ref{22h}).

 \textbf{Step 2.}
 \emph{We now give some estimates to prove} \eqref{A23}.

The idea is as follows:
\[
D_{\lambda_k}(1-\varepsilon_0)<0,
\qquad
D_{\lambda_k}(1+\varepsilon_0)>0.
\]
These inequalities hold for every \(k\).
Using a perturbation of \(D_{\lambda_k}(s_k)\) given by \(\Phi_k(\mathbf s)\),
 we obtain the desired assertion \eqref{A23}.

By noting
\[
 D_{\lambda_k}(s_k)=
 \min_{\eta\in J_{k,1}}
  \phi_{\lambda_k}\Big(s_k,\frac{\eta}{R_k}\Big)-
  \min_{\eta\in J_{k,3} }\phi_{\lambda_k}\Big(s_k,\frac{\eta}{R_k}\Big)
\]
and \[J_{k,1}=[R_k,R_k\lambda_k],\qquad
J_{k,3} =[R_k\lambda_k^2,R_k\lambda_k^3s_k],\]
we immediately have
\[ D_{\lambda_k}(s_k)=
 \min_{\eta\in [1, \lambda_k]}
  \phi_{\lambda_k}(s_k,\eta)-
  \min_{\eta\in [\lambda_k^2,\lambda_k^3s_k]}
  \phi_{\lambda_k}(s_k,\eta).\]
Hence \(D_{\lambda_k}(s_k)\) does not depend on \(R_k\).
 This observation allows us to choose \(R_k\)
freely with respect to the value \(D_{\lambda_k}(s_k)\).

Applying (\ref{eq:min-est})  on the two intervals and using (\ref{22h}),
we obtain
\begin{equation}
 \begin{aligned}
 |\Phi_k(\mathbf s)-D_{\lambda_k}(s_k)|
 &\leq  \Big|\min_{\eta\in J_{k,1}}P_{\mathbf s}(\eta)
 -\min_{\eta\in J_{k,1}}  \phi_{\lambda_k}(s_k,\frac{\eta}{R_k})\Big|\\
&\ \ \ +\Big|\min_{\eta\in J_{k,3} }P_{\mathbf s}
                       -\min_{\eta\in J_{k,3} }
                       \phi_{\lambda_k}(s_k,\frac{\eta}{R_k})\Big|\\
 &<4\beta_k\leq  \frac{\mu_k}{2}<\mu_k.\label{eq:value-approx}
 \end{aligned}
\end{equation}
 The definition of $\mu_k$ gives
\[
 D_{\lambda_k}(1-\varepsilon_0)\leq  -2\mu_k,
 \qquad D_{\lambda_k}(1+\varepsilon_0)\geq 2\mu_k.
\]
If $s_k=1-\varepsilon_0$, then by (\ref{eq:value-approx})
\[
 \Phi_k(\mathbf s)<D_{\lambda_k}(1-\varepsilon_0)+\mu_k
\leq  -\mu_k<0.
\]
If $s_k=1+\varepsilon_0$, then
\[
 \Phi_k(\mathbf s)>D_{\lambda_k}(1+\varepsilon_0)-\mu_k
\geq \mu_k>0.
\]
We thereby obtain (\ref{A23}), as desired.
\end{proof}

\subsection{The construction of $P_{\mathbf s}$}

\textbf{1. Some observations.}

We take the product topology on
 \(I^\infty =I\times I\times\cdots\), where
$
I=\left[\frac34,\frac54\right].
$
Its neighborhood basis consists of sets
\[
U_1\times U_2\times\cdots\times U_k\times I\times I\times\cdots,
\]
where \(k<\infty\) and each \(U_j\) \((1\leq  j\leq  k)\) is a relatively open subset of \(I\).
For \(\mathbf s^N\) and \(\mathbf s^*\) in \(I^\infty\),
 if \(\{\mathbf s^N\}_{N=1}^\infty\)
converges to \(\mathbf s^*\) in each coordinate $I^\infty$,
then  \(\{\mathbf s^N\}_{N=1}^\infty\)
converges to \(\mathbf s^*\) in the product topology of $I^\infty$.

For each $m<\infty$,
\[
P_{\mathbf s,m}(\eta)=
\prod_{k=1}^m \phi_{\lambda_k}
\left(s_k,\frac{\eta}{R_k}\right)
\]
as the partial product of $P_{\mathbf s}. $
In the following, we abbreviate ``Observation" as ``Ob."
\begin{itemize}
\item[ (\textbf{Ob.}1)] \emph{Let $1\leq m <\infty.$ For \(\mathbf s^N\) and \(\mathbf s^*\) in \(I^\infty\), if \(\{\mathbf s^N\}_{N=1}^\infty\) converges to \(\mathbf s^*\) in each coordinate, then
\[
P_{\mathbf s,m}(\mathbf s^N,\eta)\to P_{\mathbf s,m}(\mathbf s^*,\eta)
\]  uniformly for $\eta $ in each compact subset
 \(K\subseteq \mathbb H\).}

\ \ \ \ This follows directly from the definition
of \(P_{\mathbf s,m}\), which involves
only finitely many factors of the form
\(\phi_{\lambda_k}(s_k,\frac{\eta}{R_k})\),
 and only finitely many of
\(s_1,\cdots,s_m\) are involved.

\item[(\textbf{Ob.}2)] \emph{For each compact subset
 \(K\subseteq \mathbb H\), we have
 \[
\sup_{(\mathbf s_N,\mathbf s_{N+1},\ldots)\in I^\infty}~\sup_{\eta\in K}
\left|
\prod_{k=N+1}^{\infty}\phi_{\lambda_k}\left(s_k,\frac{\eta}{R_k}\right)-1
\right|
\to 0
\quad(N\to\infty)
\] and
\[
\sup_{\mathbf s\in I^{\infty}}\sup_{\eta\in K}
 |P_{\mathbf s}(\eta)-\prod_{j=1}^N\phi_{\lambda_j}
(s_j,\frac{\eta}{R_j})|\to 0,
\quad N\to\infty.
\]}

\ \ \ \ In fact, for a compact subset \(K\subseteq \mathbb H\),
put \(M_K=\max\{|\eta|:\eta\in K\}\).
For those \(R_j>M_K\), each zero \(A\)
in block \(j\) satisfies \(A \geq R_j\), and
\[
\left|\frac{A-\eta}{A+\eta}-1\right|\leq  \frac{2M_K}{R_j},
\qquad
\left|\frac{A-\eta}{A+\eta}\right|\leq  1.
\]
Applying \eqref{fact1} gives
\[
\left|\phi_{\lambda_j}\left(s_j,\frac{\eta}{R_j}\right)-1\right|\leq  \frac{8M_K}{R_j}.
\]
Furthermore,
\begin{align*}
\sup_{\mathbf s\in I^{\infty}}\sup_{\eta\in K}
 \bigg | P_{\mathbf s}(\eta)&-\prod_{j=1}^N\phi_{\lambda_j}
\Big(s_j,\frac{\eta}{R_j}\Big)\bigg|\\
&= \sup_{\mathbf s\in I^{\infty}}\sup_{\eta\in K}
\bigg|\prod_{j=1}^N\phi_{\lambda_j}
\Big(s_j,\frac{\eta}{R_j}\Big) \bigg|~\bigg|  \prod_{j=N+1}^\infty \phi_{\lambda_j}
\Big(s_j,\frac{\eta}{R_j}\Big)-1\bigg|\\
&\leq   \sup_{\mathbf s\in I^{\infty}}\sup_{\eta\in K}
 \bigg|\prod_{j=N+1}^\infty \phi_{\lambda_j}
\Big(s_j,\frac{\eta}{R_j}\Big)-1\bigg| \\
&\leq  8M_K\sum_{j>N}\frac1{R_j}\to 0, \ \ N\to \infty.
 \label{22d}
\end{align*}
The formula $\sum\limits_{j>N}\frac1{R_j}\to 0\
(N\to \infty)$ follows from  \eqref{22a}.

\item[(\textbf{Ob.}3)] \textbf{Claim:} \emph{For \(\mathbf s^N\) and \(\mathbf s^*\) in \(I^\infty\), if \(\{\mathbf s^N\}_{N=1}^\infty\) converges to \(\mathbf s^*\)
in each coordinate,
then
\[
c_1^k(\mathbf s^N)\to c_1^k(\mathbf s^*)
\quad\text{and}\quad
c_3^k(\mathbf s^N)\to c_3^k(\mathbf s^*),
\quad  N\to\infty.
\]
Recall that  \(c_1^k(\mathbf s)\) and \(c_3^k(\mathbf s)\) are
 the unique critical points of \(P_{\mathbf s}\)
 in \(J_{k,1}\) and \(J_{k,3}\), respectively.}

 \ \ \ \  In Subsection~\ref{orpre}, by the arguments in the paragraph above titled ``A symmetric four-zero block'', we have that \(P_{\mathbf s^N}\) and \(P_{\mathbf s^*}\)
have exactly one critical point in \(J_{k,1}\).
For any $\varepsilon'>0,$ we
now show that there exists an integer \(N_0  \geq 1\)
 such that
\[
\left|c_1^k(\mathbf s^N)-c_1^k(\mathbf s^*)\right|<
\varepsilon'
\qquad\text{for all }N \geq N_0.
\]
For this, we will apply Rouch$\acute{\mathrm{e}}$'s theorem on
a closed elliptic disk
  \(K\) $(K\subseteq \mathbb{H})$ whose interior contains \(J_{k,1}\) and \(J_{k,3}\).
There exists a $\delta>0$ such that
 if a holomorphic function \(h\) on \(K\) satisfies
\[
\sup_{\eta\in K}|h(\eta)-P_{\mathbf s^*}(\eta)|<\delta,
\]
then \(h\) has exactly one critical
point lying in an \(\varepsilon'\)-neighborhood of
 \(c_1^k(\mathbf s^*)\).
 We   take \(h=P_{\mathbf s^N}\),
 and the unique critical point is $c_1^k(\mathbf s^N) $.
We will shortly prove that
\[
\sup_{\eta\in K}
\left|P_{\mathbf s^N}(\eta)-P_{\mathbf s^*}(\eta)\right|
<\delta
\]
for sufficiently large $N.$

To be precise, by Ob.(2), for all $\eta$ on this compact set \(K\),
\[
\left|P_{\mathbf s^N,m}(\eta)-P_{\mathbf s^N}(\eta)\right|
<\frac{\delta}{4}
 ~~\mathrm{and}~~
 \left|P_{\mathbf s^*,m}(\eta)-P_{\mathbf s^*}(\eta)\right|<\frac{\delta}{4}
\]
for some $m.$
Since \(m\) is fixed, by (Ob.1), there exists an integer $N_0\geq 1$ such that
\[
\left|P_{\mathbf s^*,m}(\eta)-P_{\mathbf s^N,m}(\eta)\right|<\frac{\delta}{4},\quad \eta \in K
\]
for   \(N\geq N_0\). Hence we obtain that for those $N,$
\[
\begin{aligned}
\left|P_{\mathbf s^N}(\eta)-P_{\mathbf s^*}(\eta)\right|
&\leq
\left|P_{\mathbf s^N}(\eta)-P_{\mathbf s^N,m}(\eta)\right|+
\left|P_{\mathbf s^N,m}(\eta)-P_{\mathbf s^*,m}(\eta)\right|\\
&\ \ \  +
\left|P_{\mathbf s^*,m}(\eta)-P_{\mathbf s^*}(\eta)\right|\\
&<\frac{3\delta}{4}<\delta,
\qquad \eta\in K,
\end{aligned}
\]
as desired.
In this way, we can find an \(N_0\geq  1\) such that
\[
\left|c_1^k(\mathbf s^N)-c_1^k(\mathbf s^*)\right|<\varepsilon'
\ \text{for\  all}\ N \geq N_0.\]
A similar argument applies to \(\left|c_3^k(\mathbf s^N)-c_3^k(\mathbf s^*)\right|\).

\item[(\textbf{Ob.}4)]\emph{Remind that \[
\Phi_k(\mathbf s)=P_{\mathbf s}\bigl(c_1^k(\mathbf s)\bigr)-P_{\mathbf s}\bigl(c_3^k(\mathbf s)\bigr).
\]
 For \(\mathbf s^N\) and \(\mathbf s^*\) in \(I^\infty\), if \(\{\mathbf s^N\}_{N=1}^\infty\) converges to \(\mathbf s^*\) in each coordinate, then
 $$ \Phi_k(\mathbf s^N)\to \Phi_k(\mathbf s^*), \ \
 N\to \infty.$$}

    To see this, note that
\[
\begin{aligned}
\left|\Phi_k(\mathbf s^N)-\Phi_k(\mathbf s^*)\right|
&=\left|
P_{\mathbf s^N}\bigl(c_1^k(\mathbf s^N)\bigr)
-P_{\mathbf s^N}\bigl(c_3^k(\mathbf s^N)\bigr)
-P_{\mathbf s^*}\bigl(c_1^k(\mathbf s^*)\bigr)
+P_{\mathbf s^*}\bigl(c_3^k(\mathbf s^*)\bigr)
\right|\\
&\leq
\left|P_{\mathbf s^N}\bigl(c_1^k(\mathbf s^N)\bigr)
-P_{\mathbf s^*}\bigl(c_1^k(\mathbf s^*)\bigr)\right|
 +
\left|P_{\mathbf s^N}\bigl(c_3^k(\mathbf s^N)\bigr)-P_{\mathbf s^*}\bigl(c_3^k(\mathbf s^*)\bigr)\right|.
\end{aligned}
\]
For every $\varepsilon>0$,  it suffices to  prove
\[
\left|P_{\mathbf s^N}\bigl(c_1^k(\mathbf s^N)\bigr)-P_{\mathbf s^*}\bigl(c_1^k(\mathbf s^*)\bigr)\right|<3\varepsilon
\]
for large \(N\), and similarly,
\[
\left|P_{\mathbf s^N}\bigl(c_3^k(\mathbf s^N)\bigr)-P_{\mathbf s^*}\bigl(c_3^k(\mathbf s^*)\bigr)\right|<3\varepsilon.
\]
Recall that
\[
P_{{\mathbf s},m}(\eta)=\prod_{k=1}^m \phi_{\lambda_k}\left(s_k,\frac{\eta}{R_k}\right).
\]
Then
\[
\begin{aligned}
&\left|P_{\mathbf s^N}\bigl(c_1^k(\mathbf s^N)\bigr)-P_{\mathbf s^*}\bigl(c_1^k(\mathbf s^*)\bigr)\right|\\
&\leq
\left|P_{\mathbf s^N}\bigl(c_1^k(\mathbf s^N)\bigr)-P_{\mathbf s^N,m}\bigl(c_1^k(\mathbf s^N)\bigr)\right|\\
& \ \ \ +
\left|P_{\mathbf s^N,m}\bigl(c_1^k(\mathbf s^N)\bigr)-P_{\mathbf s^*,m}\bigl(c_1^k(\mathbf s^*)\bigr)\right|\\
& \ \ \ +
\left|P_{\mathbf s^*}\bigl(c_1^k(\mathbf s^*)\bigr)-P_{\mathbf s^*,m}\bigl(c_1^k(\mathbf s^*)\bigr)\right|.
\end{aligned}
\]
By (Ob.2), the first term
$
\left|P_{\mathbf s^N}\bigl(c_1^k(\mathbf s^N)\bigr)-P_{\mathbf s^N,m}\bigl(c_1^k(\mathbf s^N)\bigr)\right|
$
and the third term
$
\left|P_{\mathbf s^*}\bigl(c_1^k(\mathbf s^*)\bigr)-P_{\mathbf s^*,m}\bigl(c_1^k(\mathbf s^*)\bigr)\right|
$
can be made less than \(\varepsilon\) for some fixed \(m\) and all \(N \geq 1\). Since \(P_{\mathbf s,m}\) depends only on finitely many products of
$\phi_{\lambda_k}(s_k,\frac{\eta}{R_k})$,
the middle term converges uniformly for
 \(\mathbf s\in I^\infty\), essentially
 for \((\mathbf s_1,\ldots,\mathbf s_m)\in I^m\).
 Then we have  by (Ob.3) that
\[
c_1^k(\mathbf s^N)\to c_1^k(\mathbf s^*),\ \   N\to\infty,
\]
and thus the middle term is less than
 \(\varepsilon\) for sufficiently large \(N\), and
 then
 $$
 \left|P_{\mathbf s^N}\bigl(c_1^k(\mathbf s^N)\bigr)
 -P_{\mathbf s^*}\bigl(c_1^k(\mathbf s^*)\bigr)\right|<
  3\varepsilon,$$ as desired.
 \end{itemize}
In the following construction, we will specialize
\(\mathbf s^N\) and \(\mathbf s^*\) in \(I^{\infty}\).\vspace{2mm}\\
\textbf{2. Completion of the construction.}

We will seek a parameter vector \(\mathbf s^*\) in \(I^{\infty}\) such that
\begin{equation}
\Phi_k(\mathbf s^*)=0,
\qquad k=1,2,\cdots.
\label{A24}
\end{equation}

 The result used here is the
 \textbf{Poincar\'e--Miranda Theorem}
  (also called Miranda's Theorem), the multivariable
  analogue of the intermediate value theorem and
  a theorem equivalent to the Brouwer fixed point theorem
  (see \cite{Miranda1940}).
   \begin{thm}[Poincar\'e--Miranda Theorem]
\label{mira}
Let
\[
Q=\prod_{j=1}^{N}[a_j,b_j]\subset\mathbb R^N,
\qquad
F=(F_1,\cdots,F_N):Q\to\mathbb R^N
\]
be continuous. Suppose that, for each \(j=1,\cdots,N\),
\[
F_j(x)\leq  0\quad\text{on the face }x_j=a_j;
\qquad
F_j(x)\geq 0\quad\text{on the face }x_j=b_j.
\]
Then there exists \(x^*\in Q\) such that
\[
F(x^*)=0.
\]
\end{thm}
 We apply Theorem \ref{mira}  to the present construction.
Fix \(N\). Freeze the tail parameters, for instance by setting
\[
s_{N+1}=s_{N+2}=\cdots=1,
\]
and consider the finite-dimensional rectangle
\[
Q_N=\prod_{j=1}^{N}
[1-\varepsilon_0,1+\varepsilon_0],
 ~~\mathrm{with}~~\varepsilon_0=\frac{1}{4}.
\]
Define
\[
F^{(N)}(s_1,\cdots,s_N)
=
\bigl(\Phi_1(\mathbf s),\cdots,\Phi_N(\mathbf s)\bigr).
\]
By (Ob.4), each function $\Phi_1(\mathbf s),
\cdots,\Phi_N(\mathbf s)$
is continuous in $\mathbf s$, and hence $F^{(N)}$ is
continuous. By Theorem \ref{key1} 
 for every \(j\), the \(j\)-th component \(\Phi_j\)
  is negative on
the face \(s_j=1-\varepsilon_0\) and positive on the opposite face
\(s_j=1+\varepsilon_0\), uniformly with respect to the remaining
coordinates. Hence the Poincar\'e--Miranda Theorem applies and yields a point
of \(Q_N\) at which
\[
\Phi_1=\cdots=\Phi_N=0.
\]
Appending the frozen tail coordinates produces
 an infinite parameter vector \(\mathbf s^N\).
 Thus, for $N=1,2,\cdots,$ we have
 \[
\begin{aligned}
&  \Phi_1(\mathbf s^1)=0,\\
& \Phi_1(\mathbf s^2)= \Phi_2(\mathbf s^2)=0,\\
 & \Phi_1(\mathbf s^3)= \Phi_2(\mathbf s^3)=\Phi_3(\mathbf s^3)=0,\\
 & \cdots \cdots\\
& \Phi_1(\mathbf s^k)= \cdots=\Phi_k(\mathbf s^k)=0,\\
& \cdots\cdots
 \end{aligned}
\]
For each \(k\), every coordinate \(s_k^j\) lies in the compact interval \(I\). By a diagonal subsequence argument, there exists a subsequence of \(\{\mathbf s^N\}_{N=1}^\infty\) that converges in each coordinate. Without loss of generality, suppose that \(\{\mathbf s^N\}_{N=1}^\infty\) converges in each coordinate; then there exists \(s_k^*\in I\) for each \(k\) such that
\[
s_k^N \to  s_k^*,\ \     N\to \infty.
\]
Set
\[
\mathbf s^*=(s_1^*,s_2^*,\cdots)\in I^\infty.
\]
For each fixed $k$, since $\Phi_k(\mathbf s^j)=0$
holds for $ j=k,k+1,\cdots$, by (Ob.4) we immediately
conclude that $\Phi_k(\mathbf s^*)=0$.
That is, we find an $\mathbf s^* \in I^\infty$
such that
\[
\Phi_k(\mathbf s^*)=0
\qquad\text{for every }k.
\]
In other words, in each block
the infinite product $P_{\mathbf s^*}$
 has two distinct critical points \(\xi_k,\eta_k\) such that
\begin{equation}
P_{\mathbf s^*}'(\xi_k)=P_{\mathbf s^*}'(\eta_k)=0,
\qquad
P_{\mathbf s^*}(\xi_k)=P_{\mathbf s^*}(\eta_k).
\label{A25}
\end{equation}
By definition, $P_{\mathbf s^*}$ is not essentially injective on its
critical points.
\begin{rem} The Poincar\'e--Miranda Theorem itself
 presents no difficulty.
 The two points that require explicit
 justification are: first, the sign condition (\ref{A23})
  must hold uniformly with respect to all
  remaining parameters; second, when one passes
  from finite-dimensional solutions to the diagonal
  limit, the continuity of each \(\Phi_k\)
  with respect to the infinite parameter vector
  must be specified. \end{rem} 

\subsection{Verification:
 thinness, primeness, essential injectivity}
\textbf{1. Thinness  of the
 resulting infinite product}

Let \(\{A_j\}_{j=1}^\infty\) be the union of the
 four-zero blocks in \eqref{eq:A21}:
\[
R_k,\qquad R_k\lambda_k,\qquad R_k\lambda_k^2,\qquad R_k\lambda_k^3s_k.
\]
To pass from a positive real half-plane zero
 \(A\) to the corresponding disk zero, we use the Cayley transform
\[
a=h(A)=\frac{A-1}{A+1}.
\]
In what follows, we write \(A\) and
\(C\) for two distinct elements of \(\{A_j\}_{j=1}^\infty\).

For \(A,C>0\),
a direct calculation gives
\begin{equation}
d(h(A),h(C))=\frac{|A-C|}{A+C},
\label{A26}
\end{equation}
where $d(z,w)=\Big|\dfrac{z-w}{1-\overline{z}w}\Big|$.
Set
\[
r=\frac{\min(A,C)}{\max(A,C)}.
\]
Then
\begin{equation*}
d(h(A),h(C))=\frac{1-r}{1+r}.
\end{equation*}
Thus \(\rho(h(A),h(C))\) is close to \(1\) whenever \(r\) is small. For \(0\leq  r\leq  \frac{1}{2}\),
\begin{equation}
0\leq  -\log\frac{1-r}{1+r}=\int_0^r \frac{2}{1-t^2}dt\leq  \frac{8}{3}r
\leq   4r.
\label{A28}
\end{equation}
Choose \(\lambda_k\to\infty\) and let \(R_k\) grow sufficiently rapidly so that, for every zero \(A\) in the \(k\)-th block,
\begin{equation}
\sum_{\{C:~C\ne A\}}\frac{\min(A,C)}{\max(A,C)}\to 0 \ \ \ \  (k\to\infty).
\label{A29}
\end{equation}
Once (\ref{A29}) is verified,
then by (\ref{A26})--(\ref{A29}),
\begin{align*}
-\log\prod_{\{C:~C\ne A\}} d(h(A),h(C))
&=\sum_{\{C:~C\ne A\}}-\log d(h(A),h(C))\\
&\leq  4\sum_{\{C:~C\ne A\}}\frac{\min(A,C)}{\max(A,C)}
\to 0.
\end{align*}
Therefore
\[
\prod_{\{C:~C\ne A\}}\rho(h(A),h(C)) \to 1,
\]
and the resulting Blaschke product is thin.

It remains to verify (\ref{A29}) uniformly in all parameters.  Rewrite
 $$S_k(A)=\sum_{\{C:~C\ne A\}}\frac{\min(A,C)}{\max(A,C)},$$
where \(k\) indicates that \(A\) lies
 in the \(k\)-th zero block.
   For two consecutive zeros in block $k$, the smaller-to-larger ratios
are  $\frac{1}{\lambda_k}$, $\frac{1}{\lambda_k}$, and $\frac{1}{\lambda_ks_k}$.
They are all at most $\frac{4}{3\lambda_k}<1$.
  Thus the three other zeros in the
   same block give terms summing to at most
$\frac{4}{\lambda_k}$ to $S_k(A)$.

For earlier blocks, $\frac{C}{A}\leq  \frac{\Lambda_jR_j}{R_k}$
 and there are four
zeros per block. The sum of their terms is at most
\[
 \frac4{R_k}\sum_{j<k}\Lambda_jR_j<\frac{\beta_k}{2}.
\]
For later blocks, $\frac{A}{C}\leq  \frac{\Lambda_kR_k}{R_j}$,
  and their terms sum to at most
 \[
 \sum_{j>k}\frac{4\Lambda_kR_k}{R_j}
 <\frac{\beta_k}{2}\sum_{j>k}2^{-j}
 =\frac{\beta_k}{2}2^{-k}.
\]
Consequently
\[
 S_k(A)\leq   \frac4{\lambda_k}
        +\frac{\beta_k}{2}(1+2^{-k})
 \rightarrow0,
\]
which proves~\eqref{A29}.\vspace{2mm}\\
\textbf{2. Primeness}

All half-plane zeros satisfy \(A\to\infty\), and hence
\[
h(A)=\frac{A-1}{A+1}\to  1.
\]
Thus the zeros have the unique boundary
accumulation point \(1\).
The following result then shows immediately that the resulting
 Blaschke product \(B\)  on the unit disk is prime, where
 $B(z)=P_{\mathbf s^*}(\eta(z))$.

\begin{thm}\cite[Theorem 4.5]{GLMR94}
Let \(B\) be a thin Blaschke product
 whose zeros accumulate at only one point
  of \( \mathbb{T}\), say at \(1\). Then \(B\) is prime.
\end{thm}

\noindent \textbf{3. Failure of essential injectivity}

The Cayley transform satisfies
\[
\eta(z)=\frac{1+z}{1-z},
\qquad
\frac{d\eta}{dz}=\frac{2}{(1-z)^2}\ne0.
\]
Hence critical points of the disk function
\(B\) correspond exactly to critical points
 of the associated half-plane function \(P_{\mathbf s^*}\):
\[
B'(z)=0
\quad\Longleftrightarrow\quad
P_{\mathbf s^*}'(\eta(z))=0.
\]
By (\ref{A25}), for every \(k\) there exist distinct critical points \(z_k,w_k\) such that
\[
B'(z_k)=B'(w_k)=0,
\qquad
B(z_k)=B(w_k).
\]
The pairs arising from distinct blocks are disjoint.
 Thus by definition,  \(B\) is not essentially
 injective on its critical set.

 \begin{rem} Remind that $I=[\frac{3}{4}, \frac{5}{4}]$. By the construction, for a fixed increasing
 sequence \(\{\lambda_k\}\) with \(\lambda_1>2\),
 there exist a carefully chosen infinite sequence
 \(\{R_k\}\) (as in (\ref{22a})--(\ref{22c})) and at least one \(\mathbf s^*\in I^\infty\)
  such that the corresponding product \(P_{\mathbf s^*}\)
 on $\mathbb{H}$ defines a prime thin Blaschke product
 that is not essentially injective on its critical set.
 \end{rem}

 \noindent \textbf{Acknowledgements.}
 The initial construction of the unique example
 in Section~\ref{6cexample} was suggested
 by ChatGPT~5.6 Plus. The authors identified
  a gap in the proposed construction,
  which was subsequently repaired through
   further dialogue with the AI.
   The authors have verified all details.
   We note that this example
   occupies only Section \ref{6cexample}
    of the paper; 
  the remaining sections are the authors' own work.

%
%
%
%

%



\vskip3mm \noindent{\noindent{School
of Mathematical Sciences, Fudan University, Shanghai, 200433, China;
E-mail:
danniguomath@163.com

\noindent School of Mathematics, East China University of
Science and Technology, Shanghai, 200237, China;
 E-mail:
hshuang@ecust.edu.cn
\

 \noindent Department of Mathematics, Jiangsu University of
Technology, Changzhou, 273000, China;
 E-mail:
 lishan\_math@163.com
 \

 \noindent College of Mathematics and Statistics,
 Chongqing University, Chongqing, 401331, China; E-mail: xianfengzhao@cqu.edu.cn
}
\end{document}